\documentclass[a4paper,11pt]{amsart}
\usepackage{amssymb}
\usepackage{latexsym}
\usepackage{amsmath}
\usepackage{enumerate}
\usepackage{multicol}

\usepackage{amsmath, hyperref}
\newtheorem{theorem}{Theorem}[section]
\newtheorem{lemma}[theorem]{Lemma}
\newtheorem{corollary}[theorem]{Corollary}
\newtheorem{proposition}[theorem]{Proposition}

\theoremstyle{definition}
\newtheorem{definition}[theorem]{Definition}
\newtheorem{question}[theorem]{Question}
\newtheorem{fact}[theorem]{Fact}

\newcommand{\PA}{\mathbf{PA}}

\numberwithin{equation}{section}

\title{Current research on G\"{o}del's incompleteness theorems}
\author{Yong Cheng}
\address{School of Philosophy,
Wuhan University, Wuhan, Hubei Province, P.R.China, 430072}
\email{world-cyr@hotmail.com}
\thanks{Some materials of my old paper ``Note on some misinterpretations of G\"{o}del's incompleteness theorems" have been incorporated into this paper. I would like to thank Matthias Baaz, Ulrich Kohlenbach, Taishi Kurahashi, Zachiri McKenzie, Fedor Pakhomov, Michael Rathjen, Saeed Salehi,  Sam Sanders   and Albert Visser  for  their valuable comments  on this work. Especially, I would like to thank Sam Sanders  and Zachiri McKenzie for their proofreading of my English writing of this paper. I would like to thank the referees for providing detailed and  helpful comments for improvements. I thank the lamp for my feet, and the light on my path. This paper is the research result of the Humanities and Social Sciences of Ministry of Education Planning Fund project ``Research on G\"{o}del's incompleteness theorem" (project no: 17YJA72040001). I would like to thank the fund support by the Humanities and Social Sciences of Ministry of Education Planning Fund.}



\begin{document}

\begin{abstract}
We give a survey of current research on  G\"{o}del's incompleteness theorems from the following three aspects: classifications of different proofs of G\"{o}del's incompleteness theorems, the limit of the applicability of G\"{o}del's first incompleteness theorem, and the limit of the applicability of G\"{o}del's second incompleteness theorem.
\end{abstract}

\setcounter{tocdepth}{4}
\setcounter{page}{0}
\tableofcontents
\thispagestyle{empty}
\newpage

\maketitle

\section{Introduction}

G\"{o}del's first and second incompleteness theorem are some of
the most important and profound results in the foundations of mathematics and have had wide influence on the development of logic, philosophy, mathematics, computer science as well as other fields. Intuitively speaking, G\"{o}del's incompleteness
theorems express that any rich enough logical system cannot prove its own \emph{consistency}, i.e.~ that
no contradiction like $0=1$ can be derived within this system.

\smallskip

G\"{o}del \cite{Godel 1931 original proof} proves his first incompleteness theorem  $(\sf G1)$  for a certain formal
system $\mathbf{P}$ related to Russell-Whitehead's \emph{Principia Mathematica}  based on the
simple theory of types over the natural number series and the Dedekind-Peano
axioms (see \cite{Beklemishev 45}, p.3).
%SAM
%Let $T$ be the theory formulated in the language of $\mathbf{P}$ and obtained by adding a primitive recursive set of axioms to the system $\mathbf{P}$. G\"{o}del's first incompleteness theorem $(\sf G1)$ says that if $T$ is $\omega$-consistent, %then $T$ is incomplete  (see \cite[Theorem VI]{Godel 1931 original proof}).
G\"{o}del announces the second incompleteness theorem $(\sf G2)$ in an abstract published in
October 1930: no consistency proof of systems such as Principia, Zermelo-Fraenkel
set theory, or the systems investigated by Ackermann and von Neumann is possible
by methods which can be formulated in these systems (see \cite{Richard Zach}, p.431).

\smallskip

 G\"{o}del  comments in a footnote of \cite{Godel 1931 original proof} that $\sf G2$ is
corollary  of $\sf G1$ (and in fact a formalized version of $\sf G1$): if $T$ is consistent, then the consistency of $T$ is not provable in $T$ where the consistency of $T$ is formulated as the arithmetic formula which says that there exists an unprovable sentence in $T$.
G\"{o}del \cite{Godel 1931 original proof} sketches a proof of  $\sf G2$   and promises to provide full details in a subsequent publication. This promise is not fulfilled, and a detailed
proof of $\sf G2$  for first-order arithmetic only appears in a monograph by Hilbert and Bernays \cite{Hilbert-Bernays}.
Abstract logic-free formulations of G\"{o}del's incompleteness theorems
have been given by Kleene \cite{Kleene 50} (``symmetric form"), Smullyan \cite{Smullyan 94} (``representation systems"), and others. The following  is a modern reformulation of G\"{o}del's incompleteness theorems.

\begin{theorem}[G\"{o}del, \cite{Godel 1931 original proof}]
Let $T$ be a recursively axiomatized extension of $\mathbf{PA}$.
\begin{enumerate}[(1)]
  \item [$\sf G1$] If $T$ is $\omega$-consistent, then $T$ is incomplete.
  \item [$\sf G2$] If $T$ is consistent, then the consistency of $T$ is not provable in $T$.
\end{enumerate}
\end{theorem}

G\"{o}del's incompleteness theorems $\sf G1$ and $\sf G2$  are of a rather different nature and scope.   In this paper,  we will discuss different versions of $\sf G1$ and $\sf G2$, from  incompleteness for extensions of $\mathbf{PA}$ to incompleteness for systems weaker than $\mathbf{PA}$ w.r.t.~ interpretation. We will freely use  $\sf G1$ and $\sf G2$ to refer to both G\"{o}del's first and second incompleteness theorems, and their different versions. The meaning of  $\sf G1$ and $\sf G2$ will be clear from the context in which we refer to them.

\smallskip

G\"{o}del's incompleteness theorems exhibit certain weaknesses and limitations of a given formal system.
For G\"{o}del, his incompleteness theorems indicate the creative power
of human reason. In  Emil Post's celebrated words: mathematical proof is an essentially
creative activity (see \cite{metamathematics}, p.339). The impact of G\"{o}del's incompleteness theorems is not confined to the community of mathematicians and logicians;
popular accounts are well-known within the general scientific
community and beyond.
G\"{o}del's incompleteness theorems raise a number of
philosophical questions concerning the
nature of logic and mathematics as well as mind and machine. For the impact of G\"{o}del's incompleteness
theorems, Feferman said:
\begin{quote}
their relevance to mathematical
logic (and its offspring in the theory of computation) is paramount; further, their philosophical relevance is significant, but in just what way is far from
settled; and finally, their mathematical relevance
outside of logic is very much unsubstantiated but
is the object of ongoing, tantalizing efforts (see \cite{Feferman 2006}, p.434).
\end{quote}
From the literature, there are some good textbooks and survey papers  on G\"{o}del's incompleteness theorems. For textbooks, we refer to \cite{Enderton 2001, metamathematics, Aspects of incompleteness, Franzen 04, Smith 2007, Boolos 93, Smullyan 92, Smullyan 94,  Pudlak 93, Friedman 18}. For survey papers, we refer to \cite{Smorynski 1977, Beklemishev 45, Henryk Kotlarski 04, Bernd Buldt 14, Visser 16, Blanck 17, Cheng on Godel 19}. In the last twenty years, there have been a lot of advances in  the study of incompleteness. We felt that a comprehensive survey paper for the current state-of-art of this research field is missing from the literature. The motivation of this paper is four-fold:
\begin{itemize}
  \item  Give the reader an overview of the current state-of-art of research on incompleteness.
  \item Classify  these new advances on incompleteness under some important themes.
  \item Propose some new questions not covered in the literature.
  \item Set the direction for the future research of incompleteness.
\end{itemize}
Due to  space limitations and our personal taste, it is impossible to cover all research results from the literature related to incompleteness in this survey. Therefore, we will focus on three aspects of new advances in research on incompleteness:
\begin{itemize}
  \item  classifications of different proofs of G\"{o}del's incompleteness theorems;
  \item the limit of the applicability of $\sf G1$;
  \item the limit of the applicability of $\sf G2$.
\end{itemize}
We think these  are the most important three aspects of research on incompleteness and reflect the depth and breadth of the  research on incompleteness after G\"{o}del.  In this survey, we will focus on  logical and  mathematical aspects of research on incompleteness.

\smallskip

An important and interesting topic concerning incompleteness is missing in this paper: philosophy of G\"{o}del's incompleteness theorems. For us, the widely discussed   and most important philosophical questions about G\"{o}del's incompleteness theorems are: the relationship between $\sf G1$ and the mechanism thesis, the status of G\"{o}del's disjunctive thesis, and the intensionality problem of $\sf G2$. We leave a survey of philosophical  discussions of G\"{o}del's incompleteness theorems for a future philosophy paper.

\smallskip

This paper is structured as follows. In Section 1, we introduce the motivation, the main content and  the structure of this  paper. In Section 2, we list the preliminary notions and definitions used in this paper. In Section 3, we examine different proofs of G\"{o}del's incompleteness theorems and classify these proofs based on nine  criteria.   In Section 4, we examine the limit of the applicability of  $\sf G1$ both for  extensions of $\mathbf{PA}$, and for theories weaker than $\mathbf{PA}$ w.r.t.~ interpretation. In Section 5, we examine the limit of the applicability of  $\sf G2$, and discuss  sources of indeterminacy in the formulation of the consistency statement.

\section{Preliminaries}\label{Preliminaries}
\subsection{Definitions and notations}\label{defi}
We list the definitions and notations required below.
These are standard and used throughout the literature.\smallskip

\begin{definition}[Basic notions]~
\begin{itemize}
\item A \emph{language} consists of an arbitrary number of relation and function symbols of arbitrary finite arity.\footnote{We may view nullary functions as constants, and nullary relations as propositional
variables.}  For a given theory $T$, we use $L(T)$ to denote  the language of $T$, and often equate $L(T)$ with the list of non-logical symbols  of the language.
\item For a formula $\phi$ in $L(T)$,   `$T\vdash\phi$' denotes that $\phi$ is provable in $T$: i.e., there is a finite sequence of formulas $\langle\phi_0, \cdots,\phi_n\rangle$ such that $\phi_n=\phi$, and for any $0\leq i\leq n$, either $\phi_i$ is an axiom of $T$, or $\phi_i$ follows from some $\phi_j\, (j<i)$ by using one inference rule.
\item  A theory $T$ is \emph{consistent} if no contradiction is provable in $T$.
\item We say a sentence $\phi$ is \emph{independent} of $T$  if $T\nvdash \phi$ and $T\nvdash \neg\phi$.
\item A theory $T$ is \emph{incomplete} if there is a sentence $\phi$ in $L(T)$ which is independent of $T$; otherwise, $T$ is
\emph{complete} (i.e., for any sentence $\phi$ in $L(T)$, either $T\vdash\phi$ or $T\vdash \neg\phi$).
\end{itemize}
\end{definition}

In this paper, we focus on first-order theories based on a countable language, and  always assume the \emph{arithmetization} of the base theory with a recursive set of non-logical symbols. For the technical details of arithmetization, we refer to \cite{metamathematics, Buss 1998}.
Arithmetization means that any formula or finite sequence of formulas can
be coded by a natural number, called the \emph{G\"{o}del number}.  This representation of syntax was pioneered by G\"odel.
\smallskip

\begin{definition}[Basic notions following arithmetization]~
\begin{itemize}
\item We say a set of sentences $\Sigma$  is \emph{recursive} if the set of G\"{o}del numbers of sentences in $\Sigma$ is recursive.
\item A theory $T$ is \emph{decidable} if the set of sentences provable in $T$ is recursive; otherwise it is
\emph{undecidable}.
\item A theory $T$ is \emph{recursively axiomatizable} if it has a recursive set of axioms (i.e.~ the set of G\"{o}del numbers of axioms of $T$ is recursive).
\item A theory $T$ is \emph{finitely axiomatizable} if it has a finite set of axioms.
 \item A theory $T$ is \emph{locally finitely satisfiable} if every finitely axiomatized subtheory of $T$ has a finite model.
\item A theory $T$ is \emph{recursively enumerable} (r.e.) if it has a recursively enumerable set of axioms.
\item A theory $T$ is \emph{essentially undecidable} if any recursively axiomatizable consistent extension of $T$ in the same language is undecidable.
\item A theory $T$ is \emph{essentially incomplete} if  any recursively axiomatizable consistent extension of $T$ in the same language is incomplete.\footnote{The theory of completeness/incompleteness  is closely related to the theory of decidability/undecidability (see \cite{undecidable}).}
\item A theory $T$ is
\emph{minimal essentially undecidable} if $T$ is essentially undecidable, and if deleting any axiom of $T$, the remaining theory is no longer essentially undecidable.
\end{itemize}
\end{definition}

\smallskip
\begin{definition}[Basic notations]~\label{}
\begin{itemize}
  \item We denote by   $\overline{n}$ the numeral representing $n\in \omega$ in $L(\mathbf{PA})$.
  \item We denote by  $\ulcorner\phi\urcorner$ the numeral
representing the G\"{o}del number of $\phi$.
\item We denote by $\ulcorner\phi(\dot{x})\urcorner$ the numeral representing the G\"{o}del number of the sentence obtained by replacing $x$ with the
value of $x$.\footnote{Note that  the variable $x$ is free in the formula $\ulcorner\phi(\dot{x})\urcorner$ but not in $\ulcorner\phi(x)\urcorner$.}
\end{itemize}
\end{definition}
\smallskip
\begin{definition}[Representations, translations, and interpretations]~
\begin{itemize}
\item A $n$-ary relation $R(x_1, \cdots, x_n)$ on $\omega^n$ is
\emph{representable} in $T$ if there is a formula $\phi(x_1, \cdots, x_n)$ such that $T\vdash \phi(\overline{m_1}, \cdots, \overline{m_n})$ when $R(m_1, \cdots, m_n)$ holds, and $T\vdash \neg\phi(\overline{m_1}, \cdots, \overline{m_n})$ when $R(m_1, \cdots, m_n)$ does not hold.
%We can generalize representation of relations to representation of disjoint pairs of relations: such a disjoint pair $\langle P^{+}, P^{-}\rangle$ where $P^{+}, P^{-}\subseteq \mathbb{N}^n$, is representable by  a formula $\phi(x_1, \cdots, x_n)$  in $T$ if  $T\vdash\phi(\overline{a_1}, \cdots, \overline{a_{n}})$
%for $\langle a_1, \cdots, a_{n}\rangle \in P^{+}$, and
%$T\vdash\neg\phi(\overline{a_1}, \cdots, \overline{a_{n}})$ for $\langle a_1, \cdots, a_{n}\rangle \in P^{-}$.
%A theory $T$ is said to be $\omega$-consistent if there is no formula $\varphi(x)$ such that $T \vdash\exists x \varphi(x)$ and for any $n\in\mathbb{N}$, $T \vdash\neg\varphi(\overline{n})$; and T is $1$-consistent if there is no such a $\Delta^0_1$ formula $\varphi(x)$.
\item We say that a total function $f(x_1, \cdots, x_n)$ on $\omega^n$ is \emph{representable in $T$} if there is a formula $\varphi(x_1, \cdots, x_n,y)$ such  that $T\vdash \forall y(\varphi(\overline{a_1}, \cdots, \overline{a_{n}},y)\leftrightarrow y=\overline{m})$ whenever $a_1, \cdots, a_n, m\in \omega$ are such that $f(a_1, \cdots, a_n)=m$.
%$I\Sigma_{n}$ is a fragment of $\mathbf{PA}$ obtained by restricting the axiom scheme of induction to $\Sigma_{n}$ formulas.
\item Let $T$ be a theory in a language $L(T)$, and $S$ a theory in a language $L(S)$. In its simplest form,
a \emph{translation} $I$ of language $L(T)$ into language $L(S)$ is specified by the following:
\begin{itemize}
  \item an $L(S)$-formula $\delta_I(x)$ denoting the domain of $I$;
  \item for each relation symbol $R$ of $L(T)$,  as well as the equality relation =, an $L(S)$-formula $R_I$ of the same arity;
  \item for each function symbol $F$ of $L(T)$ of arity $k$, an $L(S)$-formula $F_I$ of arity $k + 1$.
\end{itemize}
\item If $\phi$ is an $L(T)$-formula, its $I$-translation $\phi^I$ is an $L(S)$-formula constructed as follows: we rewrite the
formula in an equivalent way so that function symbols only occur in atomic subformulas of the
form $F(\overline{x}) = y$, where $x_i, y$ are variables; then we replace each such atomic formula with $F_I(\overline{x}, y)$,
we replace each atomic formula of the form $R(\overline{x})$ with $R_I(\overline{x})$, and we restrict all quantifiers and
free variables to objects satisfying $\delta_I$. We take care to rename bound variables to avoid variable
capture during the process.
\item A translation $I$ of $L(T)$ into $L(S)$ is an \emph{interpretation} of $T$ in $S$ if $S$ proves the following:
\begin{itemize}
  \item for each function symbol $F$ of $L(T)$ of arity $k$, the formula expressing that $F_I$ is total on $\delta_I$:
\[\forall x_0, \cdots \forall x_{k-1} (\delta_I(x_0) \wedge \cdots \wedge \delta_I(x_{k-1}) \rightarrow \exists y (\delta_I(y) \wedge F_I(x_0, \cdots, x_{k-1}, y)));\]
  \item the $I$-translations of all axioms of $T$, and axioms of equality.
\end{itemize}
\end{itemize}
\end{definition}
The simplified picture of translations and interpretations above actually describes only \emph{one-dimensional}, \emph{parameter-free}, and \emph{one-piece}  translations. For precise
definitions of a
\emph{multi-dimensional interpretation}, an
\emph{interpretation with parameters},  and a
\emph{piece-wise interpretation}, we refer to \cite{paper on Q} \cite{Visser 11}  \cite{Visser 14}   for more details.

\smallskip

The notion of interpretation provides us with a method for comparing different theories in different languages, as follows.
\smallskip

\begin{definition}[Interpretations II]~
\begin{itemize}
\item A theory $T$ is \emph{interpretable} in a theory $S$ if there exists an
interpretation  of $T$ in $S$. If $T$ is interpretable in $S$, then all sentences provable (refutable) in $T$ are mapped, by the interpretation function, to sentences provable (refutable) in $S$.
\item We say that a theory $U$ \emph{weakly interprets} a theory $V$ (or $V$ is
\emph{weakly interpretable} in $U$) if $V$ is interpretable in some consistent extension  of $U$ in the same language (or  equivalently, for some interpretation $\tau$, the
theory $U + V^{\tau}$ is consistent).
\item Given theories $S$ and $T$, let `$S\unlhd T$' denote that $S$ is interpretable in $T$ (or $T$ interprets $S$); let `$S\lhd T$' denote that  $T$ interprets $S$ but $S$ does not interpret  $T$; we say $S$ and $T$ are
\emph{mutually interpretable} if $S\unlhd T$ and $T\unlhd S$.
\end{itemize}
\end{definition}
%We say that $I$ is a
%\emph{faithful interpretation} of $T$ in $S$ iff $I$ is an interpretation of $T$ in $S$ such that for any sentence $\phi$ %in $L(T)$, $T\vdash \phi$ iff $S\vdash \phi^{I}$.
Interpretability provides us with one measure of comparing strength of different theories. If  theories $S$ and $T$ are mutually interpretable, then $T$ and $S$ are equally strong w.r.t.~ interpretation. In this paper, whenever we say that theory $S$ is weaker than theory $T$ w.r.t.~ interpretation, this means that $S\lhd T$.

\smallskip

%Let $S\unlhd T$ denote $S$ is interpretable in $T$; let $S\lhd T$ denote $S$ is interpretable in $T$ but $T$ is not interpretable in $S$. $S$ and $T$ are mutually interpretable if $S\unlhd T$ and $T\unlhd S$.
%Let $Int(S)$ denote the degree of interpretation of theory $S$. $Int(T)<Int(S)$ means that $T$ is interpretable in $S$ but $S$ is not interpretable in $T$. $Int(T)=Int(S)$ means that $T$ and $S$ are mutually interpretable.
A general method for establishing the undecidability of theories is developed in \cite{undecidable}.
The following theorem provides us with two methods for proving the essentially undecidability of a theory respectively via interpretation and representability.

\begin{theorem}[Theorem 7, Corollary 2, \cite{undecidable}]\label{interpretable theorem}~
\begin{itemize}
  \item   Let $T_1$ and $T_2$ be two consistent theories such that $T_2$ is interpretable in $T_1$. If $T_2$ is essentially undecidable, then $T_1$ is essentially undecidable.
  \item   If all recursive functions are representable in a consistent theory $T$, then $T$ is essentially undecidable.
\end{itemize}
\end{theorem}
We shall also need some basic notions from recursion theory, as follows.
\smallskip

\begin{definition}[Basic recursion theory]~
\begin{itemize}
\item Let $\phi_0, \phi_1, \cdots$ be a list of all unary computable (partial recursive) functions such that $\phi_i(j)$, if it exists, can be computed from $i$ and $j$.
\item A recursively enumerable set (r.e. for
short) is the domain of $\phi_i$ for some $i \in\omega$, which is denoted by $W_i$.
\item The notation $\phi_i(j)\uparrow$ means that the
function $\phi_i$ is not defined at $j$, or $j \notin W_i$; and $\phi_i(j)\downarrow$ means that $\phi_i$ is defined at $j$, or $j \in W_i$.
\end{itemize}
\end{definition}

Provability logic provides us with an important tool to study the meta-mathematics of arithmetic and incompleteness. A good reference on the basics of provability logic is \cite{Boolos 93}.\smallskip

\begin{definition}[Modal logic]~
\begin{itemize}
  \item The modal system $\mathbf{K}$ consisting of the following axiom schemes:
\begin{itemize}
  \item All tautologies;
  \item $\Box(A\rightarrow B)\rightarrow (\Box A\rightarrow\Box B)$;
\end{itemize}
as well as two inference rules:
\begin{itemize}
  \item if $\vdash A$ and $\vdash A\rightarrow B$, then $\vdash B$;
  \item if $\vdash A$, then $\vdash  \Box A$.
\end{itemize}
  \item We denote by  $\mathbf{GL}$ the modal system consisting of all axioms of $\mathbf{K}$, all instances of the scheme $\Box(\Box A\rightarrow A)\rightarrow\Box A$, and the same inference rules with $\mathbf{K}$.
  \item
We denote by  $\mathbf{GLS}$ the modal system consisting of all theorems of $\mathbf{GL}$, and all instances of the scheme $\Box A\rightarrow A$. However, $\mathbf{GLS}$ has only one inference rule: Modus Ponens.
\end{itemize}
\end{definition}

\subsection{Logical systems}\label{logical system}
In this section, we introduce some well-known theories weaker than $\mathbf{PA}$ w.r.t.~ interpretation from the literature. In Section \ref{Generalizations}, we will show that these theories are essentially incomplete.
\smallskip

Robinson Arithmetic $\mathbf{Q}$ is introduced in \cite{undecidable} by Tarski, Mostowski and
Robinson as a base axiomatic theory for investigating incompleteness, and undecidability.\smallskip

\begin{definition}[Robinson Arithmetic $\mathbf{Q}$]~\label{def of Q}
Robinson Arithmetic $\mathbf{Q}$  is  defined in   the language $\{\mathbf{0}, \mathbf{S}, +, \times\}$ with the following axioms:
\begin{description}
  \item[$\mathbf{Q}_1$] $\forall x \forall y(\mathbf{S}x=\mathbf{S} y\rightarrow x=y)$;
  \item[$\mathbf{Q}_2$] $\forall x(\mathbf{S} x\neq \mathbf{0})$;
  \item[$\mathbf{Q}_3$] $\forall x(x\neq \mathbf{0}\rightarrow \exists y (x=\mathbf{S} y))$;
  \item[$\mathbf{Q}_4$]  $\forall x\forall y(x+ \mathbf{0}=x)$;
  \item[$\mathbf{Q}_5$] $\forall x\forall y(x+ \mathbf{S} y=\mathbf{S} (x+y))$;
  \item[$\mathbf{Q}_6$] $\forall x(x\times \mathbf{0}=\mathbf{0})$;
  \item[$\mathbf{Q}_7$] $\forall x\forall y(x\times \mathbf{S} y=x\times y +x)$.
\end{description}
\end{definition}

Robinson Arithmetic $\mathbf{Q}$ is  very weak: it cannot even prove that addition is associative.
\smallskip

\begin{definition}[Peano Arithmetic  $\mathbf{PA}$]~\label{}
The theory $\mathbf{PA}$ consists of the axioms $\mathbf{Q}_1$-$\mathbf{Q}_2$, $\mathbf{Q}_4$-$\mathbf{Q}_7$ in Definition \ref{def of Q} and the following axiom scheme of induction:
\begin{equation}\tag{\textsf{Induction}}
(\phi(\mathbf{0})\wedge \forall x(\phi(x)\rightarrow \phi(\mathbf{S} x)))\rightarrow \forall x  \phi(x),
\end{equation}
where $\phi$ is a formula with at least one free variable $x$.
Let $\mathfrak{N}=\langle\mathbb{N}, +, \times\rangle$ denote the standard model of arithmetic.
\end{definition}

We now introduce a well-known hierarchy of $L(\mathbf{PA})$-formulas called the \emph{arithmetical hierarchy} (see \cite{metamathematics, Pudlak 93}).\smallskip
\begin{definition}[Arithmetical hierarchy]~
\begin{itemize}
\item  \emph{Bounded formulas}   ($\Sigma^0_0$, or $\Pi^0_0$, or $\Delta^0_0$ formula) are built from atomic
 formulas using only propositional connectives and bounded quantifiers (in the form $\forall x\leq y$ or $\exists x\leq y$).
 \item A formula is  $\Sigma^0_{n+1}$ if it has the
form $\exists x\phi$ where $\phi$ is $\Pi^0_{n}$.
\item A formula is $\Pi^0_{n+1}$ if it has the form $\forall x\phi$ where $\phi$ is $\Sigma^0_{n}$. Thus, a $\Sigma^0_{n}$-formula has a block of $n$ alternating quantifiers, the first one
being existential, and this block is followed by a bounded formula. Similarly
for $\Pi^0_{n}$-formulas.
\item A formula is $\Delta^0_n$ if it is equivalent to both a $\Sigma^0_{n}$ formula and a $\Pi^0_{n}$ formula.
\end{itemize}
\end{definition}
We can now formally introduce the notion of consistency in its various guises, as well as the various fragments of Peano arithmetic $\mathbf{PA}$.\smallskip

\begin{definition}[Formal consistency and systems]~
\begin{itemize}
\item A theory $T$ is said to be  \emph{$\omega$-consistent}  if there is no formula $\varphi(x)$ such that $T \vdash\exists x \varphi(x)$, and for any $n\in\omega$, $T \vdash\neg\varphi(\bar{n})$.
\item A theory $T$ is \emph{$1$-consistent}  if there is no such  $\Delta^0_1$ formula $\varphi(x)$.
\item We say a theory $T$ is \emph{$\Sigma^0_1$-sound}  if for any $\Sigma^0_1$ sentences $\phi$, if $T\vdash \phi$, then $\mathfrak{N}\models\phi$.
 \item The collection axiom for $\Sigma^0_{n+1}$ formulas is the following principle: $(\forall x< u)(\exists y) \varphi(x,y)\rightarrow (\exists v)(\forall x< u)(\exists y< v) \varphi(x,y)$ where
$\varphi(x,y)$ is a $\Sigma^0_{n+1}$ formula possibly containing parameters distinct from $u,v$.
\item The theory $I\Sigma_n$ is $\mathbf{Q}$ plus induction for $\Sigma^0_n$ formulas, and $B\Sigma_{n+1}$ is $I\Sigma_0$ plus collection for $\Sigma^0_{n+1}$ formulas.
\item The theory $I\Delta_0$ is $\mathbf{Q}$ plus induction for $\Delta^0_0$ formulas.
\item The theory $\mathbf{PA}$ is the union of all $I\Sigma_{n}$.
\end{itemize}
\end{definition}

It is well-known that the following form a strictly increasing hierarchy:
\[
I\Sigma_0, B\Sigma_{1}, I\Sigma_1, B\Sigma_{2},\cdots, I\Sigma_n, B\Sigma_{n+1}, \cdots, \mathbf{PA}.
\]

Moreover, there are \emph{weak} fragments of $\mathbf{PA}$ that play an important role in computer science, namely in \emph{complexity theory} (\cite{Buss 86, Buss 1998}).
These systems are based on the following concept.\smallskip

By \cite[Proposition 2, p.299]{Interpretability in Robinson's Q}, there is a bounded formula ${\sf Exp(x, y, z)}$ such that $I\Sigma_0$ proves that  ${\sf Exp(x, 0, z)} \leftrightarrow z= 1$, and ${\sf Exp(x, Sy, z)} \leftrightarrow \exists t ({\sf Exp(x, y, t)} \wedge z=t\cdot x)$. However,  $I\Sigma_0$ cannot prove  the totality of ${\sf Exp(x, y, z)}$.\smallskip

\begin{definition}[Sub-exponential functions]~
\begin{itemize}
\item Let $\mathbf{exp}$ denote the statement postulating the totality of the exponential function $\forall x\forall y\exists z  {\sf Exp(x, y, z)}$.
\item Elementary Arithmetic ($\mathbf{EA}$) is  $I\Delta_0 +\mathbf{exp}$.
\item Define $\omega_1(x)=x^{\mid x\mid}$, and
$\omega_{n+1}(x)=2^{\omega_{n}(\mid x\mid)}$ where $|x|$ is the length of the binary expression of $x$.
\item Let $\Omega_{n}\equiv (\forall x)(\exists y) (\omega_{n}(x) =y)$ express that $\omega_{n}(x)$ is total.
\end{itemize}
\end{definition}

\begin{theorem}[\cite{Petr 93, Interpretability in Robinson's Q}]\label{interpretation theorem}~
\begin{itemize}
  \item  The theory $I\Sigma_0+\Omega_n$ is interpretable in $\mathbf{Q}$ for any $n\geq 1$ (see \cite[Theorem 3, p.304]{Interpretability in Robinson's Q}).
  \item The theory $I\Sigma_0+ \mathbf{exp}$ is not interpretable in $\mathbf{Q}$.\footnote{See \cite[Theorem 6, p.313]{Interpretability in Robinson's Q}. Solovay proves that $I\Sigma_0 + \neg \mathbf{exp}$ is interpretable in $\mathbf{Q}$ (see \cite[Theorem 7, p.314]{Interpretability in Robinson's Q}).}
  \item  The theory $I\Sigma_1$ is not interpretable in $I\Sigma_0+ \mathbf{exp}$ (see \cite[Theorem 1.1]{Petr 93}, p.186).
      \item The theory $I\Sigma_{n+1}$ is not interpretable in $B\Sigma_{n+1}$ (see \cite[Theorem 1.2]{Petr 93}, p.186).
      \item  The theory $B\Sigma_{1}+ exp$ is interpretable in $I\Sigma_{0}+ \mathbf{exp}$ (see \cite[Theorem 2.4]{Petr 93}, p.188).
      \item The theory $B\Sigma_{1}+\Omega_{n}$ is interpretable in $I\Sigma_{0}+\Omega_{n}$  for each $n\geq 1$  (see \cite[Theorem 2.5]{Petr 93}, p.189).
      \item   The theory $B\Sigma_{n+1}$ is interpretable in $I\Sigma_{n}$ for each $n\geq 0$ (see \cite[Theorem 2.6]{Petr 93}, p.189).
\end{itemize}
\end{theorem}

The theory $\mathbf{PA}^{-}$ is the theory of commutative, discretely
ordered semi-rings with a minimal element plus the subtraction axiom.
The theory $\mathbf{PA}^{-}$ has the following axioms, where the language $L(\mathbf{PA}^{-})$ is $L(\mathbf{PA})\cup\{\leq\}$:\smallskip

\begin{definition}[The system $\mathbf{PA}^{-}$]~
\vspace{-2mm}
\begin{multicols}{2}
\begin{itemize}
  \item $x + 0 = x$;
  \item $x + y = y + x$;
  \item $(x + y) + z = x + (y + z)$;
  \item $x \times 1 = x$;
  \item $x \times y = y \times x$;
  \item $(x \times y) \times z = x \times (y \times z)$;
  \item $x \times (y + z) = x \times y + x \times z$;
  \item $x \leq y \vee y \leq x$;
  \item $(x \leq y \wedge y \leq z)\rightarrow x \leq z$;
  \item $x + 1 \nleq x$;
  \item $x \leq y \rightarrow (x = y \vee x + 1 \leq y)$;
  \item $x \leq y \rightarrow x + z \leq y + z$;
  \item $x \leq y \rightarrow x \times z \leq y \times z$;
  \item $x \leq y \rightarrow \exists z (x + z = y)$.
  \end{itemize}
\end{multicols}
\end{definition}
The theory $\mathbf{Q}^+$ is the extension of $\mathbf{Q}$ in the language $L(\mathbf{Q}^+)=L(\mathbf{Q})\cup\{\leq\}$ with the following extra axioms:\smallskip

\begin{definition}[The system $\mathbf{Q}^{+}$]~
The system $\mathbf{Q}^{+}$ is $\mathbf{Q}$ plus
\begin{description}
  \item[$\mathbf{Q}_8$] $(x + y) + z = x + (y + z)$;
  \item[$\mathbf{Q}_9$] $x \times (y + z) = x \times y + x \times z$;
  \item[$\mathbf{Q}_{10}$] $(x \times y) \times z = x\times (y\times z)$;
  \item[$\mathbf{Q}_{11}$] $x + y = y + x$;
  \item[$\mathbf{Q}_{12}$] $x\times y = y \times x$;
  \item[$\mathbf{Q}_{13}$] $x\leq y\leftrightarrow \exists z (x + z = y)$.
\end{description}
\end{definition}

Andrzej Grzegorczyk considers a theory $\mathbf{Q}^{-}$ in which addition and multiplication
 satisfy natural reformulations of the axioms of $\mathbf{Q}$ but are possibly \emph{non-total} functions. More exactly, the language of $\mathbf{Q}^{-}$ is $\{\mathbf{0}, \mathbf{S}, A, M\}$ where $A$ and $M$
are ternary relations.\smallskip

\begin{definition}[The system $\mathbf{Q}^{-}$]
The axioms of $\mathbf{Q}^{-}$ are the axioms $\mathbf{Q}_1$-$\mathbf{Q}_3$ of $\mathbf{Q}$ plus
the following six axioms about $A$ and $M$:
\begin{description}
  \item[A] $\forall x\forall y\forall z_1\forall z_2(A(x, y, z_1)\wedge A(x, y, z_2) \rightarrow z_1 = z_2)$;
  \item[M] $\forall x\forall y\forall z_1\forall z_2(M(x, y, z_1) \wedge M(x, y, z_2) \rightarrow z_1 = z_2)$;
  \item[G4] $\forall x \, A(x, 0, x)$;
  \item[G5] $\forall x\forall y\forall z(\exists u(A(x, y, u) \wedge z = S(u)) \rightarrow A(x, S(y), z))$;
  \item[G6] $\forall x\, M(x, 0, 0)$;
  \item[G7] $\forall x\forall y\forall z(\exists u(M(x, y, u) \wedge A(u, x, z)) \rightarrow M(x, S(y), z))$.
\end{description}
\end{definition}
Samuel R. Buss   \cite{Buss 86} introduces $\mathbf{S^1_2}$, a finitely axiomatizable  theory, to study polynomial time computability. The theory $\mathbf{S^1_2}$  provides what is needed for formalizing the proof of $\sf G2$ in a natural and
effortless way:
this process is actually  easier in Buss' theory than in full
$\mathbf{PA}$, since the
restrictions present in $\mathbf{S^1_2}$ prevent one from making wrong turns and inefficient choices (see \cite{paper on Q}).
\smallskip

Next, we introduce \emph{adjunctive set theory} $\mathbf{AS}$ which has a language with only one binary relation
symbol `$\in$'.
\smallskip

\begin{definition}[Adjunctive set theory $\mathbf{AS}$, \cite{Nelson 86}]~\label{}
The axioms of $\mathbf{AS}$ consist of the following:
\begin{description}
  \item[AS1] $\exists x \forall y (y \notin x)$.
  \item[AS2] $\forall x\forall y \exists z \forall u (u \in z \leftrightarrow (u = x \vee u = y))$.
\end{description}
\end{definition}

We now consider the theory $\mathbf{R}$ introduced by A. Tarski, A. Mostowski and R. Robinson in \cite{undecidable}, and some variants of it.\smallskip

\begin{definition}[The theory $\mathbf{R}$]~
Let $\mathbf{R}$ be the theory consisting of schemes $\sf {Ax1}$-$\sf {Ax5}$ with $L(\mathbf{R})=\{\overline{0}, \cdots, \overline{n}, \cdots, +, \times, \leq\}$ where $m, n \in \omega$.
\begin{description}
  \item[Ax1] $\overline{m}+\overline{n}=\overline{m+n}$;
  \item[Ax2] $\overline{m}\times\overline{n}=\overline{m\times n}$;
  \item[Ax3] $\overline{m}\neq\overline{n}$ if $m\neq n$;
  \item[Ax4] $\forall x(x\leq \overline{n}\rightarrow x=\overline{0}\vee \cdots \vee x=\overline{n})$;
  \item[Ax5] $\forall x(x\leq \overline{n}\vee \overline{n}\leq x)$.
\end{description}
\end{definition}
%Let $\mathbf{R}$ be the system consisting of schemes $Ax1$-$Ax5$.
As it happens, the system $\mathbf{R}$ contains all key properties of arithmetic for the proof of $\sf G1$.
Unlike $\mathbf{Q}$, the theory $\mathbf{R}$ is not finitely axiomatizable.\smallskip

\begin{definition}[Variations of $\mathbf{R}$]~
\begin{itemize}
  \item Let $\mathbf{R}_0$ be $\mathbf{R}$ without $\mathbf{Ax5}$.
  \item
Let $\mathbf{R}_1$ be the system consisting of schemes $\mathbf{Ax1}, \mathbf{Ax2}, \mathbf{Ax3}$ and $\mathbf{Ax4^{\prime}}$ where the latter is as follows
\begin{description}
\item[Ax4$^{\prime}$] $\forall x(x\leq \overline{n}\leftrightarrow x=\overline{0}\vee \cdots \vee x=\overline{n}).$
\end{description}
  \item Let $\mathbf{R}_2$ be the system consisting of schemes $\mathbf{Ax2}, \mathbf{Ax3}$ and $\mathbf{Ax4^{\prime}}$.
\end{itemize}
\end{definition}

The `concatination' theory $\mathbf{TC}$ has the language $\{\frown, \alpha, \beta\}$ with a binary function
symbol and two constants.\smallskip
\begin{definition}[The system $\mathbf{TC}$]~
\begin{description}
  \item[TC1] $\forall x\forall y\forall z(x\frown (y\frown z) = (x\frown y)\frown z)$;
  \item[TC2] $\forall x\forall y\forall u\forall v(x\frown y = u\frown v \rightarrow ((x = u \wedge y = v) \vee\exists w((u = x\frown w \wedge w\frown v= y) \vee (x = u\frown w \wedge w\frown y = v))))$;
  \item[TC3] $\forall x\forall y(\alpha \neq x\frown y)$;
  \item[TC4] $\forall x\forall y(\beta \neq x\frown y)$;
  \item[TC5] $\alpha \neq\beta$.
\end{description}
\end{definition}

Primitive recursive arithmetic ($\mathbf{PRA}$) is a quantifier-free formalization of the natural numbers, and the language of $\mathbf{PRA}$ can express arithmetic statements involving natural numbers and any primitive recursive function. Weak Konig's Lemma  ($\mathbf{WKL_0}$) states that every infinite binary tree has an infinite branch. 
We refer to \cite{Pudlak 93, SOA} for the definitions of $\mathbf{PRA}$ and $\mathbf{WKL_0}$.
In a nutshell, the former system allows us to perform `iteration of functions $f:\mathbb{N} \rightarrow \mathbb{N}$', while the latter expresses a basic compactness argument for Cantor space.

\begin{theorem}[Friedman's conservation theorem, Theorem 2.1, \cite{Kikuchi-Tanaka 94}]\label{conservation thm}
If $\mathbf{WKL_0}\vdash\phi$, then $\mathbf{PRA}\vdash\phi$ for any $\Pi^0_2$ sentence $\phi$ in $L(\mathbf{PA})$.
\end{theorem}

Finally, diagnolisation, in one form or other, forms the basis for the proof of $\sf G2$.  The following lemma is crucial in this regard.\smallskip

\begin{lemma}[The Diagnolisation Lemma]~\label{diagonal lemma}
Let $T$ be a consistent r.e.~ extension of $\mathbf{Q}$. For any formula $\phi(x)$ with exactly one free variable, there exists a sentence $\theta$ such that $T\vdash\theta\leftrightarrow\phi(\ulcorner\theta\urcorner)$.
\end{lemma}
Lemma \ref{diagonal lemma} is the simplest and most often used version of the Diagnolisation Lemma. For a generalized version of the Diagnolisation Lemma, we refer to \cite{Boolos 93}. In this paper, we use the term ``Diagnolisation Lemma" to refer to Lemma \ref{diagonal lemma} and some variants of the generalized version.

\section{Proofs of G\"{o}del's incompleteness theorems}\label{diff proofs of incomplete}
\subsection{Introduction}
In this section, we discuss different proofs of G\"{o}del's incompleteness theorems from the literature, and propose nine criteria  for classifying them.

\smallskip

First of all, there are no requirements on the independent sentence in $\sf G1$.  In particular, such a sentence
need not have any mathematical meaning.  This is often the case when meta-mathematical
(proof-theoretic or recursion-theoretic or model-theoretic) methods are used to construct the
independent sentence.
In Section \ref{waf}-\ref{woef}, we will discuss proofs of  G\"{o}del's incompleteness theorems \emph{via pure logic}. In Section \ref{concrete incompleteness}, we will give an overview of the ``concrete incompleteness" research program
which seeks to identify  natural independent sentences \emph{with real mathematical meaning}.

\smallskip

Secondly, we say that a proof of $\sf G1$ is \emph{constructive} if it explicitly constructs the independent
sentence from the base theory  by algorithmic means. A non-constructive proof of $\sf G1$ only  proves the mere existence of the independent sentence
and does not show its existence algorithmically.
We say that a proof of $\sf G1$ for theory $T$ has \emph{the Rosser property} if the proof only assumes that $T$ is consistent instead of assuming that $T$ is  $\omega$-consistent or 1-consistent or $\Sigma^0_1$-sound; all these notions are introduced in Section \ref{logical system}.
%Even though the diagonalization lemma  has a constructive proof to find the fixed point sentence,  we always view the proof using the diagonalization lemma as non-constructive since the set of fixed points of a given formula with one free variable is not recursive (in fact creative, see \cite{Blanck 17}).

\smallskip

%In this paper, proof via model-theoretic method usually means the application of the Arithmetic Completeness Theorem since the model-theoretic proof of G\"{o}del's incompleteness theorems usually uses the Arithmetic Completeness Theorem.
After G\"{o}del, many different proofs of G\"{o}del's incompleteness theorems have been found. These proofs can be classified using  the following criteria:
\begin{itemize}
 \item  proof-theoretic proof;
  \item recursion-theoretic proof;
  \item model-theoretic proof;
  \item proof via arithmetization;
  \item proof via the Diagnolisation Lemma;
  \item proof based on ``logical paradox";
  \item constructive proof;
  \item  proof having the Rosser property;
  \item the independent sentence has natural and real mathematical content.\footnote{I.e.~ G\"{o}del's sentence is a pure logical construction (via the arithmetization of syntax and provability predicate) and has no relevance with classic mathematics (without any combinatorial or number-theoretic content). On the contrary, Paris-Harrington Principle is an independent arithmetic sentence from classic mathematics with combinatorial content.}
\end{itemize}
However, these aspects are not exclusive: a proof of $\sf G1$ or $\sf G2$ may satisfy several of the above criteria.
\smallskip
%\subsection{Incompleteness via pure logic}\label{proof via logic}

Thirdly, there are two kinds of proofs of G\"{o}del's incompleteness theorems via pure logic: one based on logical paradox and one not based on logical paradox. In Section \ref{waf}, we first provide an overview of the modern reformulation of proofs of  G\"{o}del's incompleteness  theorems.  We discuss  proofs of G\"{o}del's incompleteness theorems not based on logical paradox in Section~\ref{foef}. We discuss proofs of G\"{o}del's incompleteness theorems based on logical paradox in Section \ref{woef}.

\subsection{Overview and modern formulation}\label{waf}
In a nutshell, the three main ideas in the (modern/standard) proofs of $\sf G1$ and $\sf G2$ are \emph{arithmetization}, \emph{representability}, and \emph{self-reference}, as discussed in detail in Section \ref{defi2}.
Interesting properties of $\sf G1$ and $\sf G2$ are discussed in Sections \ref{waf1} and \ref{waf2}, while the formalized notions of `proof' and `truth' are discussed in Section~\ref{waf3}.
Finally, we formulate a blanket caveat for the rest of this section:
\begin{center}
\emph{Unless stated otherwise,  we will always assume that $T$ is a recursively axiomatizable consistent extension of $\mathbf{Q}$.}
\end{center}
Other sections shall contain similar caveats and we sometimes stress these.
\subsubsection{Three steps towards \textsf{\textup{G1}} and \textsf{\textup{G2}}}\label{defi2}
Intuitively speaking, G\"odel's incompleteness theorems can be proved based on the following key ingredients.
\begin{itemize}
\item \textbf{Arithmetization}: since \textsf{G1} and \textsf{G2} are theorems about properties of the syntax of logic, we need to somehow represent the latter, which is done via a coding scheme called \emph{arithmetization}.
\item \textbf{Representations}: the notion of `proof' and related concepts in \textsf{G1} and \textsf{G2} are then expressed (`represented') via arithmetization.
\item \textbf{Self-reference:}  given a representation of `proof' and related concepts, one can write down formal statements that intuitively express `self-referential' things like `this sentence does not have a proof'.
\end{itemize}
As we will see, the intuitively speaking `self-referential' statements are the key to proving \textsf{G1} and \textsf{G2}.
We now discuss these three notions in detail.

\smallskip

First of all, \textbf{arithmetization} has the following intuitive content: it establishes a one-to-one correspondence between expressions of $L(T)$ and natural numbers. Thus, we can translate metamathematical
statements about the formal theory $T$ into statements
about natural numbers. Furthermore, fundamental metamathematical relations can be translated in this way into
certain recursive relations, hence into relations representable in
 $T$. Consequently, one can speak about a formal system of arithmetic, and about its properties
as a theory in the system itself! This is the essence of G\"{o}del's
idea of arithmetization, which was revolutionary at a time when computer hardware and software did not exist yet.
%Under G\"{o}del's arithmetization, the set of G\"{o}del numbers of axioms of $T$ is recursive.

\smallskip

Secondly, in light of the previous, we can define certain relations on natural numbers that express or \textbf{represent} crucial metamathematical concepts related to the formal system $T$, like `proof' and `consistency'.
For example, modulo plenty of technical details, we can readily define a binary relation on $\omega^2$ expressing what it means to prove a formula in $T$, namely as follows:
\begin{center}
 $Proof_T(m,n)$ if and only if $n$ is the G\"{o}del number of a proof in $T$ of the formula with G\"{o}del number $m$.
 \end{center}
 %Then we need to show that  Prf_T(m,n) is recursive.
Moreover, we can show that the  relation $Proof_T(m,n)$  is
recursive. In addition, G\"{o}del proves that every recursive relation is representable
in $\PA$.

\smallskip

Next, let  $\mathbf{Proof}_T(x,y)$ be the formula which represents $Proof_T(m,n)$ in $\mathbf{PA}$.\footnote{Via  arithmetization  and representability, one can speak about the property of $T$  in $\mathbf{PA}$ itself!} From the formula $\mathbf{Proof}_T(x,y)$, we can define the `provability' predicate $\mathbf{Prov}_T(x)$ as $\exists y \mathbf{Proof}_T(x,y)$.
%Hence, if  $T$ is $\omega$-consistent, then G\"{o}del sentence is independent from $T$.
The provability predicate $\mathbf{Prov}_T(x)$ satisfies the following conditions which show that formal and intuitive provability have the same properties.
\begin{enumerate}[(1)]
  \item If $T \vdash\varphi$, then $T \vdash \mathbf{Prov}_T(\ulcorner\varphi\urcorner)$;
  \item $T \vdash \mathbf{Prov}_T(\ulcorner\varphi
      \rightarrow\psi\urcorner)\rightarrow ( \mathbf{Prov}_T(\ulcorner\varphi\urcorner)\rightarrow \mathbf{Prov}_T(\ulcorner\psi\urcorner))$;
  \item $T\vdash \mathbf{Prov}_T(\ulcorner\varphi\urcorner) \rightarrow \mathbf{Prov}_T(\ulcorner \mathbf{Prov}_T(\ulcorner\varphi\urcorner)\urcorner)$.
\end{enumerate}
For the proof of $\sf G1$, G\"{o}del defines  the \emph{G\"{o}del sentence} $\mathbf{G}$ which asserts its own unprovability in $T$ via a \textbf{self-reference} construction.
G\"{o}del shows  that  if $T$ is consistent, then $T\nvdash\mathbf{G}$, and if $T$ is $\omega$-consistent, then $T\nvdash\neg\mathbf{G}$.
One way of obtaining such a G\"odel sentence is the \emph{Diagnolisation Lemma} which intuitively speaking implies that the predicate $\neg\mathbf{Prov}_T(x)$ has a \emph{fixed point}, i.e.\ there is a sentence $\theta$ in $L(T)$ such that
\[
T \vdash \theta \leftrightarrow\neg \mathbf{Prov}_T(\ulcorner \theta \urcorner).
\]
Clearly, $T\nvdash\theta$ while $\theta$ intuitively expresses its own unprovability, i.e.\ the aforementioned \textbf{self-referential} nature.

\smallskip

For the proof of $\sf G2$, we first define the arithmetic sentence $\mathbf{Con}(T)$ in $L(T)$  as  $\neg \mathbf{Prov}_T(\ulcorner \mathbf{0}\neq\mathbf{0}\urcorner)$  which says that for all $x$, $x$ is not a code of
a proof of a contradiction in $T$. G\"{o}del's second incompleteness theorem ($\sf G2$) states that if $T$ is
consistent, then the arithmetical formula $\mathbf{Con}(T)$, which expresses the consistency of $T$, is not provable in $T$.  In Section \ref{intension problem of G2}, we will discuss some other ways of expressing   the consistency of $T$.

\smallskip

Finally, from the above conditions (1)-(3), one can show that $T\vdash \mathbf{Con}(T)\leftrightarrow \mathbf{G}$. Thus, $\sf G2$ holds: if $T$ is consistent, then $T\nvdash \mathbf{Con}(T)$.
%from any fixed point $\theta$ of $\neg\mathbf{Prov}_T(x)$ (i.e.~ $T \vdash \theta \leftrightarrow\neg \mathbf{Prov}_T(\ulcorner \theta \urcorner)$), we can show that  $T\nvdash \theta$ and $T\vdash \mathbf{Con}(T)\leftrightarrow \theta$ from conditions (1)-(3).
%For
For more details on these proofs of $\sf G1$ and $\sf G2$, we refer to Chapter 2 in \cite{metamathematics}.

%\section{Incompleteness and logical paradox}

\subsubsection{Properties of $\sf G1$}\label{waf1}
In this section, we discuss some (sometimes subtle) comments on $\sf G1$.

\smallskip

First of all, G\"{o}del's proof of $\sf G1$ is constructive as follows:  given a consistent r.e.~ extension $T$ of $\mathbf{PA}$, the proof constructs, in an algorithmic way, a true arithmetic sentence which is unprovable in $T$. In fact, one can effectively find a true $\Pi^0_1$ sentence $G_T$ of arithmetic such that $G_T$ is independent of $T$. G\"{o}del calls this the ``incompletability or inexhaustability of mathematics".

\smallskip

Secondly, for G\"{o}del's proof of $\sf G1$, only assuming that $T$ is consistent  does not suffice to show that G\"{o}del sentence is independent of $T$. In fact, the optimal condition to show that G\"{o}del sentence is independent of $T$ is:  $T+\mathbf{Con}(T)$ is consistent (see Theorems 35-36 in \cite{Isaacson 11}).\footnote{This optimal condition is much weaker than $\omega$-consistency.}

\smallskip

Thirdly, in summary, G\"{o}del's proof of $\sf G1$ has the following properties:
\begin{itemize}
  \item uses proof-theoretic method with arithmetization;
  \item does not directly use the Diagnolisation Lemma;
  \item the proof formalizes the liar paradox;
  \item the proof is constructive;
  \item the proof does not have the Rosser property;
  \item G\"{o}del's sentence has no real mathematical content.
\end{itemize}
All these characteristics of G\"{o}del's proof of $\sf G1$ are not necessary conditions for proving $\sf G1$.
For example, $\sf G1$ can be proved using recursion-theoretic or model-theoretic method, using the Diagnolisation Lemma, using other logical paradoxes, using non-constructive methods, only assuming that $T$ is consistent (i.e.~ having the Rosser property), and can be proved without arithmetization.

\smallskip

%Now, I will discuss more  different proofs of $\sf G1$ and $\sf G2$.
Fourth, $\sf G1$ does \emph{not} tell us that any consistent theory is incomplete. In fact, there are many consistent complete first-order theories. For example, the following first-order  theories are complete: the theory of dense linear orderings without endpoints ($\mathbf{DLO}$), the theory of ordered divisible groups ($\mathbf{ODG}$), the theory of algebraically closed fields of given characteristic ($\mathbf{ACF_p}$), and the theory of real  closed fields ($\mathbf{RCF}$). We refer to \cite{Epstein 2011} for details of these theories.
In fact, $\sf G1$  only  tells us that any consistent first-order  theory containing a large enough fragment of $\mathbf{PA}$ (such as $\mathbf{Q}$) is incomplete: there is then a true $\Pi^0_1$ sentence which is independent of the initial theory.
Turing's work in \cite{logic based on ordinals} shows that any  true $\Pi^0_1$-sentence of arithmetic is provable in some transfinite iteration of $\PA$. Feferman's work in \cite{recursive progressions} extends Turing's work and shows that any true sentence of arithmetic is provable in some transfinite iteration of $\mathbf{PA}$.

\smallskip

Fifth, whether a theory of arithmetic is complete depends on the language of the theory. There are respectively recursively axiomatized complete arithmetic theories in the language of $\, L(\mathbf{0}, \mathbf{S})$, $L(\mathbf{0}, \mathbf{S}, <)$  and  $L(\mathbf{0}, \mathbf{S}, <, +)$ (see Section 3.1-3.2 in \cite{Enderton 2001}).
Containing enough information of arithmetic is essential for a consistent arithmetic theory to be incomplete.  For example, Euclidean geometry is not about arithmetic but only about points, circles and lines in general; but  Euclidean geometry is complete as Tarski has proved (see \cite{Tarski 99}).
If the theory  contains only  information about the arithmetic of addition without multiplication, then it can be complete. For example, Presburger arithmetic is a complete theory of the arithmetic of addition in the language of $L(\mathbf{0}, \mathbf{S}, +)$ (see Theorem 3.2.2 in \cite{metamathematics}, p.222). Finally, containing the arithmetic of multiplication is not sufficient for a theory to be incomplete. For example, there exists a complete recursively axiomatized theory in the language of $L(\mathbf{0}, \times)$ (see \cite{metamathematics}, p.230).

\smallskip

Finally, it is well-known that $Th(\mathbb{N}, +, \times)$ is interpretable in $Th(\mathbb{Z}, +, \times)$ and $Th(\mathbb{Q}, +, \times)$.\footnote{The key point is: $\mathbb{N}$ is definable in $(\mathbb{Z}, +, \times)$ and $(\mathbb{Q}, +, \times)$. See chapter XVI in \cite{Epstein 2011}.} Since $Th(\mathbb{N}, +, \times)$ is undecidable and has a finitely axiomatizable incomplete sub-theory $\mathbf{Q}$, by Theorem \ref{interpretable theorem}, $Th(\mathbb{Z}, +, \times)$ and $Th(\mathbb{Q}, +, \times)$  are undecidable, and hence not recursively axiomatizable, but they respectively have a finitely axiomatizable incomplete sub-theory of integers and rational numbers. But $Th(\mathbb{R}, +, \times)$ is decidable and recursively axiomatizable  (even if not finitely axiomatizable). In fact, $Th(\mathbb{R}, +, \times)=\mathbf{RCF}$ (the theory of real closed field) (see \cite{Epstein 2011}, p.320-321). Note that this fact does not contradict $\sf G1$ since none of $\mathbb{N}, \mathbb{Z}$ and $\mathbb{Q}$ is definable in $(\mathbb{R}, +, \times)$.

\subsubsection{Between truth and provability}\label{waf3}

In this paper, unless stated otherwise, we equate a set of sentences with the set of G\"{o}del's numbers of these sentences.
We discuss the formalized notions of `truth' and `proof', and how they relate to incompleteness.
\smallskip

\begin{definition}\label{}
We define $\textbf{Truth}=\{\phi\in L(\mathbf{PA}): \mathfrak{N}\models\phi\}$ and $\textbf{Prov}=\{\phi\in L(\mathbf{PA}): \mathbf{PA}\vdash\phi\}$, i.e.\ the formalized notions of `proof' and `truth'.
\end{definition}

First of all, truth and provability are the same for \emph{purely existential statements}.  Put another way, incompleteness does not arise at the level of $\Sigma^0_1$ sentences.  Indeed, we have $\Sigma^0_1$-completeness for $T$: for any $\Sigma^0_1$ sentences $\phi$, $T\vdash \phi$ if and only if $\mathfrak{N}\models\phi$.
Thus, G\"{o}del's sentence is  a true $\Pi^0_1$ sentence in the form  $\forall x \phi(x)$ such that $T\nvdash \forall x \phi(x)$ but `$T\vdash\phi(\bar{n})$' holds for any $n\in\omega$.

\smallskip

Secondly,  the properties of $\textbf{Truth}$ are essentially different from that of $\textbf{Prov}$.
Before G\"{o}del's work, it was thought that $\textbf{Truth}=\textbf{Prov}$. Thus, G\"{o}del's first incompleteness theorem ($\sf G1$) reveals the difference between the notion of provability in $\PA$ and the notion of truth in the standard model of arithmetic $\mathfrak{N}$. There are some differences between  $\textbf{Truth}$ and $\textbf{Prov}$:
\begin{itemize}
  \item $\textbf{Prov}\subsetneq\textbf{Truth}$, i.e.~ there is a true arithmetic sentence   which is unprovable in $\mathbf{PA}$;
  \item Tarski proves that:  $\textbf{Truth}$ is not definable in $\mathfrak{N}$  but $\textbf{Prov}$ is definable in $\mathfrak{N}$;
  \item $\textbf{Truth}$ is not arithmetic but $\textbf{Prov}$ is recursive enumerable.
\end{itemize}
However, both $\textbf{Truth}$ and $\textbf{Prov}$  are not recursive and not representable in $\mathbf{PA}$. For more details on $\textbf{Truth}$ and $\textbf{Prov}$, we refer to \cite{metamathematics, undecidable}.

\smallskip

Thirdly, the differences between $\textbf{Truth}$ and $\textbf{Prov}$ can also be expressed in terms of \emph{arithmetical interpretations}, defined as follows.\smallskip

\begin{definition}[Arithmetical interpretations]~
A mapping from the set of all modal propositional variables to the set of $L(\mathbf{PA})$-sentences is called an \emph{arithmetical interpretation}.
\end{definition}
Every arithmetical interpretation $f$ is uniquely extended to the
mapping $f^{\ast}$ from the set of all modal formulas to the set of $L(T)$-sentences so
that $f^{\ast}$ satisfies the following conditions:
\begin{itemize}
  \item $f^{\ast}(p)=f(p)$ for each propositional variable $p$;
  \item $f^{\ast}$ commutes with every propositional connective;
  \item $f^{\ast}(\Box A)$ is $\mathbf{Prov}_T(\ulcorner f^{\ast}(A)\urcorner)$ for every modal formula $A$.
\end{itemize}
In the following, we equate arithmetical interpretations $f$ with their unique extensions $f^{\ast}$ defined on the set of all modal formulas.
In this way, Solovay's Arithmetical Completeness Theorems for $\mathbf{GL}$ and $\mathbf{GLS}$ characterize the difference between $\textbf{Prov}$ and $\textbf{Truth}$ via provability logic.\smallskip
\begin{theorem}[Solovay, \cite{Solovay 76}]~
\begin{description}
  \item[Arithmetical Completeness Theorem for $\mathbf{GL}$] Let $T$ be a $\Sigma^0_1$-sound r.e.~ extension of $\mathbf{Q}$. For any modal formula $\phi$ in $L(\mathbf{GL})$, $\mathbf{GL}\vdash\phi$ if and only if $T\vdash f(\phi)$ for every arithmetic interpretation $f$.
  \item[Arithmetical Completeness Theorem for $\mathbf{GLS}$] For any modal formula $\phi, \mathbf{GLS}\vdash\phi$ if and only if $\mathfrak{N}\models f(\phi)$ for every arithmetic interpretation $f$.
\end{description}
\end{theorem}

Finally, one can study the notion of `proof predicate' as given by $\mathbf{Proof}_T(x, y)$ in an abstract setting, namely as follows.
Recall that  $T$ is a recursively axiomatizable consistent extension of $\mathbf{Q}$.
We introduce  general notions of proof predicate and provability predicate  which generalize  the proof predicate $\mathbf{Proof}_T(x, y)$ and the provability predicate $\mathbf{Prov}_T(x)$ defined above in  G\"{o}del's proof of $\sf G1$.
\smallskip

\begin{definition}[Proof predicate]\label{def of proof predicate}
We say a formula $\mathbf{Prf}_T(x, y)$ is a \emph{proof predicate} of $T$ if it
satisfies the following conditions:\footnote{We can say that each proof predicate represents the relation ``$y$ is the code of a proof in $T$ of a formula with G\"{o}del number $x$".}
\begin{itemize}
  \item $\mathbf{Prf}_T(x, y)$ is $\Delta^0_1(\mathbf{PA})$;\footnote{We say a formula $\phi$ is $\Delta^0_1(\mathbf{PA})$ if there exists a $\Sigma^0_1$ formula $\alpha$ such that $\mathbf{PA}\vdash \phi\leftrightarrow\alpha$, and there exists a $\Pi^0_1$ formula $\beta$ such that $\mathbf{PA}\vdash \phi\leftrightarrow\beta$.}
  \item $\mathbf{PA} \vdash \forall x(\mathbf{Prov}_T(x) \leftrightarrow\exists y \mathbf{Prf}_T(x, y))$;
  \item  for any $n \in\omega$ and formula $\phi, \mathbb{N}\models \mathbf{Proof}_T(\ulcorner\phi\urcorner, \overline{n}) \leftrightarrow \mathbf{Prf}_T(\ulcorner\phi\urcorner, \overline{n})$;
      \item $\mathbf{PA} \vdash \forall x\forall x^{\prime} \forall y (\mathbf{Prf}_T(x, y) \wedge \mathbf{Prf}_T(x^{\prime}, y) \rightarrow x = x^{\prime})$.
\end{itemize}
\end{definition}

\begin{definition}[Provability and consistency]
We  define the  provability predicate $\mathbf{Pr}_T(x)$ from a proof predicate $\mathbf{Prf}_T(x,y)$ by $\exists y\, \mathbf{Prf}_T(x,y)$, and   the  consistency statement $\mathbf{Con}(T)$ from a provability predicate $\mathbf{Pr}_T(x)$ by $\neg \mathbf{Pr}_T (\ulcorner \mathbf{0}\neq\mathbf{0}\urcorner)$.
\end{definition}
The items $\mathbf{D1}$-$\mathbf{D3}$ below are called the \emph{Hilbert-Bernays-L\"{o}b derivability conditions}.
Note that $\mathbf{D1}$ holds for any provability predicate $\mathbf{Pr}_T(x)$.\smallskip
\begin{definition}[Standard proof predicate]\label{frak}
We say that provability predicate $\mathbf{Pr}_T(x)$ is  \emph{standard} if it satisfies  $\mathbf{D2}$ and $\mathbf{D3}$ as follows.
\begin{description}
  \item[D1] If $T \vdash\phi$, then $T \vdash \mathbf{Pr}_T(\ulcorner\phi\urcorner)$;
  \item[D2] If $T \vdash \mathbf{Pr}_T(\ulcorner\phi \rightarrow\varphi\urcorner) \rightarrow (\mathbf{Pr}_T(\ulcorner\phi\urcorner)\rightarrow \mathbf{Pr}_T(\ulcorner\varphi\urcorner))$;
  \item[D3] $T \vdash \mathbf{Pr}_T(\ulcorner\phi\urcorner)\rightarrow \mathbf{Pr}_T(\ulcorner \mathbf{Pr}_T(\ulcorner\phi\urcorner)\urcorner)$.
\end{description} We say that $\mathbf{Prf}_T(x, y)$ is a \emph{standard proof predicate}  if the induced provability predicate from it is standard.
\end{definition}
\noindent
The previous definition leads to another blanket caveat:

\smallskip

\emph{Unless stated otherwise,  we always assume that $\mathbf{Pr}_T(x)$ is a standard provability predicate, and $\mathbf{Con}(T)$ is the canonical consistency
statement defined as $\neg\mathbf{Pr}_T(\ulcorner \mathbf{0}\neq \mathbf{0}\urcorner)$ via the  standard provability predicate $\mathbf{Pr}_T(x)$.}

\subsubsection{Properties of $\sf G2$}\label{waf2}
In this section, we discuss some (sometimes subtle) comments on $\sf G2$.

\smallskip

First of all, we examine a somewhat delicate mistake in the argument which claims that, by an easy application of the compactness theorem, we can show that for any recursive axiomatization of a consistent theory $T$, $T$ can not prove its own consistency. Visser presents this argument in \cite{Visser 19} as an interesting dialogue between Alcibiades and Socrates:
\begin{quote}
Suppose a consistent theory $T$ can  prove its own consistency under some  axiomatization.
By compactness theorem, there must be a finitely axiomatized sub-theory $S$ of $T$ such that $S$ already proves the consistency of  $T$. Since $S$ proves the consistency
of $T$, it must also prove the consistency of $S$. So, we have a finitely
axiomatized theory which proves its own consistency. But $\sf G2$ applies  to  the finite
axiomatization  and we have a contradiction. It follows that $T$ can not prove its own consistency.
\end{quote}
The  mistake in this argument is: from the fact that $S$ can prove the consistency of $T$ we cannot infer that $S$ can prove the consistency of $S$. Some may argue that since $S$ is a sub-theory of $T$ and $S$ can prove the consistency of $T$, then of course $S$ can prove the consistency of $S$.

\smallskip

However, as Visser correctly points out in \cite{Visser 19}, we should carefully distinguish three perspectives  of the theory $T$: our external perspective, the internal perspective of
$S$, and the internal perspective of $T$. From each perspective, the consistency of the whole theory implies the consistency of its sub-theory. From $T$'s perspective, $S$ is a sub-theory of $T$. But from $S$'s perspective, $S$ may not be a sub-theory of $T$. From the fact that $T$ knows that $S$ is a sub-theory of $T$, we cannot infer that $S$ also knows that $S$ is a sub-theory of $T$ since $S$ is a finite sub-theory of $T$ and may not know any information that $T$ knows, leading to the following (dramatic) conclusion:
\begin{center}
\emph{the sub-theory relation between theories is not absolute.}
\end{center}
Similarly, the notion of consistency is not absolute. For example, let $S=\mathbf{PA}+ \neg \mathbf{Con}(\mathbf{PA})$. From $\sf G2$, $S$ is consistent from the external perspective. But since $S\vdash \neg \mathbf{Con}(S)$, the  theory $S$ is not consistent from the internal perspective of $S$. Note that $\mathbf{PA}\vdash \mathbf{Pr}_{\mathbf{PA}}(\mathbf{0}\neq\mathbf{0})\rightarrow \mathbf{Pr}_{\mathbf{PA}}(\mathbf{Pr}_{\mathbf{PA}}(\mathbf{0}\neq\mathbf{0})\rightarrow \mathbf{0}\neq\mathbf{0})$. Thus, a theory may be consistent from the external perspective but inconsistent from the internal perspective.
\smallskip

From G\"{o}del's proof of $\sf G2$, we cannot infer that if $T$ is a consistent r.e.~ extension of $\mathbf{Q}$, then $\mathbf{Con}(T)$ is independent of $T$.  The key point is: it is not enough to show that $T\nvdash\neg\mathbf{Con}(T)$  only assuming that $T$ is consistent. However, we can show that $\mathbf{Con}(T)$ is independent of $T$ assuming that $T$ is 1-consistent.\footnote{It is an easy fact that if $T$ is 1-consistent and $S$ is not a theorem of $T$, then $\mathbf{Pr}_{T}(\ulcorner S\urcorner)$ is not a theorem of $T$.} In fact, the formalized version of ``if $T$ is consistent, then $\mathbf{Con}(T)$ is independent of $T$" is not provable in $T$.\footnote{See \cite[Theorem 4, p.97]{Boolos 93}  for a modal proof in $\mathbf{GL}$ of this fact using the Arithmetic Completeness Theorem for $\mathbf{GL}$.}
\smallskip
\begin{definition}[Reflexivity]~
\begin{itemize}
\item A first-order theory $T$ containing $\mathbf{PA}$ is said
to be \emph{reflexive} if $T \vdash \mathbf{Con}(S)$ for each finite sub-theory $S$ of  $T$  where $\mathbf{Con}(S)$ is similarly defined as $\mathbf{Con}(\mathbf{PA})$.
\item We say the theory $T$ is \emph{essentially reflexive}  if any consistent extension of $T$ in $L(T)$ is reflexive.
\item Let $\mathbf{Con}(T)\!\!\upharpoonright x$ denote the finite consistency statement ``there
are no proofs of contradiction in $T$ with $\leq x$ symbols".
\end{itemize}
\end{definition}
Mostowski proves that $\mathbf{PA}$ is essentially reflexive (see \cite[Theorem 2.6.12]{metamathematics}). In fact one can show that for every $n \in \mathbb{N}$, $I\Sigma^0_{n+1}\vdash \mathbf{Con}(I\Sigma^0_{n})$.\footnote{For a proof of this result, we refer to H\'{a}jek  and Pudl\'{a}k \cite{Pudlak 93}.}
For a large class of
natural theories $U$, Pudl\'{a}k \cite{Pudlak 17} shows that the lengths of the shortest proofs
of $\mathbf{Con}(U)\!\!\upharpoonright n$ for $n\in\omega$ in the theory $U$ itself are bounded by a polynomial in $n$. Pudl\'{a}k conjectures  \cite{Pudlak 17} that $U$ does not have polynomial proofs of the finite
consistency statements $\mathbf{Con}(U + \mathbf{Con}(U))\upharpoonright n$ for $n\in\omega$.  % (see \cite{Pudlak 17}).

\smallskip

Finally, a big open question about ${\sf G2}$ is: can we find a genuinely self-reference free proof of ${\sf G2}$? As far as we know, at present there is no  convincing essentially self-reference-free proofs of either ${\sf G2}$ or of Tarski's Theorem of the Undefinability of Truth.  In \cite{Visser 19 Tarski to  Godel}, Visser gives a self-reference-free proof of ${\sf G2}$ from Tarski's Theorem of the Undefinability of Truth, which is a step in a program to find self-reference-free proofs  of both ${\sf G2}$ and  Tarski's Theorem (see \cite{Visser 19 Tarski to  Godel}). Visser's  argument in \cite{Visser 19 Tarski to  Godel} is model-theoretic and the main tool is the Interpretation Existence Lemma.\footnote{We refer to \cite{Visser 17 interpretation existence lemma} for more details about the Interpretation Existence Lemma.} Visser's proof in \cite{Visser 19 Tarski to  Godel} is not constructive. An interesting question is then whether Visser's argument can be made constructive.

\subsection{Proofs of $\sf G1$ and $\sf G2$ from mathematical logic}\label{foef}
In this section, we discuss various different proofs of $\sf G1$ and $\sf G2$.
We mention Jech's \cite{Jech 94} short proof of ${\sf G2}$ for $\mathbf{ZF}$: if $\mathbf{ZF}$ is consistent, then it is unprovable in $\mathbf{ZF}$  that there exists a model of $\mathbf{ZF}$. Jech's proof  uses the Completeness
Theorem, and also yields ${\sf G2}$ for $\mathbf{PA}$ (see \cite{Jech 94}).  Other (lengthier) proofs are discussed in Sections \ref{flk}-\ref{flk2}.

\subsubsection{Rosser's proof}\label{flk}
Rosser \cite{Rosser 36} proves a ``stronger" version of $\sf G1$, called Rosser's first incompleteness theorem, which only assumes the consistency of $T$: if $T$ is a consistent r.e.~ extension of $\mathbf{Q}$, then $T$ is incomplete. G\"{o}del's proof of $\sf G1$ assumes that $T$ is $\omega$-consistent. Note that $\omega$-consistency implies consistency. But the converse does not hold and
the notion of $\omega$-consistency is stronger than consistency since we can find examples of theories that are consistent but not $\omega$-consistent.\footnote{For example, assuming $\mathbf{PA}$ is consistent, then $\mathbf{PA} + \neg \mathbf{Con(PA)}$ is consistent, but not $\omega$-consistent.}
Rosser's proof is  constructive and algorithmically constructs the Rosser sentence that is independent of $T$. G\"{o}del's proof of $\sf G1$ uses a standard provability predicate but Rosser's proof of $\sf G1$ uses a Rosser provability predicate which is a kind of \emph{non-standard} provability predicate, giving rise to the following.
\begin{definition}\label{Rosser predicate}
Let $T$ be a recursively axiomatizable consistent extension of $\mathbf{Q}$, and $\mathbf{Prf}_T(x, y)$ be any proof predicate of $T$. Define the Rosser provability predicate $\mathbf{Pr}_T^R(x)$ to be the formula $\exists y(\mathbf{Prf}_T(x, y) \wedge \forall z \leq y \neg \mathbf{Prf}_T(\dot{\neg} x, z))$ where $\dot{\neg}$ is a function symbol expressing a primitive recursive function calculating the code of $\neg \phi$ from the code of $\phi$. The fixed point of the predicate $\neg \mathbf{Pr}_T^R(x)$ is called the Rosser sentence of $\mathbf{Pr}_T^R(x)$, i.e.~ a sentence $\theta$ satisfying $\mathbf{PA} \vdash \theta \leftrightarrow\neg \mathbf{Pr}_T^R(\ulcorner\theta\urcorner)$.
\end{definition}
In general, one can show that each Rosser sentence based on any Rosser provability predicate of $T$ is independent of $T$.  In particular, this independence does not rely on the choice of the proof predicate.

\subsubsection{Recursion-theoretic proofs}
G\"{o}del's first incompleteness theorem ($\sf G1$)  is well-known in the context of recursion
theory. Recall that $W_e=\{n\in\omega: \phi_e(n)\!\!\downarrow\}$. Let $\langle W_e: e\in\omega\rangle$ be the list of recursive enumerable subsets of $\mathbb{N}$.
The following is an example of an `effective' version of $\sf G1$:
\begin{center}
there exists a recursive function $f$ such that for any $e\in\omega$, if $W_e\subseteq \textbf{Truth}$, then $f(e)$ is defined and $f(e)\in \textbf{Truth}\setminus W_e$ (\cite{Enderton 2011}).
\end{center}
Similarly, Avigad \cite{Avigad 05} proves $\sf G1$ and $\sf G2$ in terms of the undecidability of the halting problem (see Theorem 3.1, Theorem 3.2  in \cite{Avigad 05}).
Another related result due to Kleene is as follows.
\begin{theorem}[Kleene's theorem, Theorem 2.2, \cite{Salehi-Seraji 18}]
For any consistent r.e.~ theory $T$ that contains $\mathbf{Q}$, there exists some $t \in \omega$ such that $\varphi_t(t)\!\uparrow$ holds but $T \nvdash ``\varphi_t(t)\!\uparrow"$.
\end{theorem}

Kleene's proof of his theorem uses recursion theory, and is not constructive. Salehi and Seraji \cite{Salehi-Seraji 18} show that there is a constructive proof of Kleene's theorem, but this constructive proof does not have the Rosser property. Salehi and Seraji \cite{Salehi-Seraji 18} comment that there could be a `Rosserian' version of this constructive proof of Kleene's theorem.

\subsubsection{Proofs based on Arithmetic Completeness}
Hilbert and Bernays \cite{Hilbert  Bernays 39} present the \emph{Arithmetic Completeness Theorem} expressing that any recursively axiomatizable consistent theory has an arithmetically definable model.
Later, Kreisel  \cite{Kreisel 50} and Wang \cite{Wang 55} adapt the Arithmetic Completeness Theorem and use paradoxes to obtain undecidability results.

\smallskip

Now, the Arithmetic Completeness Theorem is an important tool in model-theoretic proofs of the incompleteness theorems. For more details, we refer to \cite{Aspects of incompleteness, Kaye-Kotlarski 00, Henryk Kotlarski 04}.
Walter Dean \cite{Walter Dean} gives a detailed discussion on how the Arithmetized Completeness Theorem provides a tool for obtaining formal incompleteness results from some certain paradoxes.
\begin{theorem}[Arithmetic Completeness, Theorem 3.1, \cite{Makoto Kikuchi 94}]\label{ACT}
Let $T$ be a recursively axiomatized consistent extension  of $\mathbf{Q}$. There exists a formula
$\mathbf{Tr}_T(x)$ in  $L(\mathbf{PA})$ that defines a model of $T$ in $\mathbf{PA}+\mathbf{Con}(T)$.
\end{theorem}
Lemma \ref{key lemma for model theory} is a corollary of the Arithmetized Completeness Theorem, and is essential for model-theoretic proofs of the incompleteness theorems.
\begin{lemma}[\cite{Kikuchi-Tanaka 94, Makoto Kikuchi 94}]\label{key lemma for model theory}
Let $T$ be a recursively axiomatized consistent extension  of $\mathbf{Q}$, and $\mathbf{Tr}_T(x)$ is the formula as asserted in Theorem \ref{ACT}. For any model $M_0$ of $\mathbf{PA}+\mathbf{Con}(T)$, there exists a model $M_1$ of $T$ such that  for any sentence $\phi$,  $M_0\models \mathbf{Tr}_T(\ulcorner \phi\urcorner)$ if and only if $M_1\models\phi$.
\end{lemma}
Kreisel first applies the Arithmetized Completeness Theorem to establish model-theoretic proofs of ${\sf G2}$ (cf.~ Kreisel \cite{Kreisel 68}, Smory\'{n}ski \cite{Smorynski 1977} and Kikuchi \cite{Makoto Kikuchi 94}).
Kikuchi-Tanaka \cite{Kikuchi-Tanaka 94}, Kikuchi \cite{Makoto Kikuchi 94, Kikuchi 97} and
Kotlarski \cite{Henryk Kotlarski 94} use  the Arithmetized Completeness Theorem to give  model-theoretic proofs of  ${\sf G2}$. For example,
Kikuchi \cite{Makoto Kikuchi 94} proves ${\sf G2}$ model-theoretically via the Arithmetized Completeness Theorem (Lemma \ref{key lemma for model theory}):  if $\mathbf{PA}$ is consistent, then $\mathbf{Con}(\mathbf{PA})$ is not provable in $\mathbf{PA}$ (see Theorem 3.4, \cite{Makoto Kikuchi 94}).\footnote{The idea of the proof is: assuming that $\mathbf{PA}$ is consistent and $\mathbf{PA}\vdash\mathbf{Con(\mathbf{PA})}$, then we get a contradiction from the fact that there is a model $M$ of $\mathbf{PA}$ such that $M\models \mathbf{Con(PA)}$.}

\smallskip

Proofs of ${\sf G2}$ by Kreisel \cite{Kreisel 68} and  Kikuchi \cite{Makoto Kikuchi 94} do not directly yield the formalized version of ${\sf G2}$.
Kikuchi's proof of ${\sf G2}$ in \cite{Kikuchi 97} is not formalizable in  $\mathbf{PRA}$. Kikuchi and
Tanaka \cite{Kikuchi-Tanaka 94}  prove
in $\mathbf{WKL_0}$ that $\mathbf{Con}(\mathbf{PA})$ implies $\neg \mathbf{Pr}_{\mathbf{PA}}(\ulcorner \mathbf{Con}(\mathbf{PA})\urcorner)$, since the Completeness Theorem is provable in $\mathbf{WKL_0}$, and  the key Lemma \ref{key lemma for model theory} used in Kikuchi's proof \cite{Makoto Kikuchi 94} is provable in $\mathbf{RCA_0}$.\footnote{The theory $\mathbf{RCA_0}$ (Recursive Comprehension) is a subsystem of Second Order Arithmetic.
For the definition of $\mathbf{RCA_0}$, we refer to \cite{SOA}.}
Using Theorem \ref{conservation thm}, Tanaka \cite{Kikuchi-Tanaka 94} proves the formalized version of  ${\sf G2}$: $\mathbf{PRA}\vdash \mathbf{Con}(\mathbf{PA})\rightarrow \mathbf{Con}(\mathbf{PA}+\neg \mathbf{Con}(\mathbf{PA}))$.

\smallskip

One can give a simple proof of $\sf G1$ via the Diagnolisation Lemma (see \cite{metamathematics}). Kotlarski \cite{Henryk Kotlarski 94} proves the formalized version of  $\sf G1$ and $\sf G2$ via model-theoretic arguments (e.g.~ using the Arithmetized Completeness Theorem and some quickly growing functions).
Kotlarski \cite{Henryk Kotlarski 94} proves the following version of $\sf G1$ assuming that $\mathbf{PA}$ is $\omega$-consistent, and shows that the following sentence is provable in $\mathbf{PA}$:
\begin{center}
if $\forall\varphi, x \{[\varphi\in \Delta_0\wedge \forall y \mathbf{Pr}_{\mathbf{PA}}(\neg \varphi(S^x 0, S^y 0))]\rightarrow \neg \mathbf{Pr}_{\mathbf{PA}}(\exists y \varphi(S^x 0, y))\}$, then $\exists\varphi\in \Delta_0\exists x [\neg \mathbf{Pr}_{\mathbf{PA}}(\exists y \varphi(S^x 0, y))\wedge\neg \mathbf{Pr}_{\mathbf{PA}}(\neg\exists y \varphi(S^x 0, y))]$.
\end{center}

However, it is unknown whether the method  in \cite{Henryk Kotlarski 94} can also give a new proof of Rosser's first incompleteness theorem.
Kotlarski \cite{Henryk Kotlarski 94} proves the following formalized version of  $\sf G2$: $\mathbf{PA}\vdash\mathbf{Con(\mathbf{PA})} \rightarrow \mathbf{Con}(\mathbf{PA}+\neg \mathbf{Con}(\mathbf{PA}))$.
Later, Kotlarski \cite{Henryk Kotlarski 98} transforms the proof of the  formalized version of $\sf G2$  in \cite{Henryk Kotlarski 94} to a proof-theoretic version without the use of the Arithmetized Completeness Theorem.

\subsubsection{Proofs based on Kolmogorov complexity}
Intuitively, \emph{Kolmogorov complexity} is a measure of the quantity of information in finite
objects. Roughly speaking, the Kolmogorov complexity of a number $n$, denoted
by $K(n)$, is the size of a program which generates $n$.
\begin{definition}[Kolmogorov-Chaitin Complexity, \cite{Salehi-Seraji 18}]
For any natural number $n \in\omega$, the \emph{Kolmogorov complexity} for $n$, denoted by $K(n)$, is defined as $\min\{i \in \omega \mid \varphi_i(0)\!\downarrow = n\}$.
\end{definition}
If $n \leq K(n)$, then $n$ is called random. Kolmogorov shows in 1960's that the set of non-random numbers is
recursively enumerable but not recursive (c.f.~ Odifreddi \cite{Odifreddi 89}). Relations between $\sf G1$ and Kolmogorov complexity have been intensively discussed in the literature (c.f.~ Li and Vit\'{a}nyi \cite{Li 90}).
Chaitin \cite{Chaitin 74} gives an information-theoretic
formulation of $\sf G1$,  and proves  the following weaker version of $\sf G1$ in terms of Kolmogorov complexity.

\begin{theorem}[Chaitin \cite{Chaitin 74, Salehi-Seraji 18}]\label{Chaitin thm}
For any consistent r.e.~ extension $T$ of $\mathbf{Q}$, there exists a constant
$c_T \in \mathbb{N}$ such that for any $e \geq c_T$ and any $w \in \mathbb{N}$ we have $T \nvdash ``K(w) > e"$.
\end{theorem}
Salehi and Seraji \cite{Salehi-Seraji 18} show that  we can algorithmically construct the Chaitin constant $c_T$ in Theorem \ref{Chaitin thm}. I.e.~ for a given consistent  r.e.~ extension $T$ of $\mathbf{Q}$, one can
algorithmically construct a constant $c_T\in \mathbb{N}$ such that for all $e \geq c_T$ and all $w \in \mathbb{N}$, we have $T \nvdash ``K(w) > e"$ (see Theorem 3.4 in \cite{Salehi-Seraji 18}). From Theorem \ref{Chaitin thm}, it is not clear whether ``$K(w) > e$" holds (or whether ``$K(w) > e$" is independent of $T$). Salehi and Seraji \cite{Salehi-Seraji 18} show that Chaitin's proof of ${\sf G1}$ is non-constructive: there is no algorithm such that  given any consistent r.e.~ extension $T$ of $\mathbf{Q}$ we can compute some $w_T$ such that $K(w_T) > c_T$ holds where $c_T$ is the Chaitin constant we can compute as in Theorem \ref{Chaitin thm} (see Theorem 3.5, \cite{Salehi-Seraji 18}). If such an algorithm exists, then for any  consistent r.e.~ extension $T$ of $\mathbf{Q}$, we can compute some $c_T$ and $w_T$ such that $K(w_T) > c_T$ is true but unprovable in $T$.

\smallskip

Salehi and Seraji \cite{Salehi-Seraji 18} also strengthen Chaitin's Theorem \ref{Chaitin thm} assuming $T$ is $\Sigma^0_1$-sound:
 if $T$ is a $\Sigma^0_1$-sound r.e.~ theory extending  $\mathbf{Q}$, then there
exists some $c_T$ (which is computable from $T$) such that for any $e \geq c_T$ there are cofinitely many $w$'s such
that ``$K(w) > e$" is independent of $T$ (see Corollary 3.7, \cite{Salehi-Seraji 18}).
Using a version of the
Pigeonhole Principle in $\mathbf{Q}$, Salehi and Seraji \cite{Salehi-Seraji 18} also prove the Rosserian form of Chaitin's Theorem: for any consistent r.e.~ extension $T$ of $\mathbf{Q}$, there is a
constant $c_T$ (which is computable from $T$) such that for any $e \geq c_T$ there are cofinitely many $w$'s such that
``$K(w) > e$" is independent of $T$ (see Theorem 3.9, \cite{Salehi-Seraji 18}).

\smallskip

Kikuchi \cite{Kikuchi 97} proves the following formalized version of ${\sf G1}$ via Kolmogorov complexity for any consistent r.e.~ extension $T$ of $\mathbf{PA}$: there exists $e \in \omega$ with
\begin{itemize}
  \item $T \vdash \mathbf{Con}(T) \rightarrow \forall x(\neg \mathbf{Pr}_T(\ulcorner K (x)>e\urcorner ))$;
  \item $T \vdash \omega$-$\mathbf{Con}(T) \rightarrow \forall x(e < K(x) \rightarrow \neg \mathbf{Pr}_T(\ulcorner K (x) \leq e\urcorner))$.
\end{itemize}
However, this proof is not constructive. Moreover, Kikuchi \cite{Kikuchi 97} proves ${\sf G2}$ via Kolmogorov complexity and the Arithmetic Completeness Theorem:
if  $T$ is a consistent r.e.~ extension of $\mathbf{PA}$, then $T\nvdash\mathbf{Con}(T)$.
Kikuchi's proof of ${\sf G2}$ in \cite{Kikuchi 97} cannot be formalized in  $\mathbf{PRA}$
but can be carried out within  $\mathbf{WKL_0}$. Thus we can also obtain a formalized version of ${\sf G2}$ in  $\mathbf{WKL_0}$ by Theorem \ref{conservation thm}.

\subsubsection{Model-theoretic proofs}\label{flk2}
Adamowicz and  Bigorajska \cite{Adamowicz Bigorajska 01} prove
 ${\sf G2}$ via model-theoretic method using the notion of 1-closed models and  existentially closed models.
\smallskip

\begin{definition}~
\begin{itemize}
  \item A model $M$ of a theory $T$ is called 1-closed (w.r.t.~ $T$) if for any $a_1,\cdots, a_n$ in $M$, any $\Sigma_1$ formula $\phi$ and any $M^{\prime}$ such that $M \prec_0 M^{\prime}$ and $M^{\prime}\models T$,  we have: if $M^{\prime}\models \phi(a_1,\cdots, a_n)$, then $M \models \phi(a_1,\cdots, a_n)$. In other words, we can say
 that $M$ is 1-closed if for any $M^{\prime}$ such that $M \prec_0 M^{\prime}$, we have $M \prec_1 M^{\prime}$.
\item Let $\mathcal{K}$ be a class of structures in the same language.
A model $M\in\mathcal{K}$  is
\emph{existentially closed} in $\mathcal{K}$
if for every model $N\supseteq M$ such that $N\in\mathcal{K}$, we have $M\preceq_1 N$: every existential formula
with parameters from $M$ which is satisfied in $N$ is already satisfied in $M$.
\end{itemize}
\end{definition}

Adamowicz and  Bigorajska \cite{Adamowicz Bigorajska 01} first prove ${\sf G2}$ without the use of the Arithmetized Completeness
 Theorem:   every 1-closed
 model of any subtheory $T$ of $\mathbf{PA}$ extending $I\Delta_0 + \mathbf{exp}$ satisfies $\neg \mathbf{Con}(\mathbf{PA})$.
Then Adamowicz and  Bigorajska \cite{Adamowicz Bigorajska 01}  prove
the formalized version of ${\sf G2}$ via the idea of existentially closed models and the Arithmetized Completeness
 Theorem: $\mathbf{PA} \vdash \mathbf{Con}(\mathbf{PA}) \rightarrow \mathbf{Con}(\mathbf{PA} + \neg \mathbf{Con}(\mathbf{PA}))$ (see Theorem 2.1, \cite{Adamowicz Bigorajska 01}). This is proved by showing that an arbitrary model of $\mathbf{PA} + \mathbf{Con}(\mathbf{PA})$ satisfies $\mathbf{Con}(\mathbf{PA} + \neg \mathbf{Con(PA)})$.

\subsection{Proofs of ${\sf G1}$ and ${\sf G2}$ based on logical paradox}\label{woef}
We provide a survey of proofs of incompleteness theorems based on `logical paradox'.
\subsubsection{Introduction}
As noted in Section \ref{defi2}, G\"{o}del's incompleteness theorems are closely related to paradox and self-reference. In fact,
G\"{o}del comments in his famous paper \cite{Godel paradox} that ``any epistemological antinomy could be used for a similar proof of the existence of undecidable propositions".

\smallskip

Now, the \emph{Liar Paradox} is an old and most famous paradox in modern science.
In G\"{o}del's proof of ${\sf G1}$, we can view G\"{o}del's sentence as the formalization of the Liar Paradox.
G\"{o}del's sentence concerns the notion of provability but the liar sentence in the Liar Paradox concerns the notion of truth in the standard model of arithmetic.
There is a big difference between  the notion of provability and truth. G\"{o}del's sentence  does not lead to a contradiction as the Liar sentence does.

\smallskip

Besides the Liar Paradox, many other paradoxes have been used to give new proofs of incompleteness theorems: for example, Berry's Paradox in \cite{BOOLOS 89, Chaitin 74, Makoto Kikuchi 94, Kikuchi-Kurahashi-Sakai, Kikuchi-Tanaka 94, Vopenka 66}, Grelling-Nelson's Paradox in \cite{Heterologicality and incompleteness}, the Unexpected Examination Paradox in \cite{Fitch 64, Kritchman-Raz 10}, and Yablo's Paradox in \cite{Godelizing the Yablo sequence, Kikuchi-Kurahashi 11, Kurahashi 14, Priest 97}.  We now discuss some of these paradoxes in detail.

\subsubsection{Berry's paradox}
Berry's Paradox introduced
by Russell \cite{Russell 08} is the
paradox that ``the least integer not nameable in fewer than nineteen syllables"
is itself a name consisting of eighteen syllables.
Informally, we say that an expression names a natural number $n$ if $n$ is the
unique natural number satisfying the expression. Berry's Paradox can be
formalized in formal systems by interpreting the concept of ``name" suitably. The
following is Boolos's formulation of the concept of ``name" in \cite{BOOLOS 89}.
\smallskip

\begin{definition}[Boolos \cite{BOOLOS 89}]
Let $n \in\omega$  and $\varphi(x)$ be a formula with only one
free variable $x$. We say that $\varphi(x)$ names $n$ if $\mathfrak{N}\models\varphi(\overline{n})\wedge \forall v_0\forall v_1(\varphi(v_0)\wedge\varphi(v_1)\rightarrow v_0= v_1)$.
\end{definition}
Proofs of the incompleteness theorems based on  Berry's Paradox have been given  by Vop\v{e}nka \cite{Vopenka 66}, Chaitin \cite{Chaitin 74}, Boolos \cite{BOOLOS 89}, Kikuchi-Kurahashi-Sakai \cite{Kikuchi-Kurahashi-Sakai},  Kikuchi \cite{Makoto Kikuchi 94}, and Kikuchi-Tanaka \cite{Kikuchi-Tanaka 94}. In fact, Robinson first uses    Berry's Paradox in \cite{Robinson 63} to prove Tarski's theorem on the undefinability of truth, which anticipates the later use of  Berry's Paradox to obtain incompleteness results by Vop\v{e}nka \cite{Vopenka 66}, Boolos \cite{BOOLOS 89} and Kikuchi \cite{Kikuchi 97}.
\smallskip

Boolos \cite{BOOLOS 89} proves a weak form of ${\sf G1}$ in the 1980's by formalizing Berry's Paradox in arithmetic via considering the length of formulas that name natural numbers in the standard model of arithmetic.
Using this formulation of the concept of ``name", Boolos \cite{BOOLOS 89} first shows that Berry's Paradox leads to a proof of ${\sf G1}$ in the following form: there is no algorithm whose output
contains all true statements of arithmetic and no false ones (i.e.~ the theory of true
arithmetic is not recursively axiomatizable).
Barwise \cite{Barwise 89} praises Boolos's proof as ``very lovely and the most straightforward proof of G\"{o}del's incompleteness theorem that I have ever seen". The optimal sufficient and necessary condition for the independence of a Boolos sentence from $\mathbf{PA}$  is that $\mathbf{PA}+\mathbf{Con}(\mathbf{PA})$ is consistent (see \cite{Salehi-Seraji 18}).
\smallskip

Boolos's theorem is different from G\"{o}del's theorem in the following way:
\begin{itemize}
  \item Boolos's theorem refers to the concept of truth but G\"{o}del's theorem does not;
  \item Boolos's proof is not constructive, and we can prove that there is no algorithm for computing the true but unprovable sentence;
  \item Boolos's theorem is weaker than G\"{o}del's first incompleteness theorem, and hence we cannot obtain the second incompleteness
theorem from Boolos's theorem in the standard way (see \cite{Kikuchi-Kurahashi-Sakai}).
\end{itemize}
Boolos's proof is modified by Kikuchi and Tanaka in \cite{Kikuchi-Tanaka 94, Makoto Kikuchi 94}.
The difference between  Kikuchi's proof and Boolos's proof lies in the interpretation of the word ``name".
Kikuchi \cite{Makoto Kikuchi 94} modifies Boolos's formulation of the concept of ``name" by replacing
``truth"  with ``provability" in the definition.
\smallskip

\begin{definition}[Definition 3.1, Kikuchi \cite{Makoto Kikuchi 94}]~
Let $n \in\omega$  and $\varphi(x)$ be a formula with only one
free variable $x$. We say that $\varphi(x)$ names $n$ if $\mathbf{PA}\vdash \varphi(\overline{n})\wedge \forall v_0\forall v_1(\varphi(v_0)\wedge\varphi(v_1)\rightarrow v_0= v_1)$.
\end{definition}

Using this formulation of the concept of ``name", Kikuchi \cite{Makoto Kikuchi 94}
gives a proof-theoretic proof of ${\sf G1}$ by formalizing Berry's paradox without the use of the Diagnolisation Lemma. Kikuchi \cite{Makoto Kikuchi 94} constructs a sentence $\theta$  and shows  that  if $\mathbf{PA}$ is consistent, then $\neg \theta$ is not provable in $\mathbf{PA}$; if $\mathbf{PA}$ is $\omega$-consistent, then  $\theta$ is not provable in $\mathbf{PA}$ (see Theorem 2.2, \cite{Makoto Kikuchi 94}). Note that Kikuchi's proof of ${\sf G1}$ in \cite{Makoto Kikuchi 94} is constructive.
Kikuchi and Tanaka \cite{Kikuchi-Tanaka 94} reformulate Kikuchi's proof of ${\sf G1}$ in \cite{Makoto Kikuchi 94},  and show in $\mathbf{WKL_0}$ that if $\mathbf{PA}+\mathbf{Con}(\mathbf{PA})$ is consistent, then $\theta$ is independent of $\mathbf{PA}$.  By Theorem \ref{conservation thm}, Kikuchi and
Tanaka \cite{Kikuchi-Tanaka 94} prove the formalized version of ${\sf G1}$: $\mathbf{PRA}\vdash \mathbf{Con}(\mathbf{PA}+\mathbf{Con}(\mathbf{PA}))\rightarrow \neg \mathbf{Pr}_{\mathbf{PA}}(\ulcorner \theta\urcorner)\wedge \neg \mathbf{Pr}_{\mathbf{PA}}(\ulcorner \neg\theta\urcorner)$. An interesting question not covered in \cite{Makoto Kikuchi 94, Kikuchi-Tanaka 94}  is whether we can improve Kikuchi's proof of ${\sf G1}$ by only assuming that $\mathbf{PA}$ is consistent.
\smallskip

Vop\v{e}nka \cite{Vopenka 66} proves ${\sf G2}$ for $\mathbf{ZF}$ by formalizing Berry's Paradox, via  adopting Kikuchi's definition of the
concept of ``name" in  \cite{Makoto Kikuchi 94} over models of $\mathbf{ZF}$\footnote{I.e.~  we say that $\varphi(x)$ names $n$ in $\mathbf{ZF}$ if  $\mathbf{ZF}\vdash \varphi(\overline{n})\wedge \forall v_0\forall v_1(\varphi(v_0)\wedge\varphi(v_1)\rightarrow v_0= v_1)$ where $\varphi(x)$ is a formula with only one
free variable $x$ (see \cite{Vopenka 66}).}: $\mathbf{Con(ZF)}$ is not provable in $\mathbf{ZF}$. Vop\v{e}nka's proof uses the Completeness Theorem but does not use the Arithmetic Completeness Theorem.
Kikuchi, Kurahashi and Sakai \cite{Kikuchi-Kurahashi-Sakai} show that Vop\v{e}nka's method can be adapted to prove ${\sf G2}$ for $\mathbf{PA}$ based on Kikuchi's formalization of Berry's Paradox
in \cite{Makoto Kikuchi 94} with an application of the Arithmetic Completeness Theorem.

\smallskip

Proofs of ${\sf G1}$ and ${\sf G2}$ based on Berry's Paradox by Vop\v{e}nka \cite{Vopenka 66}, Chaitin \cite{Chaitin 74}, Boolos \cite{BOOLOS 89} and Kikuchi \cite{Makoto Kikuchi 94} do not use the Diagnolisation  Lemma. We can also prove ${\sf G1}$ based on Berry's Paradox using  the  Diagnolisation  Lemma. For example, Kikuchi, Kurahashi and Sakai \cite{Kikuchi-Kurahashi-Sakai} adopt Kikuchi's definition of the concept of ``name" in \cite{Makoto Kikuchi 94}, and show that the independent statement in Kikuchi's
proof in \cite{Makoto Kikuchi 94} can be obtained by using the  Diagnolisation  Lemma.

\smallskip

In summary, the distinctions between using and not using the Diagnolisation  Lemma, and between using
and not using the Arithmetic Completeness Theorem are not essential for proofs of ${\sf G1}$ and ${\sf G2}$ based on Berry's Paradox.
From the above discussions, we can characterize different proofs of ${\sf G1}$ and ${\sf G2}$ based on Berry's Paradox by the method of
interpreting  the word ``name": Boolos \cite{BOOLOS 89} uses
the standard model of arithmetic; Kikuchi \cite{Makoto Kikuchi 94} uses provability in arithmetic; Chaitin \cite{Chaitin 74} and Kikuchi \cite{Kikuchi 97} use Kolmogorov complexity; Kikuchi and Tanaka \cite{Kikuchi-Tanaka 94} use nonstandard
models of arithmetic; and Vop\v{e}nka \cite{Vopenka 66} uses models of $\mathbf{ZF}$ (see \cite{Kikuchi-Kurahashi-Sakai}).

\subsubsection{Unexpected Examination and Grelling-Nelson's Paradox.}
First of all, Kritchman and Raz \cite{Kritchman-Raz 10} give a new proof  of ${\sf G2}$ based on Chaitin's incompleteness theorem and an argument that resembles the Unexpected Examination Paradox\footnote{The Unexpected Examination Paradox is formulated as follows in \cite{Kritchman-Raz 10}. The teacher announces in class: ``next week
you are going to have an exam, but you will not be able to know on which day of the
week the exam is held until that day". The exam cannot be held on Friday, because
otherwise, the night before the students will know that the exam is going to be held the
next day. Hence, in the same way, the exam cannot be held on Thursday. In the same
way, the exam cannot be held on any of the days of the week.} (for more details, we refer to \cite{Kritchman-Raz 10}): for any consistent r.e.~ extension $T$ of $\mathbf{PA}$, if $T$ is consistent, then $T\nvdash \mathbf{Con}(T)$.

\smallskip

Secondly, we say  a one-place
predicate is ``heterological" if it does not apply to itself (e.g.
``long" is heterological, since
it's not a long expression). Consider the question: is the predicate ``heterological" we have
just defined heterological? If ``heterological" is heterological, then it isn't heterological; and if ``heterological" isn't heterological, then it is heterological. This contradiction is called Grelling-Nelson's Paradox.

\smallskip

Cie\'{s}li\'{n}ski \cite{Heterologicality and incompleteness} presents semantic proofs of ${\sf G2}$ for $\mathbf{ZF}$ and $\mathbf{PA}$ based on Grelling-Nelson's Paradox. For a theory $T$ containing $\mathbf{ZF}$, Cie\'{s}li\'{n}ski
defines the sentence $\mathbf{HET_T}$ which says intuitively that the predicate ``heterological" is itself heterological, and then shows that $T\nvdash\mathbf{HET_T}$  and $T\vdash \mathbf{HET_T}\leftrightarrow \mathbf{Con}(T)$. Finally, Cie\'{s}li\'{n}ski shows how to adapt the proof of ${\sf G2}$ for $\mathbf{ZF}$ to a proof of ${\sf G2}$ for $\mathbf{PA}$. In fact, Cie\'{s}li\'{n}ski \cite{Heterologicality and incompleteness} proves the semantic version of ${\sf G2}$: if $T$ has a model, then $T + \neg \mathbf{Con}(T)$ has a model (i.e. $T\nvDash \mathbf{Con}(T)$).

\subsubsection{Yablo's paradox}
We discuss proofs of ${\sf G1}$ and ${\sf G2}$ based on Yablo's Paradox in the literature.
Yablo's Paradox is an infinite version of the Liar Paradox proposed in \cite{Yablo 93}: consider an infinite sequence $Y_1, Y_2, \cdots$ of propositions such that each $Y_i$ asserts that $Y_j$ are false for all $j > i$.
Different proofs of  ${\sf G1}$ and ${\sf G2}$ based on Yablo's Paradox have been given by some authors (e.g.~ Priest \cite{Priest 97}, Cie\'{s}li\'{n}ski-Urbaniak \cite{Godelizing the Yablo sequence}  and Kikuchi-Kurahashi \cite{Kikuchi-Kurahashi 11}).
\smallskip

Recall that we assume by default that $T$ is a  consistent r.e.~ extension of  $\mathbf{Q}$. Priest \cite{Priest 97} first points out that ${\sf G1}$ can be proved by formalizing Yablo's Paradox.
 Priest
defines a formula $Y(x)$ as follows which says that
for any $y > x, Y (y)$ is not provable in $T$.\smallskip

\begin{definition}[\cite{Godelizing the Yablo sequence, Kurahashi 14}]~
A formula $Y(x)$  is called a Yablo formula of $T$ if
$T \vdash \forall x (Y(x)\leftrightarrow \forall y > x \neg \mathbf{Pr}_T(\ulcorner Y(\dot{y})\urcorner))$.
\end{definition}

Cie\'{s}li\'{n}ski and Urbaniak originally prove   the following version of ${\sf G1}$, and show that each instance $Y(\overline{n})$ of the Yablo formula is independent of $T$ if $T$ is
$\Sigma^0_1$-sound (or 1-consistent).\smallskip

\begin{theorem}[Theorem 19, \cite{Godelizing the Yablo sequence}; see also Theorem 4, \cite{Kurahashi 14}]~
Let $Y(x)$  be a Yablo formula.
\begin{itemize}
  \item  If $T$ is consistent, then $T\nvdash Y(\overline{n})$.
  \item If $T$ is $\Sigma^0_1$-sound, then $T \nvdash \neg Y(\overline{n})$.
\end{itemize}
\end{theorem}

Cie\'{s}li\'{n}ski and Urbaniak originally prove that $T \vdash \forall x(Y(x)\leftrightarrow \mathbf{Con}(T))$ (see  \cite[Theorem 21-22]{Godelizing the Yablo sequence}). As a corollary, we have $T \vdash \forall x\forall y(Y(x)\leftrightarrow Y(y))$, and  ${\sf G2}$ holds: if $T$ is consistent, then $T\nvdash \mathbf{Con}(T)$.
%It is unknown whether we can derive ${\sf G1}$ based on Yablo's Paradox only assuming that $T$ is consistent.
\smallskip
\begin{definition}[\cite{Kurahashi 14, Godelizing the Yablo sequence}]
A formula $Y^R(x)$  is called a Rosser-type Yablo formula of $\mathbf{Prf}_T(x, y)$ if
$\mathbf{PA} \vdash \forall x (Y^R(x) \leftrightarrow \forall y > x\neg \mathbf{Pr}_T^R(x)(\ulcorner Y^R(\dot{y})\urcorner))$.
\end{definition}

Theorem \ref{Rosser-type Yablo} shows that the Rosser-type Yablo formula is independent of any
$\Sigma^0_1$-sound theory $T$.\smallskip

\begin{theorem}[Theorem 10, \cite{Kurahashi 14}]~\label{Rosser-type Yablo}
Let $\mathbf{Prf}_T(x, y)$ be any standard proof predicate of $T$, and $Y^R(x)$ be any
Rosser-type Yablo formula of $\mathbf{Prf}_T(x, y)$. Given $n \in \omega$,
if $T$ is consistent, then $T \nvdash Y^R(\overline{n})$;
if $T$ is $\Sigma^0_1$-sound, then $T \nvdash \neg Y^R(\overline{n})$.
\end{theorem}

The independence of $Y^R(\overline{n})$ for $T$ which is not $\Sigma^0_1$-sound is
discussed in \cite{Kurahashi 14}. For a consistent but not $\Sigma^0_1$-sound theory,  the situation of  Rosser-type Yablo formulas is quite different from that of Rosser sentences. Kurahashi \cite{Kurahashi 14} shows that for any consistent but not $\Sigma^0_1$-sound theory, the independence
of each instance of a Rosser-type Yablo formula depends on the choice of standard proof predicates (see Theorem 12 and Theorem 25 in \cite{Kurahashi 14}). Kurahashi \cite{Kurahashi 14} shows
that for any consistent but not $\Sigma^0_1$-sound theory $T$, there is a standard proof predicate of
$T$ such that each instance $Y^R(\overline{n})$ of the Rosser-type Yablo formula $Y^R(x)$ based on this
proof predicate is provable in $T$ for any $n \in\omega$. Moreover, Kurahashi \cite{Kurahashi 14}  constructs
a standard proof predicate of $T$ and a Rosser-type Yablo formula $Y^R(x)$ based on this
proof predicate such that each instance of $Y^R(x)$ is independent of $T$.
Proofs of these results use the technique of Guaspari and Solovay in \cite{Guaspari 79}.

\smallskip

Cie\'{s}li\'{n}ski and Urbaniak \cite{Godelizing the Yablo sequence} conjecture that any  two distinct
instances $Y^R(\overline{m})$ and $Y^R(\overline{n})$ of a Rosser-type Yablo formula $Y^R(x)$ based on a
standard proof predicate
are not provably
equivalent.  Leach-Krouse \cite{Leach-Krouse 13} and Kurahashi \cite{Kurahashi 14} construct
a standard proof predicate,  and a Rosser-type Yablo formula $Y^R(x)$ based on this proof
predicate such that $T\vdash\forall x\forall y (Y^R(x) \leftrightarrow Y^R(y))$ (see  \cite[Theorem 9]{Leach-Krouse 13} and \cite[Corollary 21]{Kurahashi 14}).

\smallskip

Kurahashi \cite{Kurahashi 14} constructs a partial counterexample to Cie\'{s}li\'{n}ski and
Urbaniak's conjecture: a standard proof predicate, and a Rosser-type Yablo formula $Y^R(x)$ based on this
proof predicate such that $\forall x\forall y (Y^R(x) \leftrightarrow Y^R(y))$ is not provable in $T$ (Corollary 20, \cite{Kurahashi 14}).
Thus the provability of the sentence $\forall x\forall y (Y^R(x) \leftrightarrow Y^R(y))$ also depends on the choice of  standard proof predicates (see Corollary 20-21, \cite{Kurahashi 14}). Proofs of these results by Leach-Krouse and Kurahashi also use the technique of Guaspari and Solovay in \cite{Guaspari 79}.
An interesting open question  is: whether there is a standard proof predicate such that $Y^R(\overline{n})$ and $Y^R(\overline{n+1})$ are not provably equivalent for some $n \in\omega$ (see \cite{Kurahashi 14}).

\subsubsection{Beyond arithmetization}
All the proofs of ${\sf G1}$ we have discussed  use arithmetization.
Andrzej Grzegorczyk proposes the theory $\mathbf{TC}$ in \cite{Grzegorczyk 05} as a possible alternative theory for studying incompleteness
and undecidability, and shows that  $\mathbf{TC}$ is essentially incomplete and mutually interpretable with $\mathbf{Q}$ without arithmetization.

\smallskip

Now, in $\mathbf{PA}$ we have numbers that can be added or multiplied; while in  $\mathbf{TC}$, one has strings (or texts) that can be concatenated. In G\"{o}del's proof, the only use of numbers is coding of syntactical objects.
The motivations for accepting strings rather than numbers as the
basic notion are as follows: on metamathematical level, the notion of computability can be
defined without reference to numbers; for Grzegorczyk, dealing with texts is philosophically better justified since  intellectual activities like reasoning, communicating or even computing involve working with texts not with numbers (see \cite{Grzegorczyk 05}). Thus, it is natural to
define notions like undecidability directly in terms of texts instead of natural numbers.
Grzegorczyk only proves the incompleteness of $\mathbf{TC}$ in \cite{Grzegorczyk 05}. Later, Grzegorczyk and Zdanowski \cite{Grzegorczyk 08} prove that $\mathbf{TC}$ is essentially incomplete.

\subsection{Concrete incompleteness}\label{concrete incompleteness}
\subsubsection{Introduction}
All proofs of G\"{o}del's incompleteness theorems we have discussed above
make use of meta-mathematical or logical methods, and the independent sentence constructed has a clear meta-mathematical or logical flavour which is devoid of real mathematical content.
To be blunt, from a purely mathematical point of view, G\"{o}del's sentence is artificial and not mathematically interesting. G\"{o}del's sentence is constructed not by reflecting about arithmetical properties of natural numbers, but by
reflecting about an axiomatic system in which those properties are
formalized (see \cite{Isaacson 87}).
A natural question is then: can we find true sentences not provable in $\mathbf{PA}$ with real mathematical content?
The research program \emph{concrete incompleteness} is the search for natural independent sentences with real mathematical content.

\smallskip

This program has received a lot of attention because despite G\"{o}del's incompleteness theorems, one can still cherish the hope that all natural and mathematically interesting sentences about natural numbers are provable or refutable  in $\mathbf{PA}$, and that elementary arithmetic is complete w.r.t.~ natural and mathematically interesting sentences.
However,  after G\"{o}del, many natural independent sentences with real mathematical content have been found. These independent sentences have a clear mathematical flavor, and do not refer to the arithmetization of syntax and provability.

\smallskip

In this section, we provide an overview of the research on concrete incompleteness. The survey paper \cite{Bovykin} provides a good overview on the state-of-the-art up to Autumn 2006.
For more detailed discussions about concrete incompleteness, we refer to Cheng \cite{Cheng book 19}, Bovykin \cite{Bovykin} and Friedman \cite{Friedman 18}.

\subsubsection{Paris-Harrington and beyond}
%First of all, Gentzen \cite{Gentzen 43} shows that the mathematically natural statement expressing $\epsilon_0$-transfinite induction is independent of $\mathbf{PA}$.
Paris and Harrington \cite{Paris and Harrington} propose the first mathematically natural statement independent of $\mathbf{PA}$:  the \emph{Paris-Harrington Principle} $\textsf{PH}$ which generalizes the \emph{finite Ramsey theorem}.
G\"{o}del's sentence is a pure logical construction (via the arithmetization of syntax and provability predicate) and has no relevance with classic mathematics (without any combinatorial or number-theoretic content). On the contrary, Paris-Harrington Principle is an independent arithmetic sentence from classic mathematics with combinatorial content as we will show. We refer to \cite{Isaacson 87} for more discussions about the distinction between mathematical arithmetic sentences and meta-mathematical arithmetic sentences.

\smallskip

\begin{definition}[Paris-Harrington Principle (${\sf PH}$), \cite{Paris and Harrington}]~
\begin{itemize}
  \item For set $X$ and $n\in\omega$, let $[X]^n$ be the set of all $n$-elements subset of $X$.  We identify  $n$ with $\{0, \cdots, n-1\}$.
  \item For all  $m, n, c\in \omega$, there is $N\in \omega$ such that for all  $f : [N]^m\rightarrow c$, we have: $(\exists H \subseteq N)(|H| \geq n \wedge  \textup{$H$ is homogeneous for $f$}  \wedge |H| > \min(H))$.
\end{itemize}
\end{definition}\smallskip
%One can then prove the following ground-breaking theorem.
\begin{theorem}[Paris-Harrington, \cite{Paris and Harrington}]~
The principle $\sf PH$ is true but not provable in $\mathbf{PA}$.
\end{theorem}
Now, $\sf PH$ has a clear combinatorial flavor, and is of the form $\forall x \exists y\psi(x, y)$ where  $\psi$ is a $\Delta^0_0$ formula.
It can be shown that for any given natural number $n$, $\mathbf{PA} \vdash\exists y\psi(\overline{n}, y)$, i.e.\ all particular \emph{instances} of $\sf PH$ are provable in $\mathbf{PA}$.

\smallskip

Following $\sf PH$, many other mathematically natural statements independent of $\mathbf{PA}$ with combinatorial or number-theoretic content have been formulated: the Kanamori-McAloon principle \cite{Kanamori and McAloon}, the Kirby-Paris sentence \cite{Kirby and Paris}, the
Hercules-Hydra game \cite{Kirby and Paris}, the Worm principle \cite{Beklemishev 2003, Hamano and Okada},
the flipping principle \cite{Kirby 82}, the arboreal statement \cite{Mills 80},  the kiralic and regal principles \cite{Clote 83}, and  the Pudl\'{a}k's principle \cite{Pudlak 79, Hajek 86}  (see \cite{Bovykin}, p.40). In fact, all these principles are equivalent to $\sf PH$ (see \cite{Bovykin}, p.40).

\smallskip

An interesting and amazing fact is:  all the above mathematically natural principles are in fact provably equivalent in $\mathbf{PA}$ to a certain meta-mathematical sentence.
Consider the following reflection principle for $\Sigma^0_1$ sentences:
for any $\Sigma^0_1$ sentence $\phi$ in $L(\mathbf{PA})$, if $\phi$ is provable in $\mathbf{PA}$, then $\phi$ is true.
Using the arithmetization of syntax, one can write this principle
as a sentence of $L(\PA)$, and denote it by $\sf Rfn_{\Sigma^0_1}(\mathbf{PA})$ (see  \cite[p.301]{metamathematics}). McAloon has shown that\, $\mathbf{PA} \vdash \sf PH\leftrightarrow {\sf Rfn}_{\Sigma^0_1}(\mathbf{PA})$ (see \cite{metamathematics}, p.301), and similar equivalences can be established for the  other independent principles mentioned above. Equivalently, all these principles are equivalent to so-called 1-consistency of $\mathbf{PA}$ (see \cite[p.36]{Beklemishev 45}, \cite[p.3]{Beklemishev 2003} and \cite[p.301]{metamathematics}).

\smallskip

The above phenomenon indicates that the difference between mathematical and meta-mathematical statements is perhaps not as huge as we might have expected.
Moreover, the above principles are provable in fragments of Second-Order Arithmetic and are more complex than G\"{o}del's sentence: G\"{o}del's sentence is equivalent  to $\mathbf{Con}(\mathbf{PA})$ in $\mathbf{PA}$; but all these principles
are not only independent of $\mathbf{PA}$ but also independent of $\mathbf{PA}+ \mathbf{Con}(\mathbf{PA})$ (see \cite[p.36]{Beklemishev 45}  and \cite[p.301]{metamathematics}).

\subsubsection{Harvey Friedman's contributions}
Incompleteness would not be complete without mentioning the work of Harvey Friedman who is a central figure in  research on the foundations of mathematics after G\"{o}del.  He has made many important contributions to concrete mathematical
incompleteness. The following quote is telltale:
\begin{quote}
the long range impact
and significance of ongoing investigations in the
foundations of mathematics is going to depend greatly on
the extent to which the Incompleteness Phenomena touches
normal concrete mathematics (see \cite{Friedman 18}, p.7).
\end{quote}
In the following, we give a brief introduction to H. Friedman's work on concrete mathematical incompleteness.
In his early work, H. Friedman examines how one uses large cardinals in an essential and
natural way in number theory, as follows.
\begin{quote}
the quest for a simple meaningful finite mathematical theorem that can
only be proved by going beyond the usual axioms for mathematics has been a
goal in the foundations of mathematics since G\"{o}del's incompleteness theorems (see \cite{Friedman 98}, p.805).
\end{quote}
H. Friedman shows in \cite{Friedman 81, Friedman 98}  that there are many mathematically
natural combinatorial statements in $L(\mathbf{PA})$ that are neither provable nor refutable
in  $\mathbf{ZFC}$  or $\mathbf{ZFC}$  + large cardinals.
H. Friedman's more recent monograph \cite{Friedman 18} is a comprehensive study of concrete
mathematical incompleteness. H. Friedman studies concrete mathematical incompleteness over different systems, ranging from weak subsystems of $\mathbf{PA}$ to higher-order arithmetic and $\mathbf{ZFC}$.
H. Friedman lists many concrete mathematical statements in $L(\mathbf{PA})$ that are independent of subsystems of $\mathbf{PA}$, or stronger
theories like higher-order arithmetic and set theory.

\smallskip

The theories $\mathbf{RCA_0}$ (Recursive Comprehension), $\mathbf{WKL_0}$  (Weak Konig's Lemma), $\mathbf{ACA_0}$ (Arithmetical Comprehension), $\mathbf{ATR_0}$ (Arithmetic
Transfinite Recursion) and {$\Pi^1_1$-$\mathbf{CA_0}$} ($\Pi^1_1$-Comprehension) are the most famous five subsystems of Second-Order Arithmetic ($\mathbf{SOA}$), and are called the `Big Five'.
For the definition of $\mathbf{SOA}$ and the `Big Five', we refer to \cite{SOA}.

\smallskip

To give the reader a better sense of H. Friedman's work, we list some sections dealing with concrete mathematical
incompleteness in \cite{Friedman 18}.
\begin{itemize}
\item Section 0.5 on Incompleteness in Exponential Function Arithmetic.
\item Section 0.6 on Incompleteness in Primitive Recursive Arithmetic, Single Quantifier Arithmetic, $\mathbf{RCA_0}$, and $\mathbf{WKL_0}$.
\item Section 0.7 on Incompleteness in Nested Multiply Recursive Arithmetic and Two Quantifier Arithmetic.
\item Section 0.8 on Incompleteness in Peano Arithmetic and $\mathbf{ACA_0}$.
\item Section 0.9 on Incompleteness in Predicative Analysis and $\mathbf{ATR_0}$.
\item Section 0.10 on Incompleteness in Iterated Inductive Definitions and $\Pi^1_1$-$\mathbf{CA_0}$.
\item Section 0.11 on Incompleteness in Second-Order Arithmetic and $\mathbf{ZFC}^{-}$.\footnote{$\mathbf{ZFC}^{-}$ denotes $\mathbf{ZFC}$  with the Power Set Axiom deleted and Collection instead of Replacement.}
\item Section 0.12 on Incompleteness in Russell Type Theory and Zermelo Set Theory.
\item Section 0.13 on Incompleteness in $\mathbf{ZFC}$ using Borel Functions.
\item Section 0.14 on Incompleteness in $\mathbf{ZFC}$ using Discrete Structures.
\end{itemize}

H.  Friedman \cite{Friedman 18} provides us with examples of concrete mathematical theorems not provable in subsystems of Second-Order Arithmetic stronger than $\mathbf{PA}$, and a number of concrete mathematical statements provable in Third-Order Arithmetic but not provable in Second-Order Arithmetic.

\smallskip

Related to Friedman's work, Cheng \cite{Cheng book 19, Cheng and Schindler 2015}  gives an example of concrete mathematical theorems based on Harrington's principle which is isolated from the proof of the Harrington's Theorem (the determinacy of $\Sigma^1_1$ games implies the existence of zero sharp), and shows that this concrete theorem saying that Harrington's principle implies the existence of zero sharp is  expressible in Second-Order Arithmetic, not provable in Second-Order Arithmetic or Third-Order Arithmetic, but provable in Fourth-Order Arithmetic (i.e.~ the minimal system in higher-order arithmetic to prove this concrete theorem is Fourth-Order Arithmetic).

\smallskip

Many other examples of concrete mathematical incompleteness,  and the discussion of this  subject in  1970s-1980s  can be found in the four volumes \cite{Simpson 87, Simpson 85, Pacholski 80, Berline 81}. Weiermann's work in \cite{Weiermann 03}-\cite{Weiermann 12} provides us with more  examples  of naturally mathematical independent sentences.
We refer to \cite{Friedman 18} for new advances in Boolean Relation Theory and for more examples of concrete mathematical
incompleteness.

\section{The limit of the applicability of  ${\sf G1}$}\label{Generalizations}
%In this section, I examine generalizations of  G\"{o}del's incompleteness theorems from different perspectives. In the following, I first examine the generalizations of ${\sf G1}$ and then the generalizations of ${\sf G2}$.
%\subsection{Generalizations of ${\sf G1}$}\label{Generalizations of G1}
\subsection{Introduction}
In this section, we discuss the limit of the applicability of ${\sf G1}$ based on the following two questions.
\begin{itemize}
\item To what extent does ${\sf G1}$ apply to  extensions of $\mathbf{PA}$?
\item To what extent does ${\sf G1}$ apply to theories weaker  than $\mathbf{PA}$ w.r.t.~ interpretation?
\end{itemize}

\begin{definition}[Conservativity]~\label{}
\begin{itemize}
  \item Let $\Gamma$ denote either $\Sigma^0_n$ or $\Pi^0_n$ for some $n \geq 1$, and $\Gamma^d$ denote either $\Pi^0_n$ or $\Sigma^0_n$.
  \item We say a sentence $\varphi$ is $\Gamma$-conservative over theory $T$ if for any
$\Gamma$ sentence $\psi$, $T \vdash\psi$ whenever $T + \varphi\vdash \psi$.
\end{itemize}
\end{definition}

We list some generalizations of ${\sf G1}$ needed below.\smallskip
\begin{fact}[Guaspari \cite{Guaspari 1979}]\label{Guaspari's fact}
Let $T$ be a consistent r.e.~ extension of $\mathbf{Q}$. Then there is a  $\Gamma^d$ sentence $\phi$ such that $\phi$ is $\Gamma$-conservative over $T$ and $T\nvdash\phi$.
\end{fact}
If $T\vdash\neg\phi$, then $\phi$ is not $\Gamma$-conservative over
$T$ because $T$ is consistent. Thus, we can view Fact \ref{Guaspari's fact}  as an extension of Rosser's first incompleteness theorem. Solovay improves this fact and shows that there is a $\Gamma^d$ sentence $\phi$ such that $\phi$ is $\Gamma$-conservative over $T$, $\neg\phi$ is $\Gamma^d$-conservative over $T$, but $\phi$ is independent of $T$.
\smallskip
\begin{fact}[Mostowski \cite{Mostowski 1961}]\label{Mostowski's fact}
Let $\{T_n: n\in\omega\}$ be an r.e.~ sequence of consistent theories extending $\mathbf{Q}$. Then
there is a $\Pi^0_1$ sentence $\phi$ such that for any $n\in\omega$, $T_n\nvdash\phi$, and $T_n\nvdash\neg\phi$.
\end{fact}

\subsection{Generalizations of ${\sf G1}$ beyond $\mathbf{PA}$}

We study generalization of ${\sf G1}$ for extensions of $\mathbf{PA}$ w.r.t.~ interpretation.
We know that ${\sf G1}$ applies to all consistent r.e.~ extensions of $\mathbf{PA}$.
A natural question is then: whether ${\sf G1}$ can be extended to non-r.e.~ arithmetically definable extensions of $\mathbf{PA}$.

\smallskip

Kikuchi-Kurahashi \cite{Kikuchi 17} and  Salehi-Seraji \cite{Saeed Salehi 2017} make contributions to  generalize   G\"{o}del-Rosser's first incompleteness
theorem  to non-r.e.~ arithmetically definable extensions of $\mathbf{PA}$.\smallskip
%Proofs in \cite{Kikuchi 17} heavily use techniques from \cite{Aspects of incompleteness}.

\begin{definition}[\cite{Kikuchi 17}]
Let $T$ be a consistent  extension of $\mathbf{Q}$.
\begin{itemize}
  \item $T$ is $\Sigma^0_n$-definable if there is a  $\Sigma^0_n$ formula $\phi(x)$ such that $n$ is the G\"{o}del number of some sentence of $T$ if and only if $\mathfrak{N}\models \phi(\overline{n})$.\footnote{Recall that $\mathfrak{N}$ is the standard model of arithmetic.}
  \item $T$ is $\Sigma^0_n$-sound if for all $\Sigma^0_n$ sentences $\phi$,  $T\vdash\phi$ implies $\mathfrak{N}\models\phi$; $T$ is sound if  $T$ is $\Sigma^0_{n}$-sound for any $n \in\omega$.
%\item $T$ is $\Gamma$-consistent if and only if for all $\Gamma$ formulas $\phi(x)$, if $T \vdash\neg\phi(\overline{n})$ for all $n\in\omega$, then $T \nvdash \exists x\phi(x)$.
    \item $T$ is $\Sigma^0_n$-consistent if  for all $\Sigma^0_{n}$ formulas $\phi$ with $\phi=\exists x\theta(x)$ and $\theta\in\Pi^0_{n-1}$, if $T \vdash\neg\theta(\overline{n})$ for all $n\in\omega$, then $T \nvdash \phi$.
 \item
 $T$ is $\Pi^0_n$-decisive if  for all $\Pi^0_n$ sentences $\phi$, either $T\vdash\phi$ or $T \vdash\neg\phi$ holds.
\end{itemize}
\end{definition}

From ${\sf G1}$, we have: if $T$ is a $\Sigma^0_{1}$-definable and $\Sigma^0_{1}$-sound extension of $\mathbf{Q}$, then $T$ is not $\Pi^0_{1}$-decisive. Kikuchi and Kurahashi \cite{Kikuchi 17}   generalize ${\sf G1}$ to arithmetically definable theories via the notion of ``$\Sigma^0_{n}$-sound".
\smallskip

\begin{theorem}[Theorem 4.8 \cite{Kikuchi 17}, Theorem 2.5 \cite{Saeed Salehi 2017}]\label{generalization of G}~
If $T$ is a $\Sigma^0_{n+1}$-definable and $\Sigma^0_{n}$-sound extension of $\mathbf{Q}$, then $T$ is not $\Pi^0_{n+1}$-decisive.
\end{theorem}

Salehi and Seraji \cite{Saeed Salehi 2017} point out that  Theorem \ref{generalization of G} has a constructive proof:   given a $\Sigma^0_{n+1}$-definable and $\Sigma^0_{n}$-sound extension $T$ of $\mathbf{Q}$,  one can effectively construct a $\Pi^0_{n+1}$ sentence which is independent of $T$.
The optimality of Theorem \ref{generalization of G} is shown by Salehi and Seraji in \cite{Saeed Salehi 2017}: there exists a $\Sigma^0_{n-1}$-sound and $\Sigma^0_{n+1}$-definable complete extension of $\mathbf{Q}$ for any $n\geq 1$ (Theorem 2.6, \cite{Saeed Salehi 2017}).

\smallskip

Salehi and Seraji \cite{Saeed Salehi 2017}  generalize ${\sf G1}$ to arithmetically definable theories via the notion of ``$\Sigma^0_n$-consistent".
\begin{theorem}[Theorem 4.9 \cite{Kikuchi 17}, Theorem 4.3 \cite{Saeed Salehi 2017}]\label{generalization of G second}~
If $T$ is a $\Sigma^0_{n+1}$-definable and $\Sigma^0_n$-consistent extension of $\mathbf{Q}$, then $T$ is not $\Pi^0_{n+1}$-decisive.
\end{theorem}

Theorem \ref{generalization of G second} is also optimal: the complete $\Sigma^0_{n-1}$-sound and $\Sigma^0_{n+1}$-definable theory constructed in the
proof of Theorem 2.6 in \cite{Saeed Salehi 2017} is also $\Sigma^0_{n-1}$-consistent since if a theory is $\Sigma^0_{n}$-sound, then it is $\Sigma^0_{n}$-consistent.
The proof of Theorem~\ref{generalization of G second} cannot be
constructive as the following theorem shows.\smallskip

\begin{theorem}[Non-constructivity of $\Sigma^0_{n}$-consistency incompleteness, Theorem 4.4, \cite{Saeed Salehi 2017}]~
For $n\geq 3$, there is no (partial) recursive function $f$ (even with the oracle $0^{n}$) such that if $m$ codes (the G\"{o}del code) a $\Sigma^0_{n+1}$-formula which defines an $\Sigma^0_{n}$-consistent extension $T$ of $\mathbf{Q}$, then $f(m)$ halts and codes a  $\Pi^0_{n+1}$ sentence which is independent of $T$.\footnote{Salehi and Seraji \cite{Saeed Salehi 2017} remark that there indeed exists some $0^{n+1}$-(total) recursive function $f$  such that if $m$ codes  a $\Sigma^0_{n+1}$-formula defining an $\Sigma^0_{n}$-consistent extension $T$ of $\mathbf{Q}$, then $f(m)$ halts and codes a  $\Pi^0_{n+1}$ sentence  independent of $T$.}
\end{theorem}

In  summary, ${\sf G1}$ can be generalized to the incompleteness of $\Sigma^0_{n+1}$-definable
and $\Sigma^0_{n}$-sound extensions of $\mathbf{Q}$ constructively;  and  to the incompleteness of $\Sigma^0_{n+1}$-definable and $\Sigma^0_{n}$-consistent extensions of $\mathbf{Q}$ non-constructively (when $n>2$).

\subsection{Generalizations of ${\sf G1}$ below $\mathbf{PA}$}
We study generalizations of ${\sf G1}$ for theories weaker than $\mathbf{PA}$ w.r.t.~ interpretation.

\subsubsection{Generalizations of ${\sf G1}$ via interpretability}

We show that $\sf G1$ can be generalized to theories weaker than $\mathbf{PA}$ via interpretability.  Indeed, there exists a weak recursively axiomatizable consistent subtheory $T$ of $\mathbf{PA}$ such that each recursively axiomatizable theory $S$ in which $T$ is interpretable is incomplete (see \cite{undecidable}). To generalize this fact further, we propose a new notion ``$\sf G1$ holds for  $T$", as follows. %which proves us a route to find the limit of incompleteness.
\smallskip

\begin{definition}
Let $T$ be a consistent r.e.~  theory. We say \emph{$\sf G1$ holds for  $T$}  if for any recursively axiomatizable consistent theory $S$, if $T$ is interpretable in $S$, then  $S$ is incomplete.
\end{definition}
First of all, for a consistent r.e.~  theory $T$, it is not hard to show that the followings are equivalent (see \cite{Cheng 19}):
\begin{itemize}
  \item  $\sf G1$ holds for  $T$.
  \item $T$ is essentially incomplete.
  \item $T$ is essentially undecidable.
\end{itemize}
It is well-known that $\sf G1$ holds for many weaker theories than $\mathbf{PA}$ w.r.t.~  interpretation (e.g. Robinson arithmetic $\mathbf{Q}$).

\smallskip

Secondly, we mention theories weaker than $\mathbf{PA}$ w.r.t.~  interpretation for which $\sf G1$ holds.
We first review some  essentially undecidable theories weaker than $\mathbf{PA}$ w.r.t.~ interpretation  from the literature (i.e.~ $\sf G1$ holds for these theories).  For the definition of theory $\mathbf{Q}$, $I\Sigma_n$, $B\Sigma_n$,
$\mathbf{PA}^{-}$, $\mathbf{Q}^+$, $\mathbf{Q}^{-}$,  $\mathbf{S^1_2}$,  $\mathbf{AS}$, $\mathbf{EA}$, $\mathbf{PRA}$, $\mathbf{R}$, $\mathbf{R}_0$, $\mathbf{R}_1$ and $\mathbf{R}_2$, we refer to Section \ref{Preliminaries}.

\smallskip

Robinson shows that any consistent r.e.~ theory that interprets $\mathbf{Q}$ is undecidable, and hence
$\mathbf{Q}$ is essentially undecidable. The fact that $\mathbf{Q}$ is  essentially undecidable is very useful and can be used to prove the  essentially undecidability of other theories via Theorem \ref{interpretable theorem}. Since $\mathbf{Q}$ is finitely axiomatized, it follows that any
theory that weakly interprets $\mathbf{Q}$ is also undecidable.

\smallskip

The Lindenbaum algebras of all r.e.~ theories that interpret $\mathbf{Q}$ are recursively isomorphic (see Pour-El and Kripke \cite{Kripke 67}).
In fact, $\mathbf{Q}$ is minimal essentially undecidable in the sense that if deleting any axiom of $\mathbf{Q}$, then the remaining theory is not essentially undecidable and has a complete decidable extension (see \cite[Theorem 11, p.62]{undecidable}).

\smallskip

Thirdly, Nelson  \cite{Nelson 86} embarks  on a program   of investigating how much mathematics can
be interpreted in Robinson's Arithmetic $\mathbf{Q}$: what can be interpreted in $\mathbf{Q}$,  and what cannot  be  interpreted in $\mathbf{Q}$. In fact, $\mathbf{Q}$ represents a rich degree of interpretability since a lot of stronger theories are interpretable in it as we will show in the following passages. For example, using Solovay's
method of shortening cuts (see \cite{Guaspari 79}), one can show that $\mathbf{Q}$ interprets fairly strong theories like $I\Delta_0 + \Omega_1$ on a definable cut.

\smallskip

Fourth, we  discuss some prominent fragments of $\mathbf{PA}$ extending $\mathbf{Q}$ from the literature.
As a corollary of Theorem \ref{interpretation theorem}, we have:
\begin{itemize}
  \item The theories $\mathbf{Q},  I\Sigma_0, I\Sigma_0+\Omega_{1}, \cdots, I\Sigma_0+\Omega_{n}, \cdots, B\Sigma_1, B\Sigma_1+\Omega_{1}, \cdots$, $B\Sigma_1+\Omega_{n}, \cdots$ are all mutually interpretable;
  \item $I\Sigma_0+\mathbf{exp}$ and $B\Sigma_1+\mathbf{exp}$ are  mutually interpretable;
   \item For $n\geq 1$, $I\Sigma_n$ and $B\Sigma_{n+1}$ are  mutually interpretable;
    \item  $\mathbf{Q}\lhd I\Sigma_0+\mathbf{exp}\lhd I\Sigma_1\lhd I\Sigma_2\lhd\cdots\lhd I\Sigma_n\lhd\cdots\lhd \mathbf{PA}$.
\end{itemize}
Since any consistent r.e.~ theory which interprets $\mathbf{Q}$ is essentially undecidable, $\sf G1$ holds for all these fragments of $\mathbf{PA}$ extending $\mathbf{Q}$.

\smallskip

Fifth, we discuss some  weak theories  mutually interpretable with  $\mathbf{Q}$ from the literature.
It is interesting to compare $\mathbf{Q}$ with its bigger brother $\mathbf{PA}^{-}$.
From \cite{paper on Q}, $\mathbf{PA}^{-}$ is interpretable in $\mathbf{Q}$, and hence $\mathbf{Q}$ is mutually interpretable with $\mathbf{PA}^{-}$.
The theory $\mathbf{Q}^+$ is interpretable in $\mathbf{Q}$ (see Theorem 1 in \cite{Interpretability in Robinson's Q}, p.296), and thus mutually interpretable with $\mathbf{Q}$.
A. Grzegorczyk asks whether $\mathbf{Q}^{-}$ is essentially undecidable.
\v{S}vejdar \cite{Svejdar 07} provids a positive answer to Grzegorczyk's original question by showing that $\mathbf{Q}$ is interpretable in $\mathbf{Q}^{-}$ using the Solovay's method of shortening  cuts. Thus $\mathbf{Q}^{-}$ is essentially undecidable and mutually interpretable with $\mathbf{Q}$.

\smallskip

Sixth, by \cite{Interpretability in Robinson's Q},  $I\Sigma_0$ is interpretable in  $\mathbf{S^1_2}$,  and $\mathbf{S^1_2}$ is interpretable in $\mathbf{Q}$. Hence $\mathbf{S^1_2}$ is essentially undecidable and mutually interpretable with $\mathbf{Q}$.
The theory $\mathbf{AS}$ interprets Robinson's Arithmetic $\mathbf{Q}$, and hence is essentially undecidable.
Nelson \cite{Nelson 86} shows that $\mathbf{AS}$ is interpretable in  $\mathbf{Q}$. Thus,  $\mathbf{AS}$ is mutually interpretable with $\mathbf{Q}$.

\smallskip

Seventh, Grzegorczyk and Zdanowski \cite{Grzegorczyk 08} formulate but leave unanswered an interesting problem: are $\mathbf{TC}$ and $\mathbf{Q}$ mutually interpretable?
M. Ganea \cite{Ganea 09} provs that $\mathbf{Q}$ is interpretable in  $\mathbf{TC}$ using the detour via $\mathbf{Q}^{-}$ (i.e.~ first show that $\mathbf{Q}^{-}$ is interpretable in  $\mathbf{TC}$; since $\mathbf{Q}$ is interpretable in $\mathbf{Q}^{-}$, then we have $\mathbf{Q}$ is interpretable in  $\mathbf{TC}$). Sterken and Visser \cite{Visser 09} give a proof of the interpretability of $\mathbf{Q}$ in $\mathbf{TC}$ not using $\mathbf{Q}^{-}$.
Note that $\mathbf{TC}$ is easily interpretable in the bounded arithmetic $I\Sigma_0$. Thus,  $\mathbf{TC}$ is mutually interpretable with $\mathbf{Q}$.

\smallskip

Note that $\mathbf{R}\lhd\mathbf{Q}$ since $\mathbf{Q}$ is not interpretable in $\mathbf{R}$ (if $\mathbf{Q}$ is interpretable in $\mathbf{R}$, then $\mathbf{Q}$ is interpretable in some finite fragment of $\mathbf{R}$; however  $\mathbf{R}$ is locally finitely satisfiable and any model of $\mathbf{Q}$ is infinite).
Visser \cite{Visser 14} provides us with a unique characterization of $\mathbf{R}$.\smallskip

\begin{theorem}[Visser, Theorem 6, \cite{Visser 14}]~\label{visser thm on R}
For any consistent r.e.~ theory $T$,  $T$ is interpretable in $\mathbf{R}$  if and only if $T$ is locally finitely satisfiable.\footnote{In fact, if $T$ is locally finitely satisfiable, then $T$ is interpretable in $\mathbf{R}$ via a one-piece one-dimensional parameter-free interpretation.}
\end{theorem}

Since relational $\Sigma_2$ sentences have the finite model property, by Theorem~\ref{visser thm on R},  any consistent theory axiomatized by a recursive set of $\Sigma_2$ sentences in a finite relational language is interpretable in $\mathbf{R}$. Since all recursive functions are representable in $\mathbf{R}$ (see \cite[theorem 6]{undecidable}, p.56), as a corollary of  Theorem \ref{interpretable theorem},  $\mathbf{R}$ is  essentially undecidable. Cobham shows that $\mathbf{R}$ has a stronger property than essential undecidability. Vaught gives a proof of Cobham's Theorem \ref{Cobham theorem} via existential interpretation in \cite{Vaught 62}.\smallskip

\begin{theorem}[Cobham, \cite{Vaught 62}]~\label{Cobham theorem}
Any consistent r.e.~ theory that weakly interprets $\mathbf{R}$ is undecidable.
\end{theorem}

Eighth, we discuss some variants of $\mathbf{R}$   in the same language as $L(\mathbf{R})=\{\overline{0}, \cdots, \overline{n}, \cdots, +, \times, \leq\}$.
The theory $\mathbf{R}_0$ is no longer essentially undecidable in the same language as $\mathbf{R}$.\footnote{The theory $\mathbf{R}_0$ has a decidable complete extension given by the theory of reals with $\leq$ as the empty relation on reals.}  In fact, whether $\mathbf{R}_0$ is  essentially undecidable depends on the language of $\mathbf{R}_0$: if $L(\mathbf{R}_0)=\{\mathbf{0}, \mathbf{S}, + , \times, \leq\}$  with $\leq$  defined in terms
of $+$, then $\mathbf{R}_0$ is essentially undecidable (Cobham first observed that  $\mathbf{R}$ is interpretable in $\mathbf{R}_0$ in the same language $\{\mathbf{0}, \mathbf{S}, + , \times\}$, and hence $\mathbf{R}_0$ is essentially undecidable (see \cite{Vaught 62} and \cite{Jones 83})). The theory $\mathbf{R}_1$ is  essentially undecidable since $\mathbf{R}$ is interpretable in $\mathbf{R}_1$ (see \cite{Jones 83}, p.62).

\smallskip

However $\mathbf{R}_1$ is not minimal essentially undecidable.
From \cite{Jones 83}, $\mathbf{R}$ is interpretable in $\mathbf{R}_2$, and hence $\mathbf{R}_2$ is essentially undecidable.\footnote{Another way to show that $\mathbf{R}_2$ is essentially
undecidable is to prove that all recursive functions are representable in $\mathbf{R}_2$.} The theory $\mathbf{R}_2$  is minimal essentially undecidable in the sense that if we delete any axiom scheme of $\mathbf{R}_2$, then the remaining system is not essentially undecidable.\footnote{If we delete $\sf {Ax2}$, then the theory of natural numbers with $x\times y$ defined as $x+y$ is a complete  decidable extension; if we delete $\sf {Ax3}$, then the theory of models with only one element is a complete  decidable extension; if we delete $\sf {Ax4^{\prime}}$, then the theory of reals is a complete  decidable extension.}
By essentially the same argument as in \cite{paper on Q}, we can show that any consistent r.e.~ theory that weakly interprets $\mathbf{R}_2$ is undecidable.

\smallskip

Kojiro Higuchi and Yoshihiro Horihata introduce the theory of concatenation $\mathbf{WTC}^{-\epsilon}$, which is a weak subtheory
of Grzegorczyk's theory $\mathbf{TC}$, and show that $\mathbf{WTC}^{-\epsilon}$ is  minimal essentially undecidable and $\mathbf{WTC}^{-\epsilon}$ is mutually interpretable with $\mathbf{R}$ (see \cite{Higuchi}).

\smallskip

In summary,  we have the following pictures:
\begin{itemize}
  \item Theories $\mathbf{PA}^{-}, \mathbf{Q}^{+}, \mathbf{Q}^{-},  \mathbf{TC}, \mathbf{AS},\mathbf{S^1_2}$ and $\mathbf{Q}$  are all mutually interpretable, and hence ${\sf G1}$ holds for them;
      \item Theories $\mathbf{R}, \mathbf{R}_1, \mathbf{R}_2$ and $\mathbf{WTC}^{-\epsilon}$ are mutually interpretable, and hence ${\sf G1}$ holds for them;
\item  $\mathbf{R}\lhd\mathbf{Q}\lhd \mathbf{EA}\lhd \mathbf{PRA}\lhd\mathbf{PA}$.
\end{itemize}

\subsubsection{The limit of  ${\sf G1}$ w.r.t.~ interpretation and Turing reducibility}

We first discuss the limit of  ${\sf G1}$ for theories weaker than $\mathbf{PA}$ w.r.t.~ interpretation, i.e.~ finding  a  theory with minimal degree of interpretation for which  ${\sf G1}$ holds.

\smallskip

First of all, a natural question is: is $\mathbf{Q}$  the weakest finitely axiomatized essentially undecidable theory w.r.t.~ interpretation such that $\mathbf{R}\lhd\mathbf{Q}$? The following theorem tells us that the answer is no: for any finitely axiomatized subtheory $A$ of $\mathbf{Q}$ that extends $\mathbf{R}$, we
can find a finitely axiomatized subtheory $B$ of $A$ such that $B$ extends
$\mathbf{R}$ and $B$ does not interpret $A$.\smallskip

\begin{theorem}[Visser, Theorem 2, \cite{paper on Q}]~\label{Visser on Q}
Suppose $A$ is a finitely axiomatized consistent theory and $\mathbf{R} \subseteq A$. Then there is a finitely axiomatized theory $B$ such
that $\mathbf{R} \subseteq B \subseteq A$ and $B\lhd A$.
\end{theorem}

Define $X=\{S: \mathbf{R}\unlhd S\lhd \mathbf{Q}$ and $S$ is finitely axiomatized\}.
Theorem \ref{Visser on Q}  shows that the structure $\langle X, \lhd\rangle$ is not well-founded.\smallskip

\begin{theorem}[Visser, Theorem 12, \cite{paper on Q}]~\label{incomparable thm}
Suppose $A$ and $B$ are finitely axiomatized theories that weakly interpret   $\mathbf{Q}$. Then there are finitely axiomatized theories $\overline{A}\supseteq A$ and $\overline{B}\supseteq B$ such that $\overline{A}$ and $\overline{B}$ are incomparable (i.e. $\overline{A}\ntrianglelefteq \overline{B}$ and $\overline{B}\ntrianglelefteq \overline{A}$).
\end{theorem}
Theorem \ref{incomparable thm} shows that there are incomparable theories  extending $\mathbf{Q}$ w.r.t.~ interpretation.

\smallskip

Up to now, we do not have an example of essentially undecidable theory that is weaker than $\mathbf{R}$ w.r.t.~ interpretation.
To this end, we introduce Je\v{r}\'{a}bek's theory $\mathbf{Rep_{\sf PRF}}$.\smallskip
\begin{definition}[The system $\mathbf{Rep_{\sf PRF}}$]~
\begin{itemize}
  \item Let ${\sf PRF}$ denote the sets of all partial recursive functions.
  \item The language $L(\mathbf{Rep_{\sf PRF}})$ consists of constant symbols $\overline{n}$ for each $n \in \omega$, and function symbols $\overline{f}$ of appropriate arity for each partial recursive function $f$.
  \item The theory $\mathbf{Rep_{\sf PRF}}$ has axioms:
\begin{itemize}
  \item $\overline{n} \neq \overline{m}$ for $n \neq m \in \omega$;
  \item $\overline{f}(\overline{n_0}, \cdots, \overline{n_{k-1}}) = \overline{m}$
for each $k$-ary partial recursive function $f$ such that $f(n_0, \cdots, n_{k-1}) = m$ where $n_0, \cdots, n_{k-1}, m\in \omega$.
\end{itemize}
\end{itemize}
\end{definition}
The theory $\mathbf{Rep_{\sf PRF}}$ is essentially undecidable since all recursive functions are representable in it. Since $\mathbf{Rep_{\sf PRF}}$ is locally finitely satisfiable, by Theorem \ref{visser thm on R}, $\mathbf{Rep_{\sf PRF}}\unlhd \mathbf{R}$. Je\v{r}\'{a}bek \cite{Emil} proves that $\mathbf{R}$ is not interpretable in
$\mathbf{Rep_{\sf PRF}}$. Thus $\mathbf{Rep_{\sf PRF}}\lhd \mathbf{R}$.

\smallskip

Cheng \cite{Cheng 19} provides more examples of a theory $S$ such that $\sf G1$ holds for $S$ and $S\lhd \mathbf{R}$, and shows that we can find  many theories $T$ such that $\sf G1$ holds for $T$ and $T \lhd\mathbf{R}$  based on Je\v{r}\'{a}bek's work \cite{Emil} which uses model theory.\smallskip

\begin{theorem}[Cheng, \cite{Cheng 19}]~\label{main thm}
For any recursively inseparable pair $\langle A,B\rangle$, there is a r.e.~ theory $U_{\langle A,B\rangle}$ such that $\sf G1$  holds for $U_{\langle A,B\rangle}$,  and $U_{\langle A,B\rangle}\lhd\mathbf{R}$.
\end{theorem}

Define ${\sf D}=\{S: S\lhd \mathbf{R}$ and $\sf G1$  holds for theory $S$\}. Theorem \ref{main thm} shows that we could find many witnesses for ${\sf D}$.  Naturally, we could ask the following questions:\smallskip
\begin{question}~\label{open qn}
\begin{itemize}
  \item   Is $\langle {\sf D}, \lhd\rangle$ well-founded?
  \item Are any two elements of $\langle {\sf D}, \lhd\rangle$ comparable?
  \item  Does there exist a minimal theory w.r.t.~ interpretation  such that $\sf G1$ holds for it?
\end{itemize}
\end{question}
We conjectured the following answers to these questions: $\langle {\sf D}, \lhd\rangle$ is not well founded, $\langle {\sf D}, \lhd\rangle$ has incomparable elements, and there is no minimal theory w.r.t.~ interpretation  for which $\sf G1$ holds.

\smallskip

Finally, we discuss the limit of applicability of ${\sf G1}$ w.r.t.~ Turing reducibility. We  have discussed  the limit of applicability of ${\sf G1}$ w.r.t.~ interpretation.  A natural question is: what is the limit of applicability of ${\sf G1}$ w.r.t.~ Turing reducibility.
\smallskip

\begin{definition}[Turing reducibility, the structure ${\sf \overline{D}}$]~\label{}
\begin{itemize}
  \item Let $\mathcal{R}$ be the structure of the r.e.~ degrees with the ordering $\leq_{T}$ induced by Turing
reducibility  with the least element $\mathbf{0}$ and the greatest element $\mathbf{0}^{\prime}$.
  \item
Let ${\sf \overline{D}}=\{S: S<_{T} \mathbf{R}$ and $\sf G1$  holds for theory $S$\} where $S<_{T} \mathbf{R}$ stands for $S\leq_{T} \mathbf{R}$ but $\mathbf{R} \nleq_{T} S$.
\end{itemize}
\end{definition}

\smallskip

Cheng \cite{Cheng 19} shows that for any Turing degree $\mathbf{0}< \mathbf{d}<\mathbf{0}^{\prime}$, there is a theory $U$ such that $\sf G1$  holds for $U$, $U<_{T} \mathbf{R}$, and $U$ has Turing degree $\mathbf{d}$. As a corollary of this result and known results about the degree structure of $\langle \mathcal{R}, <_{T}\rangle$ in recursion theory, we can answer above questions  for the structure $\langle {\sf \overline{D}}, <_{T}\rangle$:  \smallskip
\begin{theorem}[Cheng, \cite{Cheng 19}]~
\begin{itemize}
  \item $\langle {\sf \overline{D}}, <_{T}\rangle$ is not well-founded;
  \item $\langle {\sf \overline{D}}, <_{T}\rangle$ has incomparable elements;
  \item There is no minimal theory w.r.t.~ Turing reducibility  such that $\sf G1$ holds for it.
\end{itemize}
\end{theorem}

Moreover, Cheng \cite{Cheng 19} shows that for any Turing degree $\mathbf{0}< \mathbf{d}<\mathbf{0}^{\prime}$, there is a theory $U$ such that $\sf G1$  holds for $U$, $U\unlhd \mathbf{R}$, and $U$ has Turing degree $\mathbf{d}$.
Thus, examining the limit of applicability of ${\sf G1}$ w.r.t.~ interpretation is much harder than that w.r.t.~ Turing reducibility. The structure of $\langle {\sf D}, \lhd\rangle$ is a deep and interesting open question for future research.

\section{The limit of the applicability of  ${\sf G2}$}\label{The limit of applicability}

\subsection{Introduction}

In our view, ${\sf G2}$ is \emph{fundamentally different} from ${\sf G1}$.
In fact, both mathematically and philosophically, $\sf G2$ is more problematic than $\sf G1$ for the following reason.  On one hand, in the case of $\sf G1$, we can construct a natural independent sentence with real  mathematical content \emph{without} referring to arithmetization and provability predicates.
On the other hand, the meaning of $\sf G2$ strongly depends on how we exactly formulate the consistency statement.

\smallskip

Similar to \cite{Halbach 2014a}, we call a result \emph{intensional} if it depends on (the details of) the representation used.  Thus, $\sf G1$ can be called extensional (that is, non-intensional), while $\sf G2$ is (highly) intensional.
We refer to Section \ref{intension problem of G2} for more discussion on the intensionality of $\sf G2$.
%However, ${\sf G2}$ is intensional, and ``whether ${\sf G2}$ holds" depends on varied factors.
In this section, we  discuss the limit of applicability of  ${\sf G2}$: under what conditions ${\sf G2}$ holds, and under what conditions ${\sf G2}$ fails. In Section \ref{Generalizations of G2}, we  discuss generalizations of ${\sf G2}$.
%In Section \ref{intension problem of G2}, we discuss the intensionality  of $\sf G2$.

\subsection{Some generalizations of ${\sf G2}$}\label{Generalizations of G2}
After G\"{o}del, generalizations of ${\sf G2}$ are the subject of extensive studies. We know that $\sf G2$ holds for any consistent r.e.~ extension of $\mathbf{PA}$. However, it is not  true that $\sf G2$ holds for any extension of $\mathbf{PA}$. For example, Karl-Georg Niebergall \cite{Niebergall}  shows that the theory ($\mathbf{PA}+\mathbf{RFN}(\mathbf{PA}))\cap (\mathbf{PA}$ + all true $\Pi^0_1$ sentences) can prove its own canonical  consistency sentence.\footnote{For the definition of $\mathbf{RFN}(\mathbf{PA})$, we refer to \cite{Aspects of incompleteness}: $\mathbf{RFN}(\mathbf{PA})=\{\forall x((\Gamma(x)\wedge \mathbf{Pr}_{\mathbf{PA}}(x))\rightarrow \mathbf{Tr}_{\Gamma}(x)): \Gamma$ arbitrary\}.}

\smallskip

Similarly to ${\sf G1}$,  one can generalise ${\sf G2}$ to  arithmetically definable non-r.e.~ extensions of $\mathbf{PA}$.
Kikuchi and Kurahashi \cite{Kikuchi 17} reformulate ${\sf G2}$ as: if $S$ is a $\Sigma^0_{1}$-definable and consistent extension of $\mathbf{PA}$, then for any $\Sigma^0_{1}$ definition $\sigma(u)$ of $S$, $S \nvdash \mathbf{Con}_{\sigma}(S)$ (see fact 5.1 \cite{Kikuchi 17}).
Kikuchi and Kurahashi \cite{Kikuchi 17} generalize ${\sf G2}$ to  arithmetically definable non-r.e.~ extensions of $\mathbf{PA}$ and prove that if $S$ is a $\Sigma^0_{n+1}$-definable and $\Sigma^0_{n}$-sound extension of $\mathbf{PA}$, then there exists a $\Sigma^0_{n+1}$ definition $\sigma(u)$ of some axiomatization of $Th(S)$ such that $\mathbf{Con}_{\sigma}(S)$ is independent of $S$. This corollary shows that the witness for the generalized version of ${\sf G1}$ can be
provided by the appropriate consistency statement.

\smallskip

Chao-Seraji \cite{Chao-Seraji 18}  and Kikuchi-Kurahashi \cite{Kikuchi 17} give another generalization of ${\sf G2}$ to  arithmetically definable non-r.e.~ extensions of $\mathbf{PA}$: for each $n \in \omega$, any $\Sigma^0_{n+1}$-definable and $\Sigma^0_{n}$-sound extension of $\mathbf{PA}$ cannot prove its own
$\Sigma^0_{n}$-soundness (see \cite[Theorem 2]{Chao-Seraji 18} and \cite[Theorem 5.6]{Kikuchi 17}).
The optimality of this generalization is shown in \cite{Chao-Seraji 18}:  there is a $\Sigma^0_{n+1}$-definable and $\Sigma^0_{n-1}$-sound extension of $\mathbf{PA}$ that proves its own $\Sigma^0_{n-1}$-soundness for $n>0$ (see \cite[Theorem 3]{Chao-Seraji 18}).

\smallskip

Let $T$ be a consistent r.e.~ extension of $\mathbf{Q}$. Kreisel \cite{Kreisel 1962} shows that $\neg \mathbf{Con}(T)$ is $\Pi^0_1$-conservative over $T$ which is a generalization of ${\sf G2}$.
We can also generalize $\sf G2$ via the notion of standard provability predicate.\smallskip

\begin{theorem}\label{standard provability predicate}
Let $T$ be any  consistent r.e.~ extension of $\mathbf{Q}$. If $\mathbf{Pr}_T(x)$ is a standard provability predicate, then $T\nvdash \mathbf{Con}(T)$.
\end{theorem}

Lev Beklemishev and Daniyar Shamkanov \cite{Beklemishev 16}
prove that in an abstract setting that presupposes the presence of G\"{o}del's fixed
point (instead of directly constructing it, as in the case of formal arithmetic), the Hilbert-Bernays-L\"{o}b conditions implies ${\sf G2}$ even with fairly minimal conditions on the underlying logic.
The following two theorems,  due to Feferman and Visser, generalize $\sf G2$ in terms of the notion of interpretation.\smallskip

\begin{theorem}[Feferman's theorem on the interpretability of inconsistency, \cite{Feferman 60}]~
If $T$ is a consistent r.e.~ extension of
$\mathbf{Q}$, then $T + \neg \mathbf{Con}(T)$ is interpretable in $T$.
\end{theorem}
\smallskip

\begin{theorem}[Pudl\'{a}k, \cite{Pudlak 85, Pudlak 93}]~\label{general G2}
There is no  consistent r.e.~ theory $S$ such that ($\mathbf{Q} + \mathbf{Con}(S)) \unlhd S$.
\footnote{Instead of
Robinson's Arithmetic $\mathbf{Q}$, we can as well have taken $\mathbf{S^1_2}$, or $\mathbf{PA}^-$, or $I\Delta_0 + \Omega_1$.
Moreover, instead of an arithmetical theory we can have employed a string theory like Grzegorczyk's
theory $\mathbf{TC}$ or adjunctive set theory $\mathbf{AS}$. All these theories are the same
in the sense that they are mutually interpretable (see \cite{Visser 16}).}
\end{theorem}
As a corollary of Theorem \ref{general G2}, for any  consistent r.e.~ theory $S$ that interprets $\mathbf{Q}$, $\sf G2$ holds for $S$: $S\nvdash \mathbf{Con}(S)$.
The  Arithmetic Completeness Theorem tells us that $S\unlhd (\mathbf{Q} +
\mathbf{Con}(S))$ (see \cite{Visser 11} for the details). As a corollary, we have the following version of $\sf G2$ which highlights the interpretability power of consistency statements.

\begin{corollary}
For any consistent r.e.~ theory $S$, we have $S\lhd (\mathbf{Q} + \mathbf{Con}(S))$.
\end{corollary}

\begin{definition}
Let $T$ be a consistent extension of $\mathbf{Q}$. A formula $I(x)$ with one free variable (understood
as a number variable) is a definable cut in $T$ (in short, a $T$-cut) if
\begin{itemize}
  \item $T\vdash I(\mathbf{0})$;
  \item $T\vdash \forall x (I(x) \rightarrow I(x+1))$;
  \item $T\vdash \forall x\forall y (y < x \wedge I(x) \rightarrow I(y))$.
\end{itemize}
\end{definition}

\begin{definition}
Let $T\supseteq I\Sigma_1$, let $J$ be a $T$-cut and let $\tau$ be a $\Sigma_0^{\mathbf{exp}}$-definition
of $T$.\footnote{We extend the language $L(\mathbf{PA})$ by a new unary function
symbol $\overline{2}^x$ for the $x$-th power of two. The extended language is denoted
$L_0(\mathbf{exp})$. A formula is $\Sigma_0^{\mathbf{exp}}$ if it results from atomic formulas of $L_0(\mathbf{exp})$ by iterated
application of logical connectives and bounded quantifiers of the form $(\forall x \leq
y)$ or $(\exists x \leq y)$ (see \cite{Pudlak 93}).}
\begin{itemize}
  \item $\mathbf{Pr}_{\tau}^I(x)$ is the formula $\exists y(I(y)\wedge \mathbf{Proof}_{\tau}^I(x,y))$
(saying that there is a $\tau$-proof of $x$ in $I$).
  \item $\mathbf{Con}_{\tau}^I$ is the formula $\neg\exists y(I(y)\wedge \mathbf{Proof}_{\tau}^I(\mathbf{0}\neq \mathbf{0},y))$.
\end{itemize}
\end{definition}
The following theorem generalizes $\sf G2$ to definable cuts.\smallskip

\begin{theorem}[Theorem 3.11, \cite{Pudlak 93}]~
Let $T\supseteq I\Sigma_1$, let $J$ be a $T$-cut and $\tau$ a $\Sigma_0^{\mathbf{exp}}$-definition of $T$. Then $T\nvdash  \mathbf{Con}_{\tau}^I$.
\end{theorem}
Next, consider a theory $U$ and an interpretation $N$ of the Tarski-Mostowski-Robinson theory $\mathbf{R}$ in $U$. A $U$-predicate $\triangle$ is an $L$-predicate for $U, N$ if it satisfies
the following L\"{o}b conditions.  We write $\triangle A$ for $\triangle (\ulcorner A \urcorner)$, where $\ulcorner A \urcorner$
is the numeral of
the G\"{o}del number of $A$ and we interpret the numbers via $N$. The G\"{o}del numbering is
supposed to be fixed and standard.\smallskip
\begin{definition}[L\"ob conditions]~
\begin{description}
  \item[L1] $\vdash A \Rightarrow  \ \vdash \triangle A$.
  \item[L2] $\triangle A, \triangle (A \rightarrow B) \vdash \triangle B$.
  \item[L3] $\triangle A \vdash \triangle \triangle A$.
\end{description}
\end{definition}
\smallskip
\begin{proposition}[L\"{o}b's theorem, Theorem 3.3.2, \cite{Visser 16}]~\label{lob thm}
Suppose that $U$ is a theory, $N$ is an interpretation  of the theory $\mathbf{R}$ in $U$, and $\triangle$ is a  $U$-predicate that is an $L$-predicate for $U, N$. Then:
\begin{itemize}
  \item For all $U$-sentences $A$ we have: if $U \vdash \triangle A \rightarrow A$, then $U \vdash A$.
  \item For all $U$-sentences $A$ we have: $U \vdash \triangle (\triangle A \rightarrow A) \rightarrow \triangle A$.
\end{itemize}
\end{proposition}

As a corollary of Proposition \ref{lob thm}, we formulate a general version of $\sf G2$  which does not mention the notion of provability predicate.\smallskip

\begin{theorem}[Visser, \cite{Visser 16}]~
For all consistent theories $U$ and all interpretations $N$ of $\mathbf{R}$ in $U$ and all $L$-predicates $\triangle$ for $U, N$, we have $U \nvdash \neg \triangle \perp$.
\end{theorem}

\subsection{The intensionality  of $\sf G2$}\label{intension problem of G2}

%In this section, I always assume that $T$ is a recursively axiomatized consistent extension of $\mathbf{PA}$.

%I have examined generalizations of  $\sf G1$ to  weak arithmetic  weaker than $\mathbf{PA}$. Now I examine generalizations of  $\sf G2$ to  weak arithmetic  weaker than $\mathbf{PA}$.

In this section, we discuss the intensionality  of $\sf G2$ which reveals the limit of the applicability of $\sf G2$.

\subsubsection{Introduction}

For a consistent theory $T$, we say that $\sf G2$ holds for $T$ if the consistency of $T$ is not provable in $T$. However, this definition is vague, and whether $\sf G2$ holds for $T$ depends on how we formulate the consistency statement. We refer to this phenomenon as the intensionality of $\sf G2$.
In fact, $\sf G2$ is essentially different from $\sf G1$ due to the intensionality of $\sf G2$: ``whether $\sf G2$ holds for the base theory"   depends on how we formulate the consistency statement in the first place.

\smallskip

Both mathematically and philosophically, $\sf G2$ is  more problematic than $\sf G1$. In the case of $\sf G1$, we are mainly interested in the fact that \emph{some} sentence is independent of $\mathbf{PA}$. We make no claim to the effect that that sentence ``really" expresses what we would express by saying ``$\mathbf{PA}$ cannot prove this sentence". But in the case of $\sf G2$, we are also interested in the content of the consistency statement.
We can say that ${\sf G1}$ is extensional in the sense that we can construct a concrete independent mathematical statement without referring to  arithmetization and provability predicate. However, ${\sf G2}$ is intensional and ``whether ${\sf G2}$ holds for $T$" depends on varied factors as we will discuss.

\smallskip

In this section, unless stated otherwise, we assume the following:
\begin{itemize}
  \item $T$ is a consistent r.e.~ extension of  $\mathbf{Q}$;
  \item the canonical arithmetic formula to express the consistency of the base theory $T$ is $\mathbf{Con}(T)\triangleq \neg \mathbf{Pr}_T(\mathbf{0}\neq \mathbf{0})$;
  %\item the base proof system is Hilbert-style system with cut elimination;
  \item the canonical numbering we use is G\"{o}del's numbering;
  \item the provability predicate we use is standard;
  \item the formula representing the set of axioms is $\Sigma^0_1$.
\end{itemize}

The intensionality of G\"{o}del sentence and the consistency statement has been widely discussed from the literature (e.g. Halbach-Visser \cite{Halbach 2014a, Halbach 2014b}, Visser \cite{Visser 11}).
Halbach and Visser examine the sources of intensionality in the construction of self-referential sentences of arithmetic  in \cite{Halbach 2014a, Halbach 2014b}, and argue that corresponding to the three stages of the construction of self-referential sentences of arithmetic, there are at least three sources
of intensionality: coding, expressing a property and  self-reference.  The three sources of intensionality are not independent of each other, and
a choice made at an earlier stage will have influences on the availability of choices
at a later stage.
Visser \cite{Visser 11} locates three
sources of indeterminacy in the formalization of a consistency statement for a
theory $T$:
\begin{itemize}
  \item the choice of a proof system;
  \item the choice of a way of numbering;
  \item the choice of a specific formula representing the set of axioms of $T$.
\end{itemize}

In summary, the intensional nature ultimately traces back to the various parameter choices that one has to make in arithmetizing the provability predicate. That is the source of both the intensional nature of the G\"{o}del sentence and the consistency sentence.

\smallskip

Based on this and other works from the literature, we argue that ``whether $\sf G2$ holds for the base theory" depends on the following factors:
\begin{enumerate}[(1)]
\item the choice of the provability predicate (Section \ref{poep1});
    \item the choice of the formula expressing consistency (Section \ref{poep2});
      \item the choice of the base theory (Section \ref{poep3});
      \item the choice of the numbering (Section \ref{poep4});
      \item the choice of the formula representing the set of axioms (Section \ref{poep5}).
\end{enumerate}
These factors are not independent, and
a choice made at an earlier stage may have effects on the choices
available at a later stage. In the following, unless stated otherwise, when we discuss how  $\sf G2$ depends on one factor,  we always assume that other factors are fixed,  and only the factor we are discussing is varied.
For example,  Visser \cite{Visser 11}  rests on fixed choices for (1)-(2) and (4)-(5) but varies the choice of (3); Grabmayr \cite{Grabmayr 18} rests on fixed choices for (1)-(3) and (5) but varies the choice of (4); Feferman \cite{Feferman 60} rests on fixed choices for (1)-(4)   but varies the choice of (5).
\subsubsection{The choice of provability predicate}\label{poep1}
In this section, we show that ``whether $\sf G2$ holds for the base theory" depends on the choice of the  provability predicate we use.

\smallskip

As Visser argues in \cite{Visser 16}, being a consistency statement is not an absolute concept but a role w.r.t.~ a choice of  provability predicate (see Visser \cite{Visser 16}). From Theorem \ref{standard provability predicate}, $\sf G2$ holds for standard provability predicates. However, $\sf G2$ may fail for non-standard provability predicates.

\smallskip

Mostowski \cite{Mostowski 65} gives an example of a $\Sigma^0_1$ provability predicate for which $\sf G2$ fails. Let $\mathbf{Pr}_T^M(x)$ be the $\Sigma^0_1$ formula ``$\exists y(\mathbf{Prf}_T(x, y) \wedge \neg \mathbf{Prf}_T(\ulcorner \mathbf{0}\neq \mathbf{0}\urcorner, y))$" where $\mathbf{Prf}_T(x, y)$ is a $\Delta^0_1$ formula saying that ``$y$ is a proof of $x$". Then $\neg \mathbf{Pr}_T^M(\ulcorner 0\neq 0\urcorner)$ is trivially provable
in $\mathbf{PA}$. We know that $\sf G2$ holds for provability predicates satisfying $\mathbf{D1}$-$\mathbf{D3}$. Since the formula $\mathbf{Pr}_T^M(x)$ satisfies $\mathbf{D1}$ and $\mathbf{D3}$, it does not satisfy $\mathbf{D2}$.
 Mostowski's example \cite{Mostowski 65} shows that $\sf G2$ may fail for
$\Sigma^0_1$ provability predicates satisfying $\mathbf{D1}$ and $\mathbf{D3}$.

\smallskip

One important non-standard provability predicate is
Rosser provability predicate $\mathbf{Pr}_T^R(x)$ introduced by
Rosser \cite{Rosser 36} to improve G\"{o}del's first incompleteness theorem.
Recall that we have defined the Rosser provability predicate in Definition \ref{Rosser predicate}. The consistency statement $\mathbf{Con}^R(T)$ based on a Rosser provability predicate $\mathbf{Pr}_T^R(x)$ is naturally defined as $\neg \mathbf{Pr}^R_{T}(\ulcorner \mathbf{0}\neq \mathbf{0}\urcorner)$.

\smallskip

It is an easy fact that for any sentence $\phi$ and Rosser provability predicate $\mathbf{Pr}^R_{T}(x)$, if $T\vdash\neg\phi$, then $T\vdash\neg\mathbf{Pr}^R_{T}(\ulcorner\phi\urcorner)$ (see \cite[Proposition 2.1]{Kurahashi 2017}).
As a corollary, the consistency statement $\mathbf{Con}^R(T)$ based on Rosser provability predicate $\mathbf{Pr}^R_{T}(x)$ is provable in $T$. In this sense, we can say that $\sf G2$  fails for the consistency statement constructed from Rosser provability predicates.

\smallskip

We can construct different Rosser provability predicates  with varied properties.
We know that each Rosser provability predicate $\mathbf{Pr}_{T}^R(x)$ does not satisfy at least one of conditions $\mathbf{D2}$ and $\mathbf{D3}$.  Guaspari and Solovay  \cite{Guaspari 79} establish a very powerful method of
constructing a new proof predicate with required properties from a given proof
predicate by reordering nonstandard proofs. Applying this tool, Guaspari and Solovay  \cite{Guaspari 79} construct   a Rosser provability predicate for which both $\mathbf{D2}$ and $\mathbf{D3}$ fail.
Arai  \cite{Arai 90} constructs   a Rosser provability predicate with condition $\mathbf{D2}$, and a Rosser provability predicate with condition $\mathbf{D3}$.

\smallskip

Slow provability, introduced by S.D. Friedman, M. Rathjen and A. Weiermann \cite{Friedman-Rathjen-Weiermann 13}, is another  notion of nonstandard provability for $\mathbf{PA}$ from the literature.  The slow consistency statement $\mathbf{Con}^{\ast}(\mathbf{PA})$ asserts that a contradiction
is not slow provable in $\mathbf{PA}$ (for the definition of $\mathbf{Con}^{\ast}(\mathbf{PA})$, we refer to \cite{Friedman-Rathjen-Weiermann 13}). In  fact, $\sf G2$ holds for  slow provability: Friedman, Rathjen and Weiermann show that  $\mathbf{PA} \nvdash \mathbf{Con}^{\ast}(\mathbf{PA})$ (see \cite[Proposition 3.3]{Friedman-Rathjen-Weiermann 13}).
Moreover,
Friedman, Rathjen and Weiermann \cite{Friedman-Rathjen-Weiermann 13} show that  $\mathbf{PA} +  \mathbf{Con}^{\ast}(\mathbf{PA}) \nvdash \mathbf{Con}(\mathbf{PA})$ (see \cite[Theorem 3.10]{Friedman-Rathjen-Weiermann 13}), and the logical strength of
the theory $\mathbf{PA}+\mathbf{Con}^{\ast}(\mathbf{PA})$ lies strictly between $\mathbf{PA}$ and $\mathbf{PA}+ \mathbf{Con}(\mathbf{PA})$: $\mathbf{PA} \varsubsetneq \mathbf{PA}+ \mathbf{Con}^{\ast}(\mathbf{PA}) \varsubsetneq \mathbf{PA}+ \mathbf{Con}(\mathbf{PA})$. Henk and Pakhomov \cite{Henk-Pakhomov 17} study three variants of
slow provability, and show that the associated consistency statement of each of these notions of provability
yields a theory that lies strictly between $\mathbf{PA}$ and $\mathbf{PA}+\mathbf{Con}(\mathbf{PA})$ in terms of
logical strength.

\subsubsection{The choice of the formula expressing consistency}\label{poep2}
We show that ``whether $\sf G2$ holds for the base theory" depends on the choice of the arithmetic formula used to express consistency.
In the literature, an arithmetic formula  is usually used to express the consistency statement.
Artemov \cite{Sergei Artemov} argues  that in Hilbert's consistency program, the original formulation of consistency ``no sequence of formulas  is a derivation of a contradiction" is about finite sequences
of formulas, not about arithmetization, proof codes, and internalized quantifiers.

\smallskip

The canonical consistency statement, the arithmetical formula  $\mathbf{Con}(\mathbf{PA})$, says that for all $x$, $x$ is not a code of
a proof of a contradiction in $\mathbf{PA}$. In a nonstandard model of $\mathbf{PA}$, the universal quantifier ``for all $x$" ranges over both standard and nonstandard numbers, and hence $\mathbf{Con}(\mathbf{PA})$ expresses the consistency of both standard and
nonstandard proof codes (see \cite{Sergei Artemov}). Thus, $\mathbf{Con}(\mathbf{PA})$ is stronger than the original formulation of consistency which only talks about sequences  of formulas and such sequences have only standard codes.
Hence, Artemov \cite{Sergei Artemov} concludes that $\sf G2$, saying that $\mathbf{PA}$ cannot prove $\mathbf{Con}(\mathbf{PA})$, does not actually exclude finitary consistency proofs of the original formulation of consistency (see \cite{Sergei Artemov}).

\smallskip

Artemov shows that the original formulation of consistency admits a direct proof in informal arithmetic, and this proof is formalizable  in $\mathbf{PA}$ (see \cite{Sergei Artemov}).\footnote{Informal arithmetic is the theory of informal elementary number theory containing recursive identities of addition and multiplication as well as the induction principle. The formal arithmetic $\mathbf{PA}$ is just the conventional formalization of the informal arithmetic (see \cite{Sergei Artemov}).}
Artemov's work establishes the consistency of $\mathbf{PA}$ by finitary means, and vindicates Hilbert's consistency program  to some extent.

\smallskip

In the following, we use a single arithmetic sentence to express the consistency statement.
Among consistency statements defined via arithmetization, there are three candidates of arithmetic formulas to express consistency as follows:
\begin{itemize}
  \item $\mathbf{Con}^0(T) \triangleq \forall x(\mathbf{Fml}(x) \wedge \mathbf{Pr}_T(x) \rightarrow \neg \mathbf{Pr}_T(\dot{\neg} x))$;\footnote{$\mathbf{Fml}(x)$ is the formula which represents the relation that $x$ is a code of a formula.}
      \item  $\mathbf{Con}(T) \triangleq \neg \mathbf{Pr}_T(\ulcorner \mathbf{0}\neq \mathbf{0}\urcorner)$;
  \item $\mathbf{Con}^1(T) \triangleq \exists x(\mathbf{Fml}(x) \wedge \neg \mathbf{Pr}_T(x))$.
\end{itemize}

G\"{o}del originally formulates $\sf G2$ with the
consistency statement $\mathbf{Con}^1(T)$: if $T$ is consistent,  then $T\nvdash \mathbf{Con}^1(T)$. From the literature, $\mathbf{Con}(T)$ is the widely used canonical consistency statement.  Note that $\mathbf{Con}^0(T)$ implies $\mathbf{Con}(T)$, and $\mathbf{Con}(T)$ implies $\mathbf{Con}^1(T)$. However the converse
implications do not hold in general (see \cite{note on derivability conditions}).
Kurahashi \cite{note on derivability conditions}  proposes different sets of derivability conditions (local version, uniform version and global version), and examines whether they are
sufficient to show the unprovability of these consistency statements (e.g.~ $\mathbf{Con}(T), \mathbf{Con}^0(T)$ and $\mathbf{Con}^1(T)$).

%Both of the Hilbert-Bernays derivability condition and the
%Hilbert-Bernays-L\"{o}b derivability condition are not sufficient conditions to prove G\"{o}del's original
%statement of $\sf G2$.

%In \cite{Hilbert-Bernays}, David Hilbert and Paul Bernays have developed a more general variant
%of the theorem that covered far weaker systems and was applicable to rather broad
%class of formalizations of the consistency assertion, rather than just the fixed G\"{o}del's
%formalization.

\begin{description}
  \item[HB1] If $T \vdash\phi \rightarrow\varphi$, then $T \vdash \mathbf{Pr}_T(\ulcorner\phi\urcorner) \rightarrow \mathbf{Pr}_T(\ulcorner\varphi\urcorner)$.
  \item[HB2] $T \vdash \mathbf{Pr}_T(\ulcorner\neg\phi(x)\urcorner) \rightarrow \mathbf{Pr}_T(\ulcorner\neg\phi(\dot{x})\urcorner)$.
  \item[HB3] $T \vdash f(x) = 0 \rightarrow \mathbf{Pr}_T(\ulcorner f(\dot{x}) = 0\urcorner)$ for every primitive recursive term
$f(x)$.
\end{description}
$\mathbf{HB1}$-$\mathbf{HB3}$ is called the Hilbert-Bernays derivability conditions.
If a provability predicate $\mathbf{Pr}_T(x)$ satisfies   $\mathbf{HB1}$-$\mathbf{HB3}$, then $T\nvdash\mathbf{Con}^0(T)$ (see \cite{Hilbert-Bernays}).
Kurahashi \cite{Rosser provability and G2} constructes two Rosser provability predicates satisfying $\mathbf{HB1}$-$\mathbf{HB3}$. Thus, $\mathbf{HB1}$-$\mathbf{HB3}$ is not  sufficient to prove that $T\nvdash \mathbf{Con}(T)$.

\smallskip

L\"{o}b \cite{Lob 55} proves that if $\mathbf{Pr}_T(x)$ satisfies  the Hilbert-Bernays-L\"{o}b derivability conditions $\mathbf{D1}$-$\mathbf{D3}$ (see Definition \ref{frak}), then L\"{o}b's theorem holds: for any sentence $\phi$, if $T \vdash \mathbf{Pr}_T(\ulcorner\phi\urcorner) \rightarrow\phi$,
then $T\vdash\phi$.
It is well-known that L\"{o}b's theorem implies $\sf G2$: $T\nvdash \mathbf{Con}(T)$  (see \cite{metamathematics}). Thus,
if a provability predicate $\mathbf{Pr}_T(x)$ satisfies $\mathbf{D1}$-$\mathbf{D3}$, then $T\nvdash \mathbf{Con}(T)$. Kurahashi \cite[Proposition 4.11]{note on derivability conditions} constructes a provability predicate $\mathbf{Pr}_T(x)$
with conditions $\mathbf{D1}$-$\mathbf{D3}$, but $T\vdash \mathbf{Con}^1(T)$. Thus,  $\mathbf{D1}$-$\mathbf{D3}$ is not  sufficient to prove that $T\nvdash \mathbf{Con}^1(T)$.

\smallskip

%Robert Jeroslow \cite{Jeroslow 73} proved that even under the assumption of certain weakening
%of the Hilbert-Bernays-L\"{o}b condition, ${\sf G2}$ holds for a large variety of deductive systems that have some means of
%representation of their own syntax.
Montagna \cite{Montagna 79} proves that if a provability predicate $\mathbf{Pr}_T(x)$ satisfies the following two conditions, then  $T\nvdash \mathbf{Con}^1(T)$:
\begin{itemize}
  \item $T \vdash \forall x$(``x is a logical axiom" $\rightarrow \mathbf{Pr}_T(x)$);
  \item $T \vdash \forall x\forall y(\mathbf{Fml}(x) \wedge \mathbf{Fml}(y) \rightarrow (\mathbf{Pr}_T(x\rightarrow y) \rightarrow  (\mathbf{Pr}_T(x) \rightarrow \mathbf{Pr}_T(y))))$.
\end{itemize}

\subsubsection{The choice of base theory}\label{poep3}
We show that ``whether $\sf G2$ holds for the base theory" depends on the base theory we choose.
A foundational question about $\sf G2$ is: how much  information about arithmetic is required for the proof of $\sf G2$. If the base system does not contain enough information about arithmetic,  then $\sf G2$ may fail.
The widely used notion of consistency  is consistency in proof systems with cut elimination. However, notions like cutfree consistency, Herbrand consistency, tableaux consistency, and restricted consistency  for different base theories behave
differently (see \cite{Visser 11}). We do have proof systems that prove their own cutfree consistency: for example, finitely axiomatized sequential theories prove their own cut-free consistency on a definable cut (see \cite{Visser 19}, p.25).

\smallskip

A natural question is: whether $\sf G2$ can be generalized to base  systems  weaker than $\mathbf{PA}$ w.r.t.~ interpretation. As a corollary  of Theorem \ref{general G2}, we have $\mathbf{Q}\nvdash \mathbf{Con(Q)}$ and hence $\mathbf{Q}\nvdash \mathbf{Con}^0(\mathbf{Q})$.  Bezboruah and Shepherdson \cite{Bezboruah 76} define the consistency of $\mathbf{Q}$ as the sentence $\mathbf{Con}^0(\mathbf{Q})$\footnote{This sentence says that for any $x$, if $x$ is the code of a formula $\phi$ and $\phi$ is provable in $\mathbf{Q}$, then $\neg\phi$ is not provable in $\mathbf{Q}$.}, and prove  that $\sf G2$ holds for $\mathbf{Q}: \mathbf{Q}\nvdash \mathbf{Con}^0(\mathbf{Q})$. However, the method used by Bezboruah and  Shepherdson in \cite{Bezboruah 76} is quite different from Theorem \ref{general G2}.  Bezboruah-Shepherdson's proof  depends
on some specific assumptions about the coding, does
not easily generalize to stronger theories, and tells us
nothing about the question whether $\mathbf{Q}$ can prove its consistency on some definable cut (see \cite{paper on Q}).
The next  question is: whether $\sf G2$  holds for other theories weaker than $\mathbf{Q}$ w.r.t.~ interpretation (e.g.~  $\mathbf{R}$). In a forthcoming paper, we will show that $\sf G2$  holds for $\mathbf{R}$ via the canonical consistency statement. However, we can find weak theories mutually interpretable with  $\mathbf{R}$ for which $\sf G2$  fails.

\smallskip

Willard \cite{Willard 06}  explores the generality and boundary-case
exceptions of $\sf G2$ over some base theories.
Willard constructs examples of r.e.~ arithmetical theories that cannot prove the totality of their successor functions but can prove
their own canonical consistencies  (see \cite{Willard 01}, \cite{Willard 06}).
However, the theories Willard constructs are not
completely natural since some  axioms are constructed using G\"{o}del's
Diagnolisation Lemma.
Pakhomov \cite{Fedor Pakhomov} constructs  a more
natural example of this kind.
Pakhomov \cite{Fedor Pakhomov} defines a theory $H_{<\omega}$, and shows that it proves its own canonical consistency.
Unlike Willard's theories, $H_{<\omega}$ isn't an arithmetical theory but a theory formulated in the
language of set theory with an additional unary function.
From \cite{Fedor Pakhomov}, $H_{<\omega}$  and $\mathbf{R}$ are mutually interpretable. Hence, the theory $H_{<\omega}$  can  be regarded as the set-theoretic analogue of
$\mathbf{R}$ from the interpretability theoretic point of view.
%Willard \cite{Willard 01} \cite{Willard 06} shows that $\sf G2$ does not hold for some axiom systems which fail to recognize successor as a total function and treat addition and multiplication as 3-way relations.

\smallskip

From Theorem \ref{general G2}, $\sf G2$ holds for any consistent r.e.~  theory interpreting $\mathbf{Q}$.
However, it is not true that $\sf G2$ holds for any consistent r.e.~ theory interpreting $\mathbf{R}$ since $H_{<\omega}$ interprets $\mathbf{R}$, but $\sf G2$ fails for $H_{<\omega}$.
We know that if $S\unlhd T$ and $\sf G1$ holds for $S$, then $\sf G1$ holds for $T$.
However, it is not true that if $S\unlhd T$ and $\sf G2$ holds for $S$, then $\sf G2$ holds for $T$ since $\mathbf{R}\unlhd H_{<\omega}$, $\sf G2$ holds for $\mathbf{R}$ but $\sf G2$ fails for $H_{<\omega}$. This shows the difference between $\mathbf{Q}$ and $\mathbf{R}$, and the difference between  $\sf G1$ and $\sf G2$.

\smallskip

%Visser  \cite{Visser 19} gives an example of a finitely axiomatized theory $B$ and a $\Sigma^0_1$-formula $\gamma$
%that numerates the axiom set in $N_0 : \mathbf{EA} \unlhd B$ but not in $N_1: \mathbf{EA} \unlhd B$ such that  in
%the case of $N_0$, $B$ cannot prove the canonical consistency of the theory axiomatized by
%$\gamma$, but in the case
%of $N_1$, $B$ proves the canonical consistency of the theory axiomatized by
%$\gamma$ (see Example 6.3 in \cite{Visser 19}).
One way to eliminate the intensionality of $\sf G2$ is to uniquely characterize the consistency statement.
In \cite{Visser 11}, Visser  proposes   the interesting question of  a coordinate-free formulation of $\sf G2$ and a unique characterization of  the consistency statement.
Visser \cite{Visser 11} shows that consistency for finitely axiomatized sequential theories can be uniquely characterized modulo $\mathbf{EA}$-provable equivalence (see \cite{Visser 11}, p.543). But characterizing the consistency of infinitely axiomatized r.e.~ theories is more delicate and a big open problem in the current research on the intensionality of $\sf G2$.
%Willard \cite{Willard 06} aims to explore the generality and boundary-case
%exceptions of $\sf G2$ under a class of axiom systems that fail to recognize successor as a total function and
%instead treat addition and multiplication as 3-way relations.
%Willard \cite{Willard 06} introduces the notion of ``Naming Convention" to refer to a particular scheme for assigning integer values to named
%constant symbols.
%It will turn out that our ability to either generalize the Second Incompleteness Theorem or to find
%boundary-case exceptions to it will depend on the choice of naming method.

\smallskip

After G\"{o}del, Gentzen constructs a theory $\mathbf{T}^{\ast}$ (primitive recursive arithmetic with the additional principle of  quantifier-free
transfinite induction up to the ordinal $\epsilon_0$)\footnote{$\epsilon_0$ is the first ordinal  $\alpha$  such that  $\omega^{\alpha }=\alpha$.}, and  proves  the consistency of $\mathbf{PA}$  in   $\mathbf{T}^{\ast}$.
Gentzen's theory $\mathbf{T}^{\ast}$ contains $\mathbf{Q}$ but does not contain $\mathbf{PA}$ since $\mathbf{T}^{\ast}$ does not prove the ordinary mathematical induction for all formulas.
 By the Arithmetized Completeness Theorem, $\mathbf{Q}+\mathbf{Con}(\mathbf{PA})$ interprets $\mathbf{PA}$.  Since Gentzen's theory $\mathbf{T}^{\ast}$ contains $\mathbf{Q}$ and $\mathbf{T}^{\ast}\vdash\mathbf{Con}(\mathbf{PA})$, Gentzen's theory $\mathbf{T}^{\ast}$ interprets $\mathbf{PA}$.
By Pudl\'{a}k's result that no consistent r.e.~ extension $T$ of $\mathbf{Q}$ can interpret $\mathbf{Q}+\mathbf{Con}(T)$,  $\mathbf{PA}$ does not interpret Gentzen's theory $\mathbf{T}^{\ast}$. Thus   $\mathbf{PA}\lhd  \mathbf{T}^{\ast}$.
Gentzen's work has opened a productive new direction in proof theory: finding the means necessary to
prove the consistency of a given theory. More powerful subsystems of Second-Order Arithmetic have been given consistency proofs by Gaisi Takeuti and others, and  theories that have been proved consistent by these methods are quite strong and include most ordinary mathematics.

\subsubsection{The choice of numbering}\label{poep4}
We show that ``whether $\sf G2$ holds for the base theory" depends on the choice of the numbering encoding the language.

\smallskip

For the influence of different numberings on $\sf G2$, we refer to \cite{Grabmayr 18}.
Any injective function $\gamma$
from a set of $L(\mathbf{PA})$-expressions to $\omega$ qualifies as a numbering. G\"{o}del's numbering is a special kind of numberings under which the G\"{o}del number of the set of axioms of $\mathbf{PA}$ is recursive. In fact, $\sf G2$ is sensitive to the way of numberings.
Let  $\gamma$ be a numbering and $\ulcorner\varphi^{\gamma}\urcorner$ denote $\overline{\gamma(\varphi)}$, i.e., the standard numeral of the $\gamma$-code of $\varphi$.\smallskip

\begin{definition}[Relativized L\"ob conditions]
A formula $\mathbf{Pr}_T^{\gamma}(x)$ is said to satisfy L\"{o}b's conditions relative to $\gamma$ for the base theory $T$ if  for all $L(\mathbf{PA})$-sentences   $\varphi$ and $\psi$ we have that:
\begin{description}
  \item[$\mathbf{D1^{\ast}}$] If $T \vdash\varphi$, then $\mathbf{PA} \vdash \mathbf{Pr}_T^{\gamma}(\ulcorner\varphi^{\gamma}\urcorner)$;
      \item[$\mathbf{D2^{\ast}}$] $T \vdash \mathbf{Pr}_T^{\gamma}(\ulcorner(\varphi
      \rightarrow\psi)^{\gamma}\urcorner)\rightarrow (\mathbf{Pr}_T^{\gamma}(\ulcorner\varphi^{\gamma}\urcorner)\rightarrow \mathbf{Pr}_T^{\gamma}(\ulcorner\psi^{\gamma}\urcorner))$;
  \item[$\mathbf{D3^{\ast}}$] $T \vdash \mathbf{Pr}_T^{\gamma}(\ulcorner\varphi^{\gamma}\urcorner)\rightarrow \mathbf{Pr}_T^{\gamma}(\ulcorner(\mathbf{Pr}_T^{\gamma}(\ulcorner\varphi^{\gamma}\urcorner))^{\gamma}\urcorner)$.
\end{description}
\end{definition}
Grabmayr \cite{Grabmayr 18} examines different criteria  of acceptability, and proves the invariance of $\sf G2$ with regard to acceptable numberings (for the definition of acceptable numberings, we refer to \cite{Grabmayr 18}).

\begin{theorem}[Invariance of $\sf G2$ under acceptable numberings, Theorem 4.8, \cite{Grabmayr 18}]~\label{acceptable numberings}
Let $\gamma$ be an acceptable numbering and $T$ be a consistent r.e.~ extension of  $\mathbf{Q}$. If   $\mathbf{Pr}_T^{\gamma}(x)$ satisfies L\"{o}b's conditions $\mathbf{D1^{\ast}}$-$\mathbf{D3^{\ast}}$ relative to $\gamma$ for $T$, then  $T \nvdash \neg \mathbf{Pr}_T^{\gamma}(\ulcorner(\mathbf{0}\neq \mathbf{0})^{\gamma}\urcorner)$.
\end{theorem}
Theorem \ref{acceptable numberings} shows that $\sf G2$ holds for acceptable numberings. But $\sf G2$ may fail for non-acceptable numberings. Grabmayr \cite{Grabmayr 18} gives some examples of deviant numberings $\gamma$ such that $\sf G2$ fails w.r.t.~ $\gamma$: $T\vdash\mathbf{Pr}_T^{\gamma}(\ulcorner(\mathbf{0}\neq \mathbf{0})^{\gamma}\urcorner)$.
\smallskip

\begin{definition}\label{def of numeration}
We say that $\alpha(x)$ is a numeration of $T$ if for any $n$, we have $\mathbf{PA} \vdash \alpha(\overline{n})$ if and only if $n$ is the G\"{o}del number of some $\phi\in T$.
\end{definition}

\subsubsection{The choice of the formula representing the set of axioms}\label{poep5}
We show that ``whether $\sf G2$ holds for $T$" depends on the way the axioms of $T$ are represented.

\smallskip

First of all, Definition \ref{def of provability predicate} gives a more general definition of provability predicate and consistency statement for $T$ w.r.t.~ the numeration of $T$.\smallskip
\begin{definition}\label{def of provability predicate}
Let $T$ be any  consistent r.e.~ extension of $\mathbf{Q}$ and $\alpha(x)$ be a formula in $L(T)$.
\begin{itemize}
  \item  Define the formula $\mathbf{Prf}_{\alpha}(x,y)$ saying ``$y$ is the G\"{o}del number of a proof of the formula with G\"{o}del number $x$
from the set of all sentences satisfying $\alpha(x)$".
\item  Define the provability predicate $\mathbf{Pr}_{\alpha}(x)$ of $\alpha(x)$ as $\exists y \mathbf{Prf}_{\alpha}(x,y)$ and the consistency statement
$\mathbf{Con}_{\alpha}(T)$ as $\neg \mathbf{Pr}_{\alpha}(\ulcorner \mathbf{0}\neq \mathbf{0}\urcorner)$.
\end{itemize}
\end{definition}
\noindent
%Let $T$ be a recursively axiomatized consistent extension of $\mathbf{PA}$.
For each formula $\alpha(x)$, we have:
\[\mathbf{D2^{\prime}}\quad \mathbf{PA} \vdash \mathbf{Pr}_{\alpha}(\ulcorner\varphi
      \rightarrow\psi\urcorner)\rightarrow (\mathbf{Pr}_{\alpha}(\ulcorner\varphi\urcorner)\rightarrow \mathbf{Pr}_{\alpha}(\ulcorner\psi\urcorner)).\]
If $\alpha(x)$ is a  numeration of $T$, then $\mathbf{Pr}_{\alpha}(x)$ satisfies
the following properties (see \cite[Fact 2.2]{Kurahashi 17}):
%The proof of $\mathbf{G_2}$ essentially uses the following drivability conditions:
\begin{description}
  \item[$\mathbf{D1^{\prime}}$] If $T \vdash\varphi$, then $\mathbf{PA} \vdash \mathbf{Pr}_{\alpha}(\ulcorner\varphi\urcorner)$;
  \item[$\mathbf{D3^{\prime}}$] If $\varphi$ is $\Sigma^0_1$, then $\mathbf{PA} \vdash \varphi\rightarrow \mathbf{Pr}_{\alpha}(\ulcorner\varphi\urcorner)$.
\end{description}

Now we give a new reformulation of $\sf G2$ via numerations.
\begin{theorem}\label{G2 under numeration}
Let $T$ be any  consistent r.e.~ extension of $\mathbf{Q}$. If $\alpha(x)$ is any $\Sigma^0_1$ numeration of $T$, then $T\nvdash \mathbf{Con}_{\alpha}(T)$.
\end{theorem}

In fact, $\sf G2$ holds for any $\Sigma^0_1$ numeration of $T$, but fails for some $\Pi^0_1$ numeration of $T$.
Feferman \cite{Feferman 60} constructs    a $\Pi^0_1$ numeration $\pi(x)$ of $T$ such
that $\sf G2$ fails, i.e. $\mathbf{Con}_{\pi}(T)\triangleq \neg\mathbf{Pr}_{\pi}(\ulcorner \mathbf{0}\neq \mathbf{0}\urcorner)$ is provable in $T$.
Feferman's construction keeps the proof predicate and its numbering fixed but varies the formula  representing the
set of axioms.
Notice that Feferman's
predicate satisfies $\mathbf{D1}$ and $\mathbf{D2}$, but does not satisfy $\mathbf{D3}$. Feferman's example  shows that  $\sf G2$ may fail for provability predicates satisfying $\mathbf{D1}$ and  $\mathbf{D2}$.

\smallskip

Generally, Feferman \cite{Feferman 60} shows   that if $T$ is a $\Sigma^0_{1}$-definable extension of $\mathbf{Q}$, then there is a $\Pi^0_{1}$ definition $\tau(u)$ of $T$ such that $T \vdash \mathbf{Con}_{\tau}(T)$.
In  summary, $\sf G2$ is not coordinate-free (it is dependent on numerations of $\mathbf{PA}$).
An important question is  how to formulate $\sf G2$ in a general way such that it is coordinate-free (independent of numerations of $T$).

\smallskip

The properties  of  the provability predicate are intensional and depend on the numeration of  the theory.  I.e.,~ under different numerations of $T$, the  provability predicate  may have different properties. It may happen that  $T$ has two numerations $\alpha(x)$ and $\beta(x)$ such that $\mathbf{Con}_{\alpha}(T)$ is not equivalent to $\mathbf{Con}_{\beta}(T)$. For example, under G\"{o}del's recursive numeration $\tau(x)$ and Feferman's $\Pi^0_1$ numeration $\pi(x)$ of $T$, the corresponding consistency statement $\mathbf{Con}_{\tau}(T)$ and  $\mathbf{Con}_{\pi}(T)$ are not equivalent. But $\PA$ does not know this fact, i.e. $\PA\nvdash \neg(\mathbf{Con}_{\tau}(T)\leftrightarrow \mathbf{Con}_{\pi}(T))$ since $\PA\nvdash \neg\mathbf{Con}_{\tau}(T)$.

\smallskip

Generally, Kikuchi and Kurahashi prove in \cite[Corollary 5.11]{Kikuchi 17} that if $T$ is $\Sigma^0_{n+1}$-definable and not $\Sigma^0_{n}$-sound, then there are $\Sigma^0_{n+1}$
definitions $\sigma_1(x)$ and $\sigma_2(x)$ of $T$ such that $T\vdash \mathbf{Con}_{\sigma_1}(T)$ and $T \vdash\neg \mathbf{Con}_{\sigma_2}(T)$.
\smallskip

%\subsubsection{Incompleteness and provability logic}\label{poep6}

Provability logic is an important tool for the study of incompleteness and meta-mathematics of arithmetic.
The origins of provability logic (e.g.~ Henkin's problem, the isolation of derivability conditions, L\"{o}b's theorem) are all closely tied to G\"{o}del's incompleteness theorems historically. In this sense, we can say that G\"{o}del's incompleteness theorems play a unifying role between first order arithmetic and provability logic.

\smallskip

Provability logic is the logic of properties of provability predicates. Note that $\sf G2$ is very sensitive to the properties of the provability predicate used in its formulation.  Provability logic provides us  with a new viewpoint and an important tool that can be used to understand incompleteness. Provability logic based on different provability predicates reveals the intensionality of provability predicates which is one source of the intensionality of $\sf G2$.

\smallskip

Let $T$ be a consistent r.e.~ extension of $\mathbf{Q}$, and  $\tau(u)$ be any numeration of $T$. Recall that an arithmetical interpretation $f$ is a mapping from the set of all modal propositional variables to the set of $L(T)$-sentences.
Every arithmetical interpretation $f$ is uniquely extended to the
mapping $f_{\tau}$ from the set of all modal formulas to the set of $L(T)$-sentences such
that $f_{\tau}$ satisfies the following conditions:
\begin{itemize}
  \item $f_{\tau}(p)$ is $f(p)$ for each propositional variable $p$;
  \item $f_{\tau}(\bot)$ is $\mathbf{0} \neq \mathbf{0}$;
  \item $f_{\tau}$ commutes with every propositional connective;
  \item $f_{\tau}(\Box A)$ is $\mathbf{Pr}_{\tau}(\ulcorner f_{\tau}(A)\urcorner)$ for every modal formula $A$.
\end{itemize}
Provability logic provides us with a new way of examining the intensionality of provability predicates. Under different numerations of $T$, the provability predicate may have different properties, and hence may correspond to different modal principles under different arithmetical interpretations.
\smallskip

\begin{definition}
Given a numeration $\tau(u)$ of $T$, the provability logic $\mathbf{PL}_{\tau}(T)$ of $\tau(u)$ is defined to be the set  of modal formulas $A$ such that $T \vdash f_{\tau}(A)$ for all arithmetical interpretations $f$.
\end{definition}
Note that the provability logic $\mathbf{PL}_{\tau}(T)$ of a $\Sigma^0_n$ numeration $\tau(x)$ of
$T$ is a normal modal logic.  A natural and interesting question is: which normal modal logic can be realized as  a provability logic $\mathbf{PL}_{\tau}(T)$ of some $\Sigma^0_n$
numeration $\tau(x)$ of $T$?
An interesting research program is to classify the provability logic $\mathbf{PL}_{\alpha}(T)$ according to the numeration $\alpha(x)$ of $T$. We first discuss $\Sigma^0_1$ numerations  of $T$.

\begin{theorem}[Generalized Solovay's Arithmetical Completeness Theorem, Theorem 2.5, \cite{Kurahashi 17}]~
Let $T$ be any consistent r.e.~ extension of $\mathbf{PA}$. If
$T$ is $\Sigma^0_1$-sound, then for any $\Sigma^0_1$ numeration $\alpha(x)$ of $T$, the provability logic $\mathbf{PL}_{\alpha}(T)$ is precisely $\mathbf{GL}$.\footnote{It is a big open problem that whether Solovay's arithmetical completeness theorem holds for weak arithmetic.}
\end{theorem}

Moreover, Visser \cite{Visser 84} examines all  possible provability logics for $\Sigma^0_1$ numerations of $\Sigma^0_1$-unsound
theories. To state Visser's result, we need some definitions.\smallskip

\begin{definition}[Definition 3.5-3.6, \cite{Kurahashi numerations 17}]~
 We define the sequence $\{\mathbf{Con}_{\tau}^n: n\in\omega\}$  recursively as follows: $\mathbf{Con}_{\tau}^0$ is $\mathbf{0} = \mathbf{0}$, and  $\mathbf{Con}_{\tau}^{n+1}$ is $\neg \mathbf{Pr}_{\tau}(\ulcorner \neg \mathbf{Con}_{\tau}^{n}\urcorner)$.
The height of $\tau(u)$ is the least natural number $n$ such that $T \vdash \neg \mathbf{Con}_{\tau}^{n}$ if such an $n$ exists. If not, the height of $\tau(u)$ is $\infty$.
\end{definition}

For $\Sigma^0_1$-unsound theories, Visser proves that $\mathbf{PL}_{\tau}(T)$ is determined by the
height of the numeration $\tau(u)$. Visser \cite[Theorem 3.7]{Visser 84} shows that the height of $\tau(u)$ is $\infty$ if and only if $\mathbf{PL}(\tau) = \mathbf{GL}$; and the height of $\tau(u)$ is $n$ if and only if $\mathbf{PL}(\tau) = \mathbf{GL} + \Box^n \bot$. Beklemishev \cite[Lemma 7]{Beklemishev 90} shows that if $T$ is $\Sigma^0_1$-unsound, then the
height of $\Sigma^0_1$ numerations of $T$ can take any values except $0$.

\smallskip

Let $U$ be any consistent theory of arithmetic. Based on the previous work by
Artemov, Visser and Japaridze, Beklemishev
\cite{Beklemishev 90}  proves that for $\Sigma^0_1$ numeration $\tau$ of $U$, $\mathbf{PL}_{\tau}(U)$ coincides with one of
the logics $\mathbf{GL}_{\alpha}, \mathbf{D}_{\beta}, \mathbf{S}_{\beta}$ and $\mathbf{GL}^-_{\beta}$ where $\alpha$ and $\beta$ are subsets of $\omega$ and $\beta$ is
cofinite (for definitions of $\mathbf{GL}_{\alpha}, \mathbf{D}_{\beta}, \mathbf{S}_{\beta}$ and $\mathbf{GL}^-_{\beta}$, we refer to \cite{Beklemishev 90, Artemov-Beklemishev 05}).

\smallskip

%For each $n \geq 1$, let $\mathbf{PL}_n(T)$ be the class $\{\mathbf{PL}_{\tau}(T) : \tau$ is a $\Sigma^0_n$ numeration of $T$\} of all possible provability logics for $\Sigma^0_n$ numerations
%of $T$ (see \cite[Definition 3.2]{Kurahashi numerations 17}).  As a corollary, if $T$ is $\Sigma^0_1$-unsound, then $\mathbf{PL}_1(T)$ is equal to the set $\{\mathbf{GL} + \Box^n \bot : n \geq 1\} \cup \{\mathbf{GL}\}$.
Feferman \cite{Feferman 60} constructs a $\Pi^0_1$ numeration $\pi(x)$ of $T$ such that the consistency statement $\mathbf{Con}_{\pi}(T)$ defined via $\mathbf{Pr}_{\pi}(x)$ is provable
in $T$. Thus, the provability logic $\mathbf{PL}_{\pi}(T)$ of $\mathbf{Pr}_{\pi}(x)$ contains the formula $\neg\Box\bot$, and is different from $\mathbf{GL}$.
However, the exact
axiomatization of the provability logic $\mathbf{PL}_{\pi}(T)$ under Feferman's numeration $\pi(x)$ is
not known.
Kurahashi \cite{Kurahashi 17} proves   that  for any recursively axiomatized consistent
extension $T$ of $\mathbf{PA}$, there exists a $\Sigma^0_2$  numeration $\alpha(x)$ of $T$ such that the
provability logic $\mathbf{PL}_{\alpha}(T)$ is the modal system $\mathbf{K}$. As a corollary, the
modal principles commonly contained in every provability logic $\mathbf{PL}_{\alpha}(T)$ of
$T$ is just $\mathbf{K}$.

\smallskip

It is often thought that  a provability predicate satisfies $\mathbf{D1}$-$\mathbf{D3}$ if and only if $\sf G2$ holds  (i.e.~ for the induced consistency statement $\mathbf{Con}(T)$ from the provability predicate, $T\nvdash \mathbf{Con}(T)$). But this is not true.
From Definition \ref{def of provability predicate}, conditions $\mathbf{D1}$-$\mathbf{D2}$ hold for any numeration of $T$. Whether the provability predicate satisfies condition $\mathbf{D3}$ depends on the numeration of $T$.
For any $\Sigma^0_1$-numeration $\alpha(x)$ of $T$,
$\mathbf{D3}$ holds for $\mathbf{Pr}_{\alpha}(x)$. From Kurahashi \cite{Kurahashi 17}, there is a $\Sigma^0_2$-numeration $\alpha(x)$ of  $T$ such that the provability logic for that numeration is precisely $\mathbf{K}$. Since $\mathbf{K}\nvdash \neg\square\bot$, as a corollary,
$\sf G2$ holds for $T$, i.e.~ $\mathbf{Con}_{\alpha}(T)$ defined as $\neg\mathbf{Pr}_{\alpha}(\ulcorner \mathbf{0}\neq \mathbf{0}\urcorner)$ is not provable in $T$.
But the L\"{o}b condition $\mathbf{D3}$ does not hold since $\mathbf{K}\nvdash \square A\rightarrow \square\square A$. This gives us an example of a $\Sigma^0_2$ numeration $\alpha(x)$ of $T$ such that $\mathbf{D3}$ does not hold for $\mathbf{Pr}_{\alpha}(x)$ but $\sf G2$ holds for $T$.
Thus, $\sf G2$ may hold for a provability predicate which does not satisfy the L\"{o}b condition $\mathbf{D3}$.

\smallskip

Moreover,
Kurahashi \cite{Kurahashi numerations 17} proves  that for each $n \geq 2$, there exists
a $\Sigma^0_2$  numeration $\tau(x)$ of $T$ such that the provability logic $\mathbf{PL}_{\tau}(T)$ is just the modal logic $\mathbf{K} + \square(\square^n p \rightarrow p) \rightarrow \square p$.
Hence there are infinitely many normal modal logics that are provability logics for some $\Sigma^0_2$
numeration of $T$.
A good question from Kurahashi \cite{Kurahashi numerations 17} is:  for $n \geq 2$, is the class of provability logics $\mathbf{PL}_{\tau}(T)$ for $\Sigma^0_n$
numerations $\tau(x)$ of $T$ the same as the class of provability logics $\mathbf{PL}_{\tau}(T)$ for $\Sigma^0_{n+1}$
numerations $\tau(x)$ of $T$? However,  this question is still open as far as we know.
Define that $\mathbf{KD}=\mathbf{K}+\neg\Box\bot$.    A natural and interesting question, which is also open as far as we know, is: can we find a numeration $\tau(x)$ of $T$ such that $\mathbf{PL}_{\tau}(T)=\mathbf{KD}$?

\smallskip

In  summary, $\sf G2$ is intensional with respect to the following parameters: the formalization of consistency, the base theory, the method of numbering, the choice of a provability predicate, and the representation of the set of axioms.   Current research on incompleteness reveals that  $\sf G2$ is a deep and profound  theorem both mathematically and philosophically in the foundations of mathematics, and there is a lot more to be explored about the intensionality of $\sf G2$.

\section{Conclusion}
We conclude this paper with some personal comments. To the author, the research on concrete incompleteness is very deep and important.

\smallskip

After G\"{o}del, people have found many different proofs of incompleteness theorems via pure logic, and many concrete independent statements with real mathematical contents. As Harvey Friedman comments, the research on concrete  mathematical incompleteness shows how the Incompleteness Phenomena touches
normal concrete mathematics, and reveals the impact
and significance of  the foundations of mathematics.

\smallskip

Harvey Friedman's research project on concrete   incompleteness plans to show  that we will be able to find, in just about
any subject of mathematics, many natural looking statements that are independent of $\mathbf{ZFC}$. Harvey Friedman's work is very profound and promising, and will reveal that incompleteness is everywhere in mathematics, which, if  it is true, may be one of the most important discoveries after G\"{o}del in the foundations of mathematics.

%\subsection{Incompleteness and provability logic}\label{provability logic}

%\subsection{Provability logic under numerations}

%First, we examine the characterization of truth and provability in provability logic.

%We consider the following conditions for Rosser provability predicates.

%\begin{definition}
%(\cite[Definition 3.1]{Kurahashi 2017}) For all formulas $\varphi$ and $\psi$:
%\begin{description}
  %\item[$D1^{R}$] If $T \vdash \varphi$, then $PA \vdash Pr_{T}^R(\ulcorner\varphi\urcorner)$.
  %\item[$D2^{R}$] $PA \vdash Pr_{T}^R(\ulcorner\varphi\rightarrow\psi\urcorner) \rightarrow (Pr_{T}^R(\ulcorner\varphi\urcorner) \rightarrow Pr_{T}^R(\ulcorner\psi\urcorner))$.
  %\item[$D3^{R}$] $PA \vdash Pr_{T}^R(\ulcorner\varphi\urcorner) \rightarrow \vdash Pr_{T}^R(\ulcorner Pr_{T}^R(\ulcorner\varphi\urcorner)\urcorner)$.
  %\item[
%\end{description}
%\end{definition}.

{100}

\end{document}